\documentclass[12pt]{amsart}
\usepackage{amssymb,amscd}
\usepackage{graphicx, color}

\newtheorem*{Whitney towers}{Theorem~\ref{Whitney towers}}
\newtheorem*{h-towers}{Theorems ~\ref{half} \& \ref{$(n)$-solvable}}

\newtheorem*{surgery curves}{Theorem~\ref{surgery curves}}
\newtheorem*{cg=0}{Theorem~\ref{vanish}}

\newtheorem{thm}{Theorem}[section]

\newtheorem{prop}[thm]{Proposition}
\newtheorem{cla}[thm]{Claim}
\theoremstyle{definition}
\newtheorem{defn}[thm]{Definition}

\newtheorem{prob}[thm]{Problem}
\newtheorem{note}[thm]{Note}

\numberwithin{equation}{section}
\numberwithin{figure}{section}

\newcommand{\x}{\times}
\newcommand{\np}{\newpage}
\newcommand{\Z}{\mathbb{Z}}
\newcommand{\N}{\mathbb{N}}

\newcommand{\Q}{\mathbb{Q}}

\newcommand{\R}{\mathbb{R}}

\def\yen{{\setbox0=\hbox{Y}Y\kern-.97\wd0\vbox{hrule height.lex width.98%
\wd0\kern.33ex\hrule height.lex width.98\wd0\kern.45ex}}}

\def\np{\newpage}

\begin{document}

\np

\title{Local moves on knots and products of knots}
\author{Louis H. Kauffman and  Eiji Ogasa}      

\thanks{\hskip-4mm E-mail: kauffman@uic.edu\quad  
ogasa@mail1.meijigakuin.ac.jp \newline
Keywords: 
local moves on 1-knots, 
local moves on high dimensional knots, 
crossing-changes on 1-links, 
pass-moves on 1-links, 
products of knots, 
pass-moves on high dimensional links, 
twist-moves on high dimensional links,   
branched cyclic covering spaces,  
Seifert hypersurfaces, 
Seifert matrices. }

\date{}


\begin{abstract} 
We show a relation between 
products of knots, 
which are generalized from the theory of 
isolated singularities of complex hypersurfaces,  
and local moves on knots in all dimensions.  
We discuss the following problem. 
Let $K$ be a 1-knot which is obtained from another 1-knot $J$ by a crossing change (resp. pass-move). 
For a given knot $A$, 
what kind of relation do the products of knots,  
$K\otimes A$ and  $J\otimes A$, have? 
We characterize these kinds of relation between $K\otimes A$ and  $J\otimes A$
by using local moves on high dimensional knots. 
We also discuss a connection between local moves and knot invariants.  
We show a new form of identities for knot polynomials associated with a local move. 
\end{abstract}

\maketitle

\section{Introduction}\label{Introduction} 





Let $f: C^n \longrightarrow C$ be a (complex) polynomial mapping with an isolated singularity at the 
origin of $C^n$. That is, $f(0) = 0$ and the complex gradient $df$ 
has an isolated zero at the origin.
The {\em link} of this singularity is defined by the formula $L(f) = V(f) \cap S^{2n-1}.$ Here the symbol
$V(f)$ denotes the variety of $f$,       
and $S^{2n-1}$ is a sufficiently small sphere about the origin of $C^{n}.$
\bigbreak

Given another polynomial $g:C^{m} \longrightarrow C,$ form $f + g$ with domain $C^{n+m} = 
C^n \times C^m$ and consider $L(f + g) \subset S^{2n + 2m + 1}.$
\bigbreak

We use a topological construction for $L(f + g) \subset S^{2n + 2m + 1}$ in terms of 
$L(f) \subset S^{2n+1}$ and $L(g) \subset S^{2m+1}.$ The construction generalizes the algebraic
situation. Given nice (to be specified below) codimension-two embeddings $K \subset S^n$ and 
$L \subset S^m $, we form a product $K \otimes L \subset S^{n+m+1}.$ Then 
$L(f) \otimes L(-g) \simeq L(f + g).$
\bigbreak

We will recall and use in this paper a product operation on knots in all dimensions that generalizes this 
result about singularities  
 \cite{Kauffman, KauffmanNeumann, Kauffmanon}   
We will also associate geometric equivalence relations,     
{\it crossing changes} and {\it pass equivalence} \cite{Kauffmanon}  
of classical knots,  
with local moves on high dimensional knots and links, 
which were defined and 
have been  researched 
in 
\cite{Ogasa98n, Ogasa02, Ogasa04, Ogasa07, Ogasa09, OgasaT3, OgasaIH}, 
and relate this to the knot product construction and to the Arf invariant, 
the signature, and knot polynomials  in higher dimensions. 

\smallbreak
Furthermore we show a new form of identities for knot polynomials associated with a local move 
(Theorem \ref{pepper}
). 

\smallbreak
We will show examples of twist-moves on high dimensional knots 
after Theorem \ref{Tokyo} in \S\ref{thecr}, 
after Theorem \ref{skein2} in \S\ref{thepol}, and 
after Theorem \ref{pepper} in \S\ref{thepol}. 
We will show examples of pass-moves on high dimensional knots 
in \S\ref{LMH} and in \S\ref{overview}.

\subsection{Construction of Products}   
In this subsection and the next, we describe the results in references 
 \cite{Kauffman, KauffmanNeumann, Kauffmanon}.  
All manifolds will be smooth. Each ambient sphere $S^n$ comes equipped with an orientation.
A {\em knot} in $S^n$ is any closed oriented codimension-two submanifold $K.$ Given a knot
$K \subset S^n$ we may write $S^n = E_{K} \cup (K \times D^2)$ where $E_{K}$ is a manifold with 
boundary equal to $K \times S^1.$ If $n$ is larger than $3,$ we assume that $K$ is connected. Thus, by
Alexander duality, $H^{1}(E_{K}) \simeq Z.$ One may choose $\phi : E_{K} \longrightarrow S^1$ 
representing the generator of $H^{1}(E_{K})$ so that $\phi$ is differentiable and $\phi|\partial E_{K}$
is a projection on the second factor. If $n=3,$ then $K$ may consist of a collection of disjointly embedded circles. A choice of orientations for these circles determines $\phi$ so that $\phi^{-1}$
applied to a regular value is an oriented spanning surface for $K$ which induces the chosen orientations on each component.
\bigbreak

A knot is said to be {\em spherical} if it is homeomorphic to  sphere.
A knot is said to be {\em fibered} if there is a choice of $\phi$ as above so that 
$\phi: E_{K} \longrightarrow S^1$ is a locally trivial smooth fibration.
\bigbreak

Now suppose that we are given knots $K \subset S^n$ and $K \subset S^m$ and corresponding 
maps $\phi: E_{K} \longrightarrow S^1$ and $\psi:E_{L} \longrightarrow S^1.$ If one knot is fibered,
then $E_{K} \times_{S^{1}}E_{L} = \{(x,y) \in E_{K} \times E_{L} | \phi(x) = \psi(y) \}$ is a well-defined
smooth manifold with boundary. Henceforth, when dealing with a pair of knots, we shall assume that at least one knot is fibered. We now define a manifold $K \otimes L$ and, using its properties, obtain the product knot $K \otimes L \subset S^{n+m+1}.$
\bigbreak

\noindent {\bf Definition.} Given knots $K$ and $L$ as above, we define the manifold
$$K \otimes L = (K \times D^{m+1}) \cup (E_{K} \times_{S^{1}}E_{L})\cup(D^{n+1} \times L).$$
These three pieces are attached according to the following description: Note that 
$$\partial(E_{K} \times_{S^{1}} E_{L}) = (K \times E_{L}) \cup (E_{K} \times L)$$ and 
$$\partial(K \times D^{m+1}) = (K \times D^2 \times L)\cup(K \times E_{L}),$$
$$\partial(D^{n+1} \times L) = (K \times D^2 \times L) \cup(E_{K} \times L).$$
Using these boundary identifications, glue the three pieces together to form a closed manifold.
The manifold $K \otimes L$ is independent of the choices of maps $\phi$ and $\psi$ used in its
construction. 
\bigbreak

Now given $\phi:E_{L} \longrightarrow S^{1},$ there is an embedding 
$\hat{\phi}:E_{L} \longrightarrow  D^{m+1} \times S^{1}$ given essentially by 
$\hat{\phi}(x) = (x, \phi(x)).$ This induces an embedding $K \otimes L \subset K \otimes S^{m} 
\simeq S^{n+m+1}.$ This embedding is well-defined up to ambient isotopy and commutative in
$K$ and $L.$ In this way, we obtain a differential topological generalization of the link of the sum of two isolated singularities. In the next section we will make clear how this generalizes the  
links of singularities.
\bigbreak

\subsection{The Pullback Description for Knot Products}
In the previous section we gave a description of the knot product construction in terms of the 
map $\phi:E_{K} \longrightarrow S^{1}$ to the circle associated with the complement of the tubular neighborhood of a knot. In the case of fibered knots this map is a fibration over the circle. For an arbitrary knot we will call $\phi$ the {\em classifying map} for the knot $K \subset S^{n}.$ In this section we use the classifying maps to construct maps of
balls to the $2$-disk that can participate in a pull-back construction for the knot product.
\bigbreak

Given any map $f:S^{n} \longrightarrow D^{2},$ we can extend it to a map, the {\em cone} on $f$,
$$cf: D^{n+1} = CS^{n} = \{tu | 0\le t \le 1, u \in S^{n} \} \longrightarrow D^{2}$$ 
defined by the formula $cf(tu) = tf(u)$ where $0\le t \le 1$ and 
$u$ is a unit vector in $R^{n+1}.$  
\bigbreak

Let $\phi:E_{K} \longrightarrow S^{1}$  be a classifying map for the knot $K \subset S^{n}.$  
Extend $\phi$ to a map $\phi_{1}: S^{n} = E_{K} \cup (K \times D^{2})  \longrightarrow D^{2}$ by defining it on $K \times D^{2}$ to be the cartesian projection to $D^{2}.$ Now extend $\phi_{1}$ to the cone 
and call this map $\hat{\phi}_{K}$, {\em the cone map for $K$}. 
$$\hat{\phi}_{K} = c\phi_{1}: D^{n+1} \longrightarrow D^{2}.$$
\bigbreak

\begin{figure}
     \begin{center}
     \includegraphics[width=8cm]{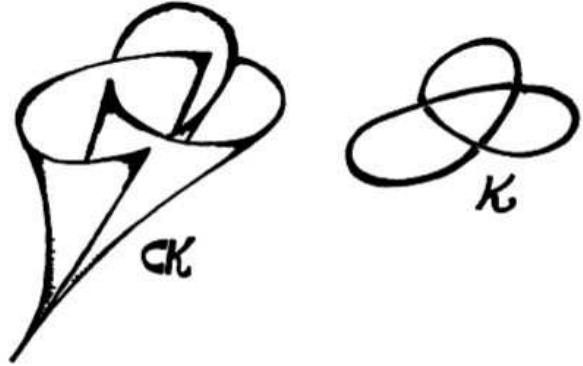}
     \caption{Cone on K}
     \label{cone}
\end{center}
\end{figure}

The point about the construction of the cone map for a given knot $K \subset S^{n}$ is that it produces 
a differential topological analog of an algebraic singularity whose link is this knot. In particular, we have
that $\hat{\phi}_{K}^{-1}(0) = CK \subset D^{n+1} = CS^{n},$ and this mimics the topology of an isolated singularity. See Figure~\ref{cone}.
\bigbreak

Let $\phi:E_{K} \longrightarrow S^{1}$ and $\psi:E_{L} \longrightarrow S^{1}$ be classifying maps for the knots $K \subset S^{n}$ and $L \subset S^{m}.$ Assume that $\psi$ is a fibration over the circle, giving a fibered structure for $L \subset S^{m}.$ Let $X[K,L] \subset D^{n} \times D^{m}$ be the pull-back as shown below.
\[ 
\begin{CD}
X[K,L]            @>j>>          D^{m+1}\\
@ViVV                       @V\hat{\psi}_{L}VV\\
D^{n+1}             @>\hat{\phi}_{K}>>         D^{2}
\end{CD}
\]
The pull-back $X[K,L]$ is the following subset of $D^{n+1} \times D^{m+1}:$
$$X[K,L] = \{ (x,y) \in D^{n+1} \times D^{m+1} |  \hat{\phi}_{K}(x) - \hat{\psi}_{L}(y) = 0 \}.$$
Thus $X[K,L]$ is the differential topological analog of the variety of the sum or difference of two polynomials. Just as with a variety with an isolated singularity, $X[K,L]$ has a singularity at the origin, but the boundary $$\partial X[K,L] \subset \partial(D^{n+1} \times D^{m+1}) \simeq S^{n+m+1}$$
is a smooth submanifold of the $n+m+1$-sphere and this embedding 
$$\partial X[K,L] \subset  S^{n+m+1}$$
is the same as the knot product defined in the previous section. That is, we have that 
$$\partial X[K,L] \simeq K \otimes L$$ and the embeddings are equivalent. It is by way of this pull-back construction that one can prove that indeed the knot product does generalize the link of the sum of isolated complex hypersurface singularities.
\bigbreak

The simplest example of the pull-back construction is given by the following diagram.
\[ 
\begin{CD}
X[a,b]            @>j>>          D^{2}\\
@ViVV                       @V[b]VV\\
D^{2}             @>[a]>>         D^{2}
\end{CD}
\]

In this diagram we have indicated the knot product construction in its lowest dimensional case.
The maps on the disks are of the form $[n]: D^{2} \longrightarrow D^{2}$ where $[n](z) = z^{n}$ where
$n$ is a natural number and $z$ is a complex variable. We take $D^{2}$ as the unit disk in the complex 
plane. Then the maps on spheres in this case are maps of degree $a$ and degree $b$ from circles to themselves. The individual knots are empty and the spanning manifolds consist in $a$ and $b$ points
respectively. We refer to $[a]$ and $[b]$ (regarding the restriction to the circles as defining the maps) as the {\it empty knots of degree $a$ and degree $b$.} We see that 
$$\partial(X[a,b]) = [a] \otimes [b] \subset S^{3}$$ is the corresponding knot product and it is easy to see that {\it $[a] \otimes [b]$ is a torus link of type $(a,b)$.} Continuing in this vein one discovers that the 
Brieskorn manifolds \cite{Kauffmanon}  
the links of singularities 
$$\Sigma(a_{1}, \cdots, a_{n}) = L(z_{1}^{a_{1}} + \cdots + z_{n}^{a_{n}}) \subset S^{2n - 1},$$
are given by the formula 
$$\Sigma(a_{1}, \cdots, a_{n}) = [a_{1}] \otimes \cdots \otimes [a_{n}] \subset S^{2n - 1}.$$
in other words the Brieskorn manifolds and their embeddings in spheres are constructed as products of empty knots of chosen degrees. This completes our description of the elements of the knot product 
construction.
\bigbreak

\subsection{Passing Bands in Low and High Dimensions}\label{PB}
In three dimensions a {\it bandpass} is a replacement of one band crossing over another band by that 
band crossing underneath the other band. 
See the following figure 
for an illustration.

\smallbreak
\includegraphics[width=10cm]{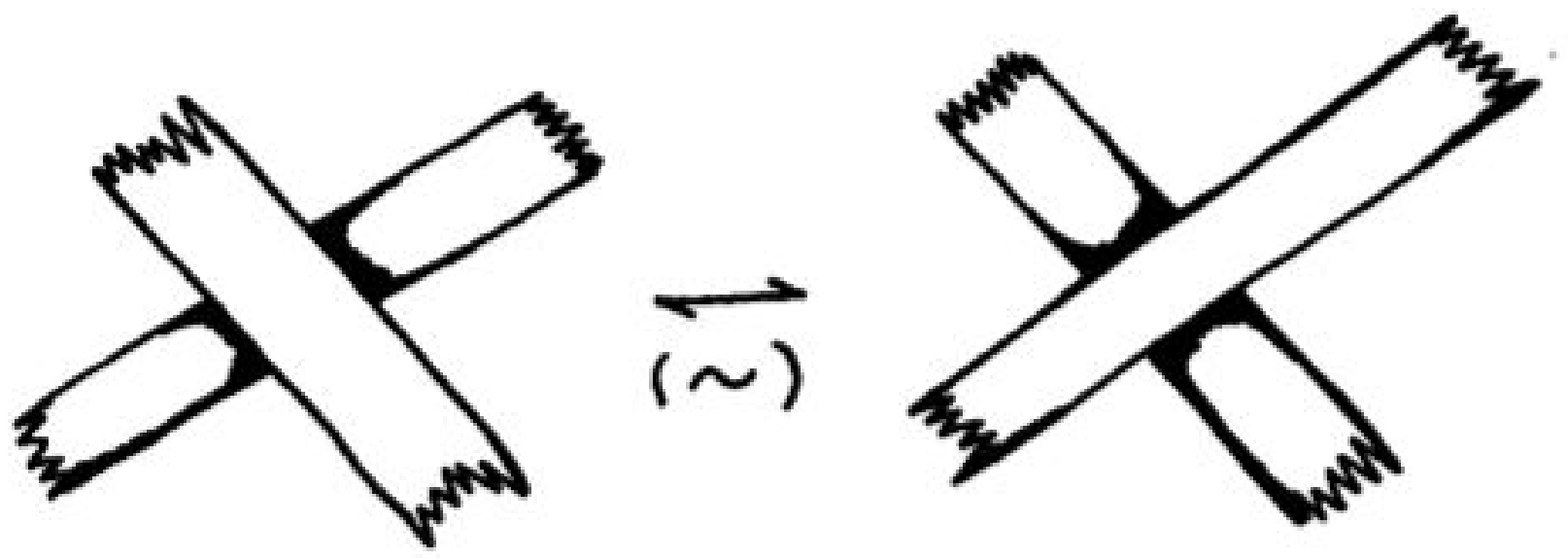}  
\smallbreak \hskip4cm Figure \ref{Introduction}.2 Band Pass
\smallbreak

\noindent 
We usually 
assume that the bands are part of oriented surfaces spanning a link. This means that the local orientation on the two edges of each band are in opposite directions. We say that two oriented knots 
or links are {\it pass equivalent} if one can be obtained from another by a sequence of ambient isotopies and band passes. It is not necessary to construct spanning surfaces for the links in order to perform the band passes, since this is a local operation on the diagrams. The surface interpretation is 
useful for proving facts about pass equivalence. One can show that any classical knot is pass equivalent to either the unknot or the trefoil knot. One can also show that two classical knots are pass equivalent
if and only if they have the same Arf invariant. See 
\S\ref{local1} and  
\cite{Kauffmanon} 
for more information on this topic.

In this paper we will 
relate crossing changes and  
pass equivalence to 
local moves on higher dimensional knots and links and interrelate them with the knot product.
Furthermore we show a connection between the local moves and invariants and polynomials of high dimensional knots.

\bigbreak
We end this subsection with an example.

\hskip3cm     \includegraphics[width=9cm]{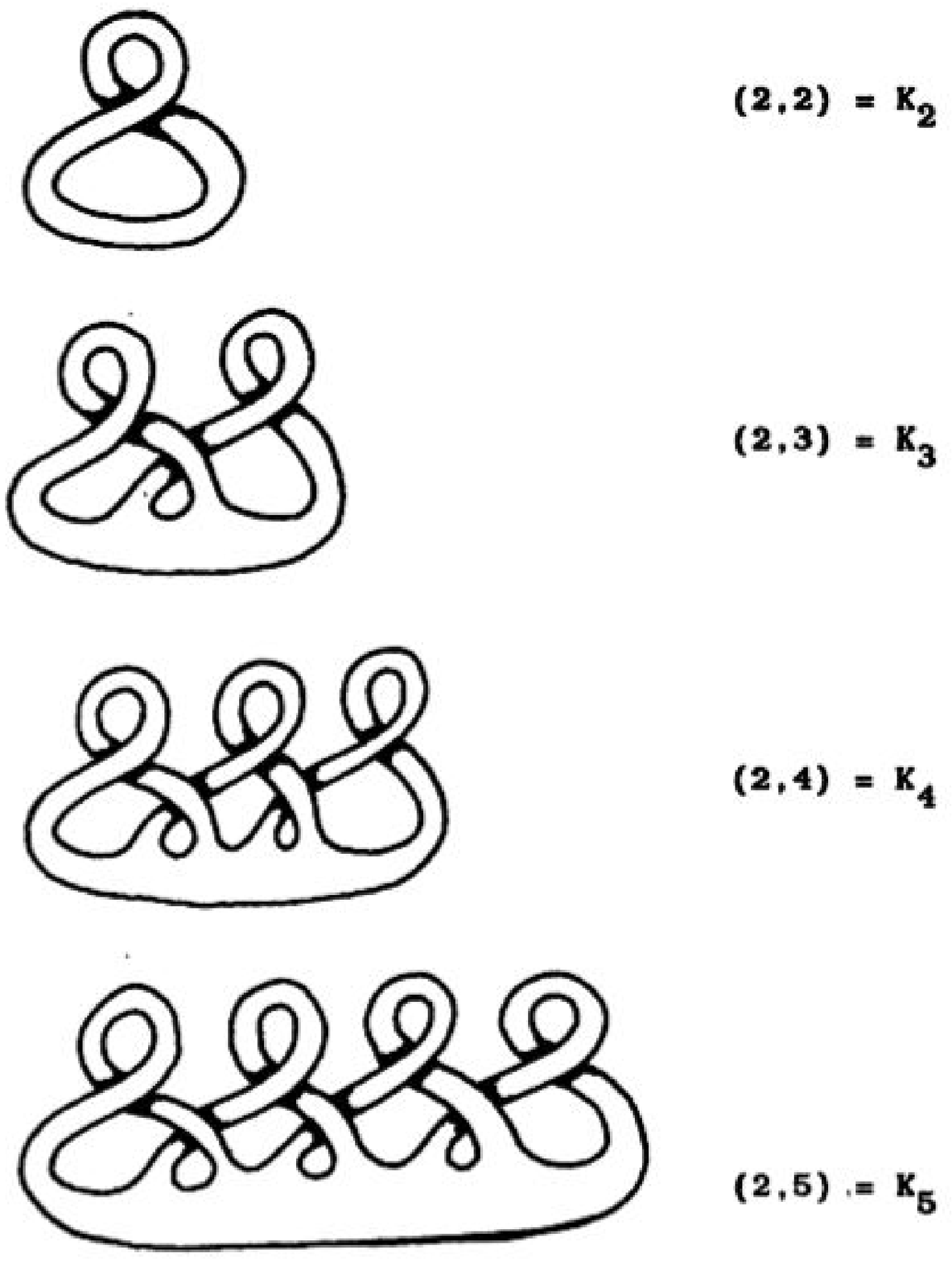}

\hskip2cm   Figure \ref{Introduction}.3. 
\quad  The {$(2,k)$ Torus Knots in Band Representation}
\bigbreak
    
\noindent 
This example is given in more detail in \cite{Kauffmanon} but here we can point to our results in this paper that make the low-dimensional band-passing that we are about to discuss, actual high-dimensional band-passing that accomplishes these results in high-dimensional
manifolds. The result is an $8$-fold periodicity in the list of Brieskorn manifolds 
$\Sigma(k,2,2,2,\cdots,2)$ where there are an odd number of $2$'s. Let $\Sigma_{k}^{4n+1}$ denote
such a Brieskorn manifold with $2n+1$ symbols that are $2$'s. Then $\Sigma_{k}^{4n+1}$
bounds a handle-body whose structure is analogous to 
the spanning surface for a $(2,k)$ torus link, and 
the operation of band-exchange results in a diffeomorphism of this handle-body, hence a diffeomorphism of its boundary. 
See the above figure 
for an illustration of the $(2,k)$ torus links, here called $K_{k}$ and the banded surfaces that bound these links.  
In \cite{Kauffmanon} we exploit this relationship with the low-dimsional topology 
to prove by band-passing that $K_{k+8}$ is pass-equivalent to $K_{k}$, 
and so prove, up to diffeomorphism that the list of manifolds $\Sigma_{k}^{4n+1}$ is 
periodic of period $8$ in $k.$ By applying the results of this paper, 
we can make this conclusion 
directly by using the higher dimensional versions of pass-moves. 
 (Outline of the proof: 
A $(4n+1)$-submanifold $\mathcal K_k=K_{k}\otimes^{n}\text{(the Hopf link)}$ in $S^{4n+3}$
is diffeomorphic to $\Sigma_{k}^{4n+1}$. 
$\mathcal K_k$ is high dimensional pass-move-equivalent to $\mathcal J$ with the following properties: 
A Seifert matrix associated with a Seifert hypersurface $V_{\mathcal J}$ for $\mathcal J$ 
is the same as 
a Seifert matrix associated 
with a Seifert hypersurface $V_{\mathcal K_{k+8}}$ for $\mathcal K_{k+8}$.  
$V_{\mathcal J}$ and $V_{\mathcal K_{k+8}}$ consist of 
a $(4n+2)$-dimensional 0-handle and 
$(4n+2)$-dimensional $(2n+1)$-handles. 
Of course $V_{\mathcal J}$ and $V_{\mathcal K_{k+8}}$  are compact oriented parallelizable and 
have the same intersection matrix on the $(2n+1)$-th homology groups.    
\quad 
Therefore, by surgery theory, $\Sigma_{k}^{4n+1}$ is diffeomorphic to $\Sigma_{k+8}^{4n+1}$. )
The details of this band exchange, illustrated in three-dimensions are interesting, and we refer the reader to \cite{Kauffmanon} for more about this aspect of the example.
We could investigate $\Sigma(a,b,2,2,2,\cdots,2)$, where there are an even number of $2$'s, 
by using high dimensional pass-moves 
because the $(a,b)$ torus knots are classified by pass-equivalence. 

\bigbreak

We will show an example of pass-moves on high dimensional knots in \S\ref{LMH}. 
We will show examples of twist-moves on high dimensional knots 
after Theorem \ref{Tokyo} in \S\ref{thecr}, 
after Theorem \ref{skein2} in \S\ref{thepol}, and 
after Theorem \ref{pepper} in \S\ref{thepol}.

\subsection{Organization and the Main Problem}
It is natural to consider the following problem: 
Let $K$ be a 1-knot which is obtained from another 1-knot $J$ by a crossing change (resp. pass-move). 
For a given knot $A$, 
what kind of relation do $K\otimes A$ and  $J\otimes A$ have? 
In this paper we characterize these kinds of relation between $K\otimes A$ and  $J\otimes A$
by using local moves on high dimensional knots. 
This paper is organized as follows. 
\smallbreak

\noindent
\S\ref{Introduction}  Introduction

\noindent
\S\ref{Mainresults} Main results 

\noindent
\S\ref{local1} {Local moves on  classical links}

\noindent
\S\ref{localn} {Local moves on $n$-knots}

\noindent \S\ref{knotproduct} {Products of knots}

\noindent
\S\ref{Alex}  
{Review of the $\Q[t, t^{-1}]$-Alexander polynomials for $n$-knots and 
$n$-dimensional closed 

oriented submanifolds
}

\noindent
\S\ref{RIL}  
{Some results on invariants of $n$-knots and local moves on $n$-knots}

\noindent
\S\ref{thecr} {Theorems on relations between crossing-changes and  knot products}

\noindent
\S\ref{thepa} {Theorems on relations between pass-moves and knot products}

\noindent
\S\ref{thepol}
 {Theorems on relations between local move identities of a knot  polynomial and knot 

products}

\noindent
\S\ref{specialX}
A remark on the $\Z[t, t^{-1}]$ case

\noindent
\S\ref{procr} {Proof of Theorems in \S \ref{thecr}}

\noindent
\S\ref{propa} {Proof of Theorems in \S \ref{thepa}

\noindent
\S\ref{propol} {Proof of Theorems in \S \ref{thepol}}

\noindent
\S\ref{open} {A problem}

\smallbreak
By considering this problem of the effect on higher dimensional knots 
of changes from lower dimensions, 
via knot products, 
we raise many questions that deserve further investigations.

\subsection{Local Moves on High Dimensional Knots}\label{LMH}
Local moves on high dimensional knots were defined and 
have been  researched 
in 
\cite{Ogasa98n,  Ogasa02,  Ogasa04, Ogasa07,   Ogasa09, OgasaT3, OgasaIH}. 
We review the definition of local moves on high dimensional knots after showing an example.

We will show other examples 
 in \S\ref{overview},  
after Theorem \ref{Tokyo} in \S\ref{thecr}, 
after Theorem \ref{skein2} in \S\ref{thepol}, and 
after Theorem \ref{pepper} in \S\ref{thepol}.

\bigbreak
\noindent{\bf Lemma.} {\it  
Letting $B^p$ denote a $p$-dimensional ball, 
we can write 
$$S^p=B^p_u\cup B^p_d$$ 
$$S^p\x S^q
=(B^p_u\cup B^p_d)\x(B^q_u\cup B^q_d).$$ 

Thus 
$$S^p\x S^q=(B^p_u\x B^q_u)\cup(B^p_u\x B^q_d)\cup(B^p_d\x B^q_u)\cup(B^p_d\x B^q_d).$$
}

\noindent{\bf Proof.}  
Use the fact that 
$$(X\cup Y)\times Z=(X\times Z)\cup(Y\times Z).$$\qed

\bigbreak 
Now let 
$$F= (S^p\x S^q)-\text{Int}(B^p_u\x B^q_u).$$  
We indicate $F$ in the figure below and abbreviate $B^\sharp_\star$ to  $B_\star$.





\vskip5mm
\includegraphics[width=12cm]{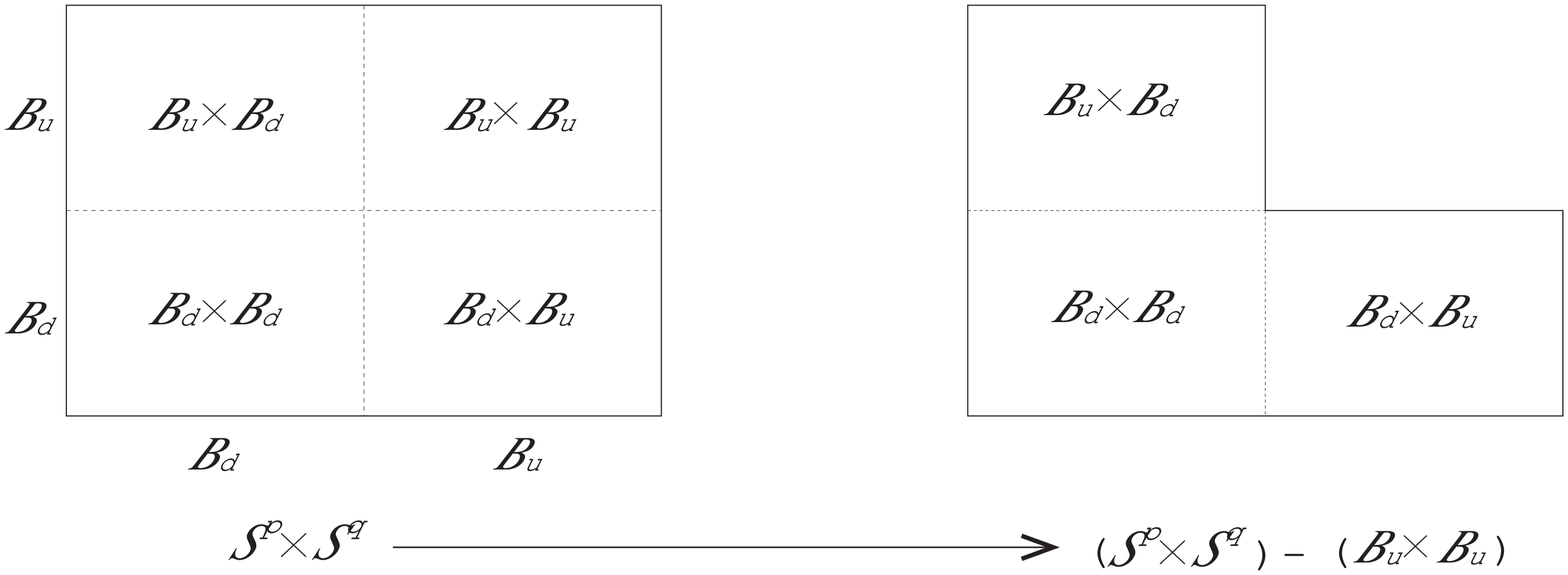}  
\vskip5mm

$F$ is drawn 
in another way as below. 
Note that we can bend the corner of $B^p_u\x B^q_u$ 
and change it into the $(p+q)$-dimensional ball. 
Let $p+q=n+1$. 
Hence the boundary of $F$ is $S^n$.

\vskip5mm
\includegraphics[width=12cm]{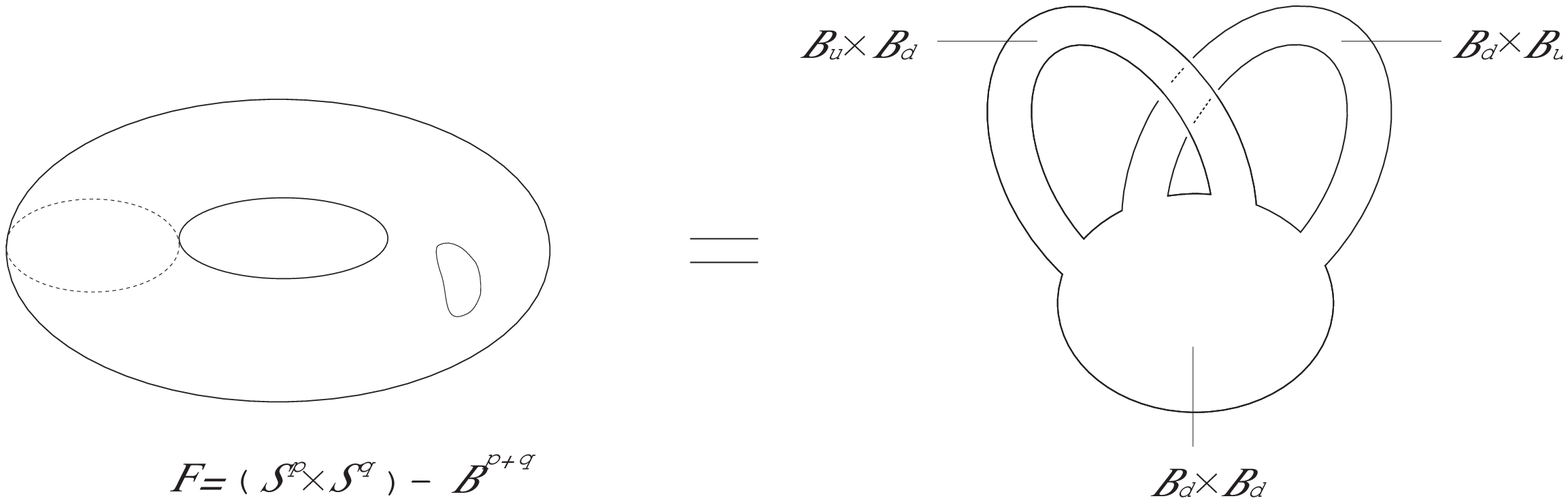}  

\smallbreak
\hskip4cm   Figure \ref{LMH}.1: $(S^p\x S^q)-B^{p+q}$

\vskip5mm

We can regard 
$(B^p_d\x B^q_d)$ as a $(p+q)$-dimensional 0-handle, \newline
$(B^p_u\x B^q_d)$ as a $(p+q)$-dimensional $p$-handle, and  \newline
$(B^p_d\x B^q_u)$ as a $(p+q)$-dimensional $q$-handle.

\noindent
Take $F$ in $S^{n+2}$. 
This is indicated in the figure below. 

\bigbreak
\hskip3cm\includegraphics[width=6cm]{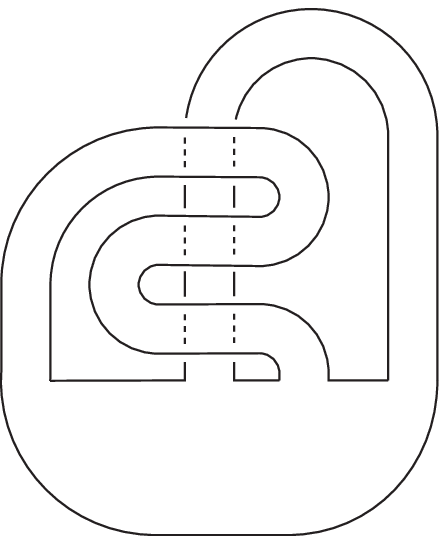}   
\smallbreak
\hskip4cm   Figure \ref{LMH}.2: A trivial $n$-knot 

\bigbreak

The boundary of $F$ in  $S^{n+2}$ is an $n$-knot.  
Furthermore it is the trivial $n$-knot.

\noindent
Carry out a `local move' on this $n$-knot 
in an $(n+2)$-ball, which is denoted by a dotted circle  
in the following figure.

\bigbreak
\hskip4cm\includegraphics[width=6cm]{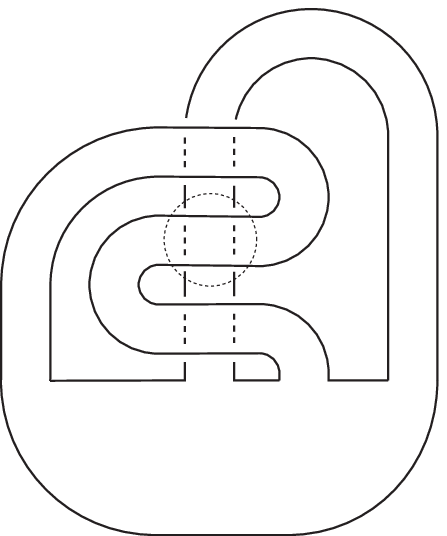}  
\smallbreak
\hskip2cm   Figure \ref{LMH}.3: A local move will be carried out in the dotted $(n+2)$-ball.   

\hskip45mm The resulting $n$-knot is a nontrivial $n$-knot. 
\bigbreak


\bigbreak
\hskip4cm\includegraphics[width=6cm]{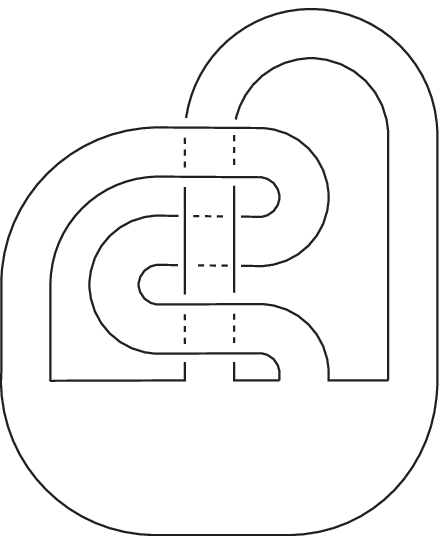}  
\smallbreak
\hskip3cm   Figure \ref{LMH}.4: A nontrivial $n$-knot

\bigbreak

We can prove this fact by using Seifert matrices and the Alexander polynomial. 
In the proof,  we use the fact that $S^p$ and $S^q$ can be `linked' in $S^{p+q+1}$. 
Recall that $p+q+1=n+2$. Note that $S^q$ and $S^p$ are included in $F$ as shown below.

\bigbreak
\hskip4cm\includegraphics[width=7cm]{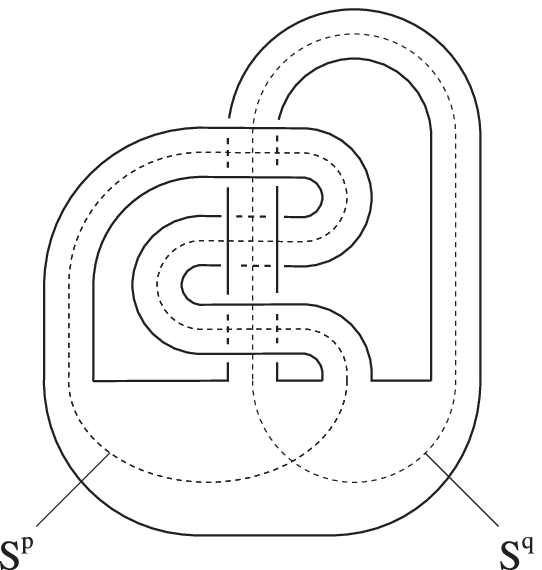}  
\smallbreak
\hskip2cm   Figure \ref{LMH}.5: $S^p$ and $S^q$ in $F$ whose boundary is the $n$-knot

\bigbreak

Note that the above operation is done only in an $(n+2)$-ball.  
This operation is an example of $(p, q)$-pass-moves.

Local moves on high dimensional knots 
(or codimension two submanifolds in a high dimensional standard sphere) that we discuss are 
twist-moves,  $(p, q)$-pass-moves, and ribbon-moves.  
These are defined and are discussed in 
\cite{Ogasa98n, Ogasa02,Ogasa04, Ogasa07, Ogasa09, OgasaT3, OgasaIH} 
(see \S\ref{localn}). 
These moves can change the trivial knot into nontrivial knots. 
Conversely they can make  some nontrivial knots into the trivial knot.  
We will show examples of twist-moves on high dimensional knots 
after Theorem \ref{Tokyo} in \S\ref{thecr}, 
after Theorem \ref{skein2} in \S\ref{thepol}, and 
after Theorem \ref{pepper} in \S\ref{thepol}.

\subsection{$(p, q)$-Pass-Moves}\label{pqPM}
The $(p, q)$-pass-moves on 
$n$-dimensional knots 
(or $n$-dimensional submanifolds in $S^{n+2}$), 
where $p+q=n+1$ and $p, q, n\in\N$,  
are as follows:  
See \S\ref{localn} for detail.   
Two $n$-dimensional knots (or $n$-dimensional submanifolds in $S^{n+2}$) 
$K_+$ and $K_-$ 
differ only in a $(n+2)$-ball $B^{n+2}$ as shown in the following figures. 
($B$ sometimes denotes $B^{n+2}$.)
$B\cap K_+$ (resp. $B\cap K_-$) is  
diffeomorphic to the disjoint union of 
$S^{n-p}\x D^p$  and 
$S^{p-1}\x D^{n+1-p}$. 
Take an $(n+1)$-dimensional $p$-handle  
(resp. an $(n+1)$-dimensional $(n+1-p)$-handle) which is 
embedded trivially in $B^{n+2}$ and attached to $\partial B^{n+2}$ trivially. 
Suppose that $S^{n-p}\x D^p$  
(resp. $S^{p-1}\x D^{n+1-p}$) 
 is 
$\overline{\partial h^p - \partial B^{n+2}}$   
(resp. $\overline{\partial h^{n+1-p} - \partial B^{n+2}}$). 
We pass    
the handle $h^p$ across $h^{n+1-p}$as in the following figures, 
 move 
$S^{n-p}\x D^p$  and  $S^{p-1}\x D^{n+1-p}$ together 
and change $K_+$ into $K_-$.

We showed an example of pass-moves on high dimensional knots 
in \S\ref{LMH} and will show one in \S\ref{overview}.

\np
Figure \ref{pqPM}.1 is a diagram 
of a $(p, n+1-p)$-pass-move.
Here, $h^p_+$ denotes $h^p$.  
Figure \ref{pqPM}.2, which consists of the two figures (1) (2),  
is another diagram 
of a $(p, n+1-p)$-pass-move.  
Note that $(1,1)$-pass-moves are same as pass-moves on 1-knots (1-links).

\smallbreak
\unitlength 0.1in
\begin{picture}(57.35,27.04)(1.70,-34.67)
%
\special{pn 8}%
\special{ar 4968 1942 937 937  0.8932296 6.2831853}%
\special{ar 4968 1942 937 937  0.0000000 0.8617069}%
%
\special{pn 8}%
\special{pa 4236 1373}%
\special{pa 5415 2755}%
\special{fp}%
%
\special{pn 8}%
\special{pa 5709 2500}%
\special{pa 4541 1118}%
\special{fp}%
%
\special{pn 8}%
\special{pa 5415 1118}%
\special{pa 4998 1515}%
\special{fp}%
\special{pa 4612 1922}%
\special{pa 4134 2369}%
\special{fp}%
\special{pa 5700 1382}%
\special{pa 5242 1789}%
\special{fp}%
\special{pa 4835 2216}%
\special{pa 4388 2643}%
\special{fp}%
%
\special{pn 8}%
\special{ar 1817 1980 937 937  5.5898109 6.2831853}%
\special{ar 1817 1980 937 937  0.0000000 5.5597350}%
%
\special{pn 8}%
\special{pa 1258 2722}%
\special{pa 2622 1521}%
\special{fp}%
%
\special{pn 8}%
\special{pa 2365 1230}%
\special{pa 1000 2420}%
\special{fp}%
%
\special{pn 8}%
\special{pa 986 1546}%
\special{pa 1390 1957}%
\special{fp}%
\special{pa 1802 2337}%
\special{pa 2255 2808}%
\special{fp}%
\special{pa 1246 1258}%
\special{pa 1661 1709}%
\special{fp}%
\special{pa 2093 2109}%
\special{pa 2527 2550}%
\special{fp}%
%
\special{pn 8}%
\special{pa 963 3251}%
\special{pa 1135 2336}%
\special{dt 0.045}%
\special{sh 1}%
\special{pa 1135 2336}%
\special{pa 1103 2398}%
\special{pa 1125 2388}%
\special{pa 1142 2405}%
\special{pa 1135 2336}%
\special{fp}%
\special{pa 1004 3230}%
\special{pa 1379 2611}%
\special{dt 0.045}%
\special{sh 1}%
\special{pa 1379 2611}%
\special{pa 1327 2658}%
\special{pa 1351 2657}%
\special{pa 1362 2678}%
\special{pa 1379 2611}%
\special{fp}%
%
\special{pn 8}%
\special{pa 2650 3027}%
\special{pa 2457 2509}%
\special{dt 0.045}%
\special{sh 1}%
\special{pa 2457 2509}%
\special{pa 2462 2578}%
\special{pa 2476 2559}%
\special{pa 2499 2564}%
\special{pa 2457 2509}%
\special{fp}%
\special{pa 2640 3007}%
\special{pa 2040 2539}%
\special{dt 0.045}%
\special{sh 1}%
\special{pa 2040 2539}%
\special{pa 2080 2596}%
\special{pa 2082 2572}%
\special{pa 2105 2564}%
\special{pa 2040 2539}%
\special{fp}%
\put(24.4700,-33.2200){\makebox(0,0)[lb]{$S^{p-1}\x D^{n+1-p}$}}%
\put(3.2200,-34.3300){\makebox(0,0)[lb]{$S^{n-p}\x D^p$}}%
\put(24.0600,-35.8700){\makebox(0,0)[lb]{$=\overline{\partial h^{n+1-p}-\partial B}$}}%
\put(13.0800,-9.7400){\makebox(0,0)[lb]{$B\cap K_+$}}%
\put(45.0000,-9.3300){\makebox(0,0)[lb]{$B\cap K_-$}}%
\put(1.7000,-36.3700){\makebox(0,0)[lb]{$=\overline{\partial h^p_+-\partial B}$}}%
\end{picture}%
\vskip1cm
\hskip7cm Figure \ref{pqPM}.1

\noindent 
{\bf A $(p, n+1-p)$-pass-move on an $n$-dimensional submanifold  
$\subset S^{n+2}$. 
Note $B=B^{n+2}=D^{n+2}\subset S^{n+2}$. 
This figure includes $h^p_+$ and $h^{n+1-p}$.} 

\vskip1cm

\np

\hskip1mm\vskip10mm
\hskip18mm
\includegraphics[width=10cm]{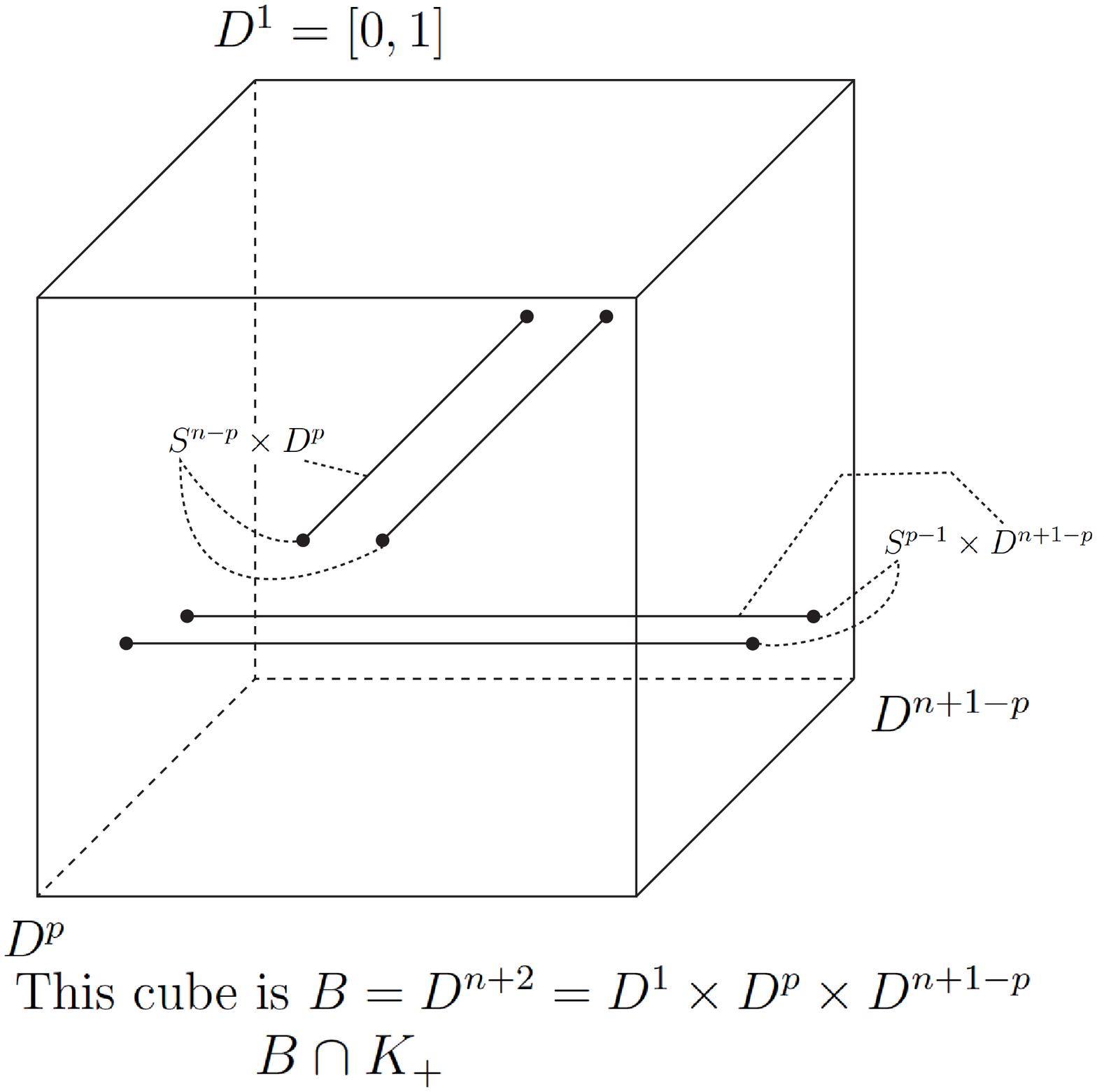}

\noindent

\vskip10mm
\hskip3cm Figure \ref{pqPM}.2.(1)

{\bf A $(p, n+1-p)$-pass-move on an $n$-dimensional submanifold 
$\subset S^{n+2}$.

Note $B=B^{n+2}=D^{n+2}\subset S^{n+2}$.}

\np

\hskip8mm
\includegraphics[width=10cm]{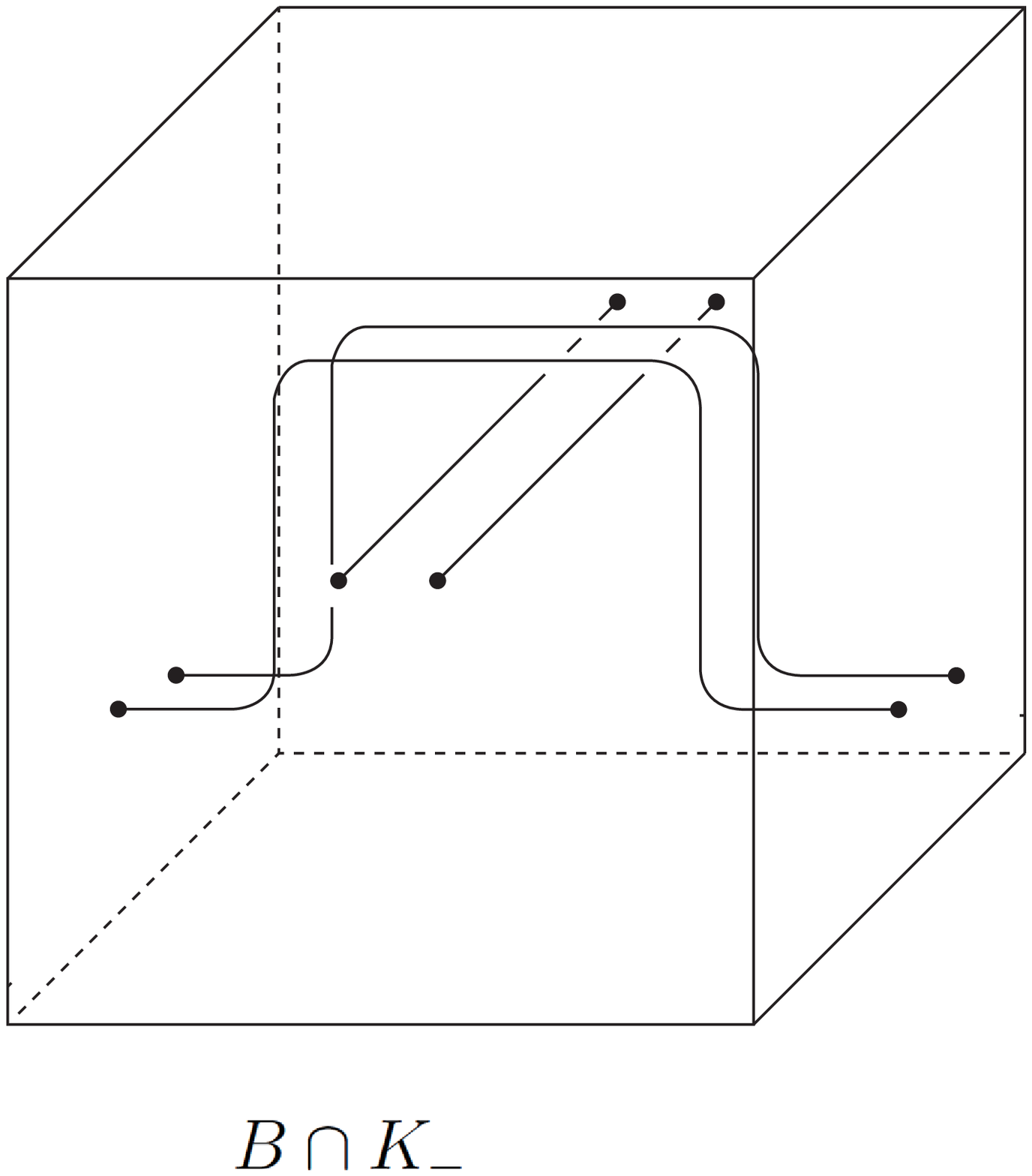}

\vskip-20mm\hskip3cm Figure \ref{pqPM}.2.(2)

{\bf A $(p, n+1-p)$-pass-move on an $n$-dimensional submanifold 
$\subset S^{n+2}$.

Note $B=B^{n+2}=D^{n+2}\subset S^{n+2}$.} 

\np

\subsection{Twist-Moves}\label{TM}
The twist-moves on 
$(2p+1)$-dimensional knots (or codimension two submanifolds in $S^{2p+3}$)
are 
as follows, where $p\in\N\cup\{0\}$. 
See \S\ref{localn} for detail.
Two $(2p+1)$-dimensional knots (or codimension two submanifolds in $S^{2p+3}$) 
$K_+$ and $K_-$ 
differ only in a $(2p+3)$-ball $B^{2p+3}$ as shown in the following figures. 
($B$ sometimes denotes $B^{n+2}$.)
$B\cap K_+$ (resp. $B\cap K_-$) is 
diffeomorphic to 
$S^p\x D^{p+1}$.  
Take a $(2p+2)$-dimensional $(p+1)$-handle  which is 
embedded trivially in $B^{2p+3}$ and attached to $\partial B^{2p+3}$ trivially. 
Suppose that $S^p\x D^{p+1}$  
 is 
$\overline{\partial h^{p+1} - \partial B^{2p+3}}$.    
Note that we can twist  the handle  
$h^{p+1}$ in $B$ once (see \S\ref{localn} for detail), 
 move 
$S^p\x D^{p+1}$  together 
and change $K_+$ into $K_-$. 

\bigbreak
Figure \ref{TM}.1, which consists of the two figures (1) (2), 
is a diagram 
of a twist-move.  
The upper half of Figure \ref{TM}.2 is 
another diagram 
of a twist-move.
Compare 
the upper half of 
of Figure \ref{TM}.2 
and 
the lower half. 
If $p=0$ (i.e., $n=2p+1=1$), 
the left figure 
in the upper half 
and 
that in the lower half 
are same. 
That is, if $p=0$, twist-moves are crossing-changes on 1-links.  
We will show examples of twist-moves on high dimensional knots 
after Theorem \ref{Tokyo} in \S\ref{thecr}, 
after Theorem \ref{skein2} in \S\ref{thepol}, and 
after Theorem \ref{pepper} in \S\ref{thepol}. 
We review another local-move on high dimensional knots (resp. submanifolds), 
ribbon-moves in \S\ref{localn}.

\np

\hskip3cm\includegraphics[width=9cm]{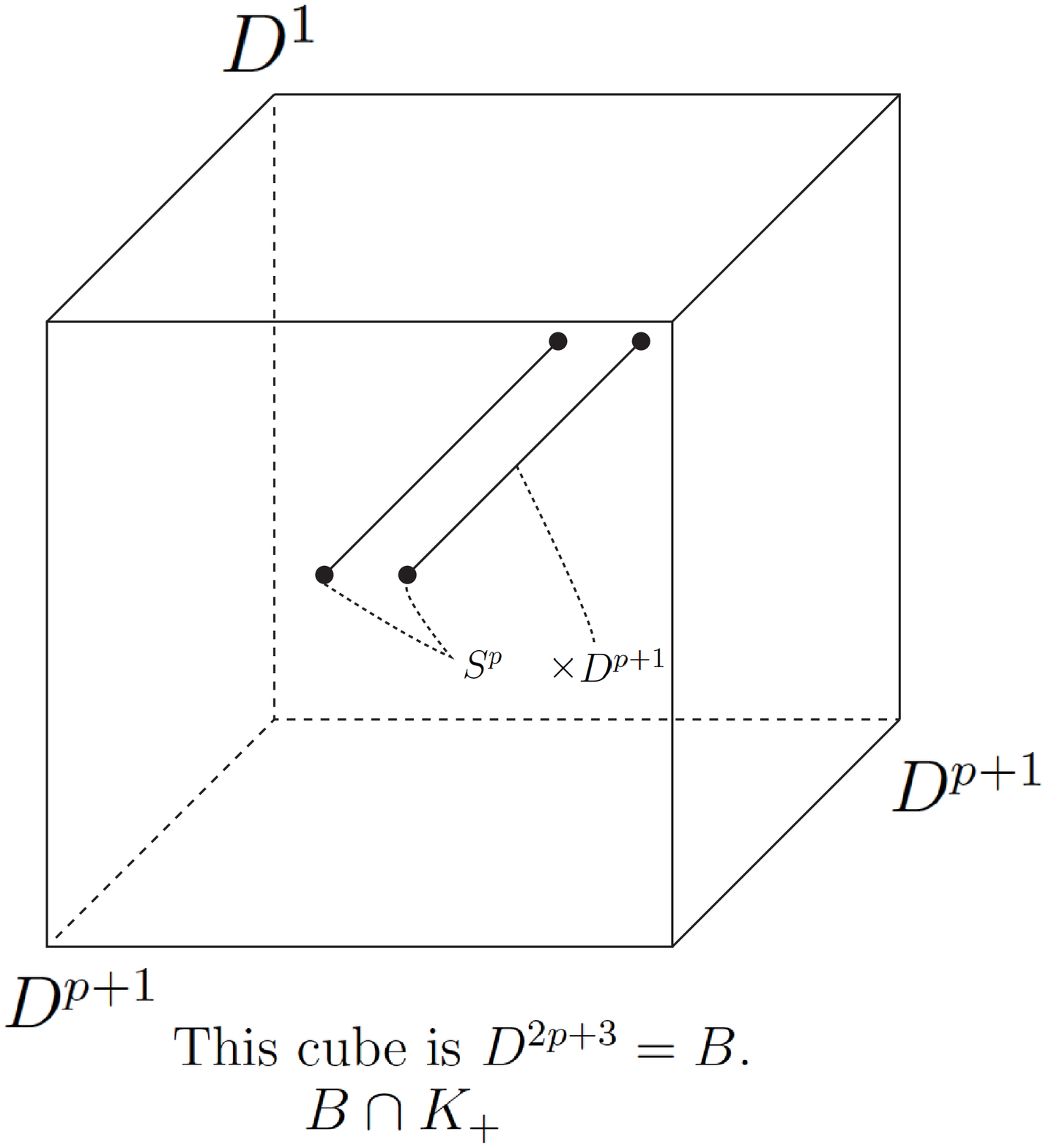}

\vskip6mm
\hskip5cm Figure \ref{TM}.1.(1)

{\bf A twist-move on $(2p+1)$-dimensional submanifold 
$\subset S^{2p+3}$. 

Note $B=D^{2p+3}\subset S^{2p+3}$. }

\np

\hskip3cm\includegraphics[width=9cm]{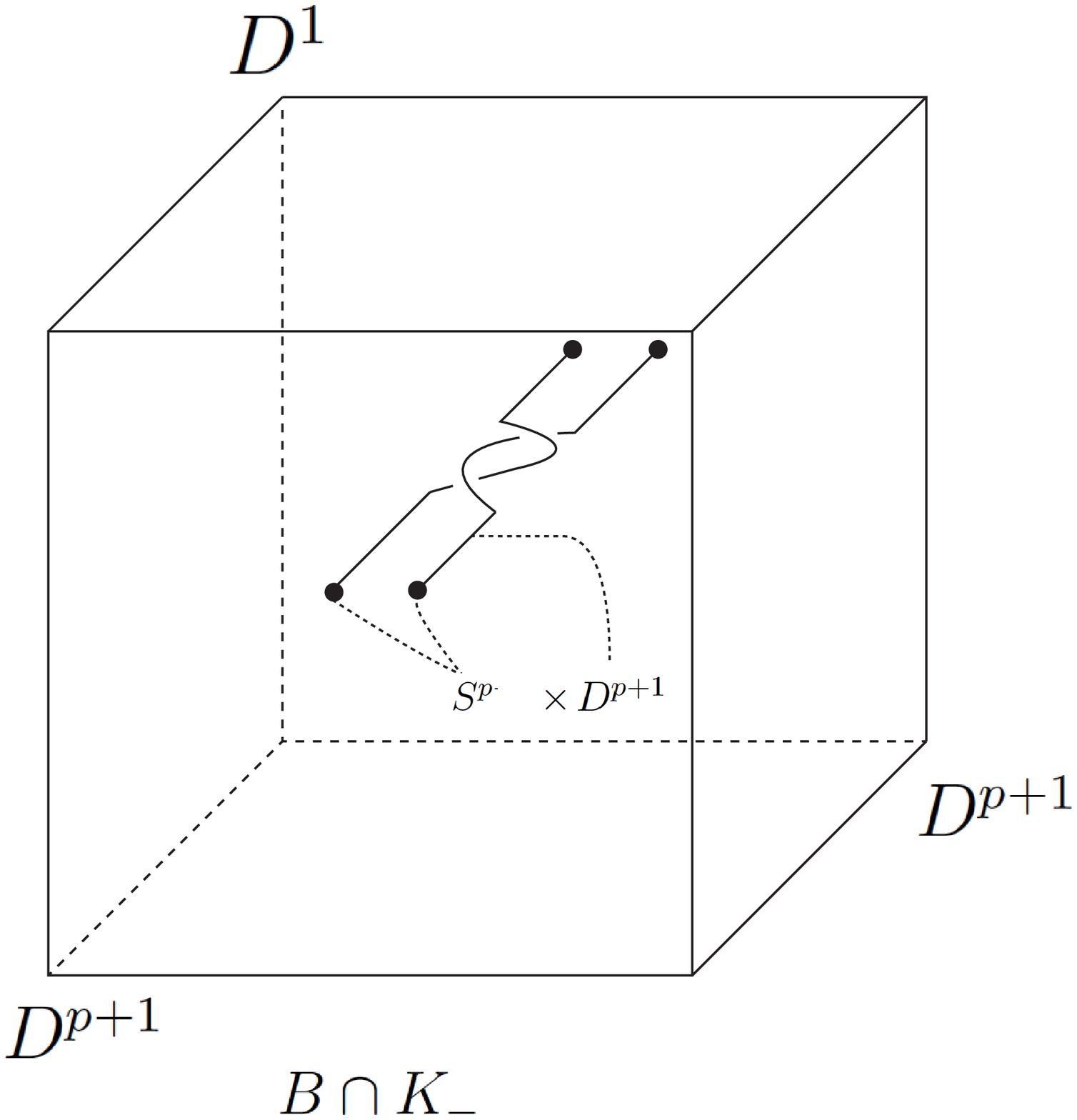}

\vskip10mm
\hskip5cm Figure \ref{TM}.1.(2)

{\bf A twist-move on $(2p+1)$-dimensional submanifold 
$\subset S^{2p+3}$. 

Note $B=D^{2p+3}\subset S^{2p+3}$. }

\np

\hskip4cm\includegraphics[width=8cm]{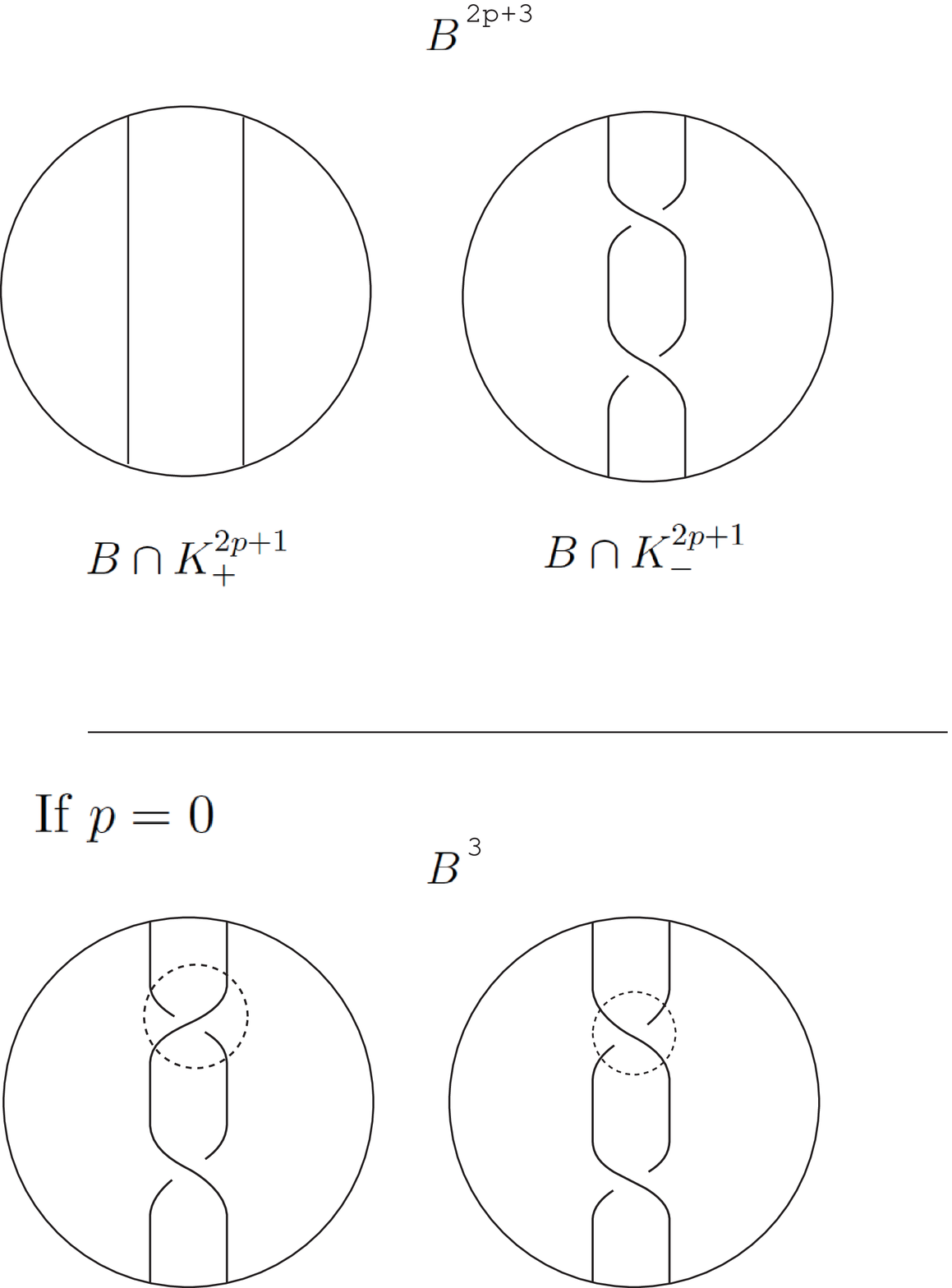}

\vskip2mm
\hskip15mm  The pair of two 
\quad
\unitlength 0.1in
\begin{picture}(2.40,2.40)(0.60,-2.60)
%
\special{pn 8}%
\special{ar 180 140 120 120  0.0000000 0.1000000}%
\special{ar 180 140 120 120  0.4000000 0.5000000}%
\special{ar 180 140 120 120  0.8000000 0.9000000}%
\special{ar 180 140 120 120  1.2000000 1.3000000}%
\special{ar 180 140 120 120  1.6000000 1.7000000}%
\special{ar 180 140 120 120  2.0000000 2.1000000}%
\special{ar 180 140 120 120  2.4000000 2.5000000}%
\special{ar 180 140 120 120  2.8000000 2.9000000}%
\special{ar 180 140 120 120  3.2000000 3.3000000}%
\special{ar 180 140 120 120  3.6000000 3.7000000}%
\special{ar 180 140 120 120  4.0000000 4.1000000}%
\special{ar 180 140 120 120  4.4000000 4.5000000}%
\special{ar 180 140 120 120  4.8000000 4.9000000}%
\special{ar 180 140 120 120  5.2000000 5.3000000}%
\special{ar 180 140 120 120  5.6000000 5.7000000}%
\special{ar 180 140 120 120  6.0000000 6.1000000}%
\end{picture}%
 \quad
makes a crossing-change on a 1-dimensional link.
 
 \vskip10mm 
 \hskip7cm Figure \ref{TM}.2

\hskip25mm
{\bf A twist-move on a 1-link is  a crossing-change on a 1-link.}

\np
\subsection{An overview of the Main Results}\label{overview}
One of our main theorems is the following.  
If a 1-link $K$ is obtained from a 1-link $K'$ by a crossing-change,  
then  
the knot product, $K\otimes$(the Hopf link),   
is obtained from 
the knot product,,  $K'\otimes$(the Hopf link),   
by a twist-move 
(see  Theorem \ref{summerc}, \ref{Illinois}, 
\ref{Tokyo}
).  
Other results in this paper are as follows: 
If a 1-knot $K$ is obtained from a 1-knot $K'$ by a pass-move,  
then the knot product, $K\otimes$(the Hopf link),   
is obtained from 
the knot product, $K'\otimes$(the Hopf link),   
by a $(3,3)$-pass-move 
(see  Theorem \ref{Waltham}). 
Let $K$ and $K'$ be 1-knots. 
The 1-knot $K$ is pass-move-equivalent to the 1-knot $K'$ if and only if 
the knot product, $K\otimes$(the Hopf link),   
is  $(3,3)$-pass-move-equivalent to 
the knot product, $K'\otimes$(the Hopf link),  
(see  Theorem \ref{Waltham} 
and \ref{ki}).  
Of course we show more results in other high dimensional cases.

\section{Main results  ----- 
Technical Statements}\label{Mainresults} 


We work in the smooth category. 
Let $L=(K_1,...,K_m)$ be an $m$-component $n$-(dimensional) submanifold $\subset S^{n+2}$. 
If  $m=1$ and if $L$ is homeomorphic 
to the standard sphere, 
then $L$ is called an {$n$-dimensional spherical knot}. 
(Note that some literature defines spherical knots are PL homeomorphic to the standard sphere.)
If each $K_i$ is a spherical knot, 
then $L$ is called an {$n$-dimensional spherical link}. 
We say that n-submanifolds $L$ and $L'$ are {\it identical} 
if there exists an orientation preserving identity map   \newline  
$id:S^{n+2}$ $\rightarrow$ $S^{n+2}$ such that $id(L)$=$L'$  and 
$id\vert_{L}:L\rightarrow L'$ is 
an orientation and order preserving identity map.    
We say that n-submanifolds $L$ and $L'$ are {\it equivalent} 
if there exists an orientation preserving diffeomorphism 
$f:S^{n+2}$ $\rightarrow$ $S^{n+2}$ such that $f(L)$=$L'$  and 
$f\vert_{L}:L\rightarrow L'$ is 
an orientation and order preserving diffeomorphism.    
An $m$-component $n$-submanifold $L=(L_1,...,L_m)$ is called 
a {\it trivial }($n$-){\it link} 
if each $L_i$ bounds an $(n+1)$-ball $B_i$ trivially embedded in $S^{n+2}$ 
and if $B_i\cap B_j=\phi (i\neq j)$. 
If $m=1$, then $L$ is called a 
{\it trivial }($n$-){\it knot}.

\vskip3mm
The following theorems are special cases of our results.

\begin{thm}\label{summerc}  
Suppose that two 1-links $J$ and $K$ differ by one crossing-change.
Then the knot products, 
$J \otimes^\mu($the Hopf link$)$ 
and   \newline 
$K\otimes^\mu($the Hopf link$)$,   
differ by one twist-move, 
where $\mu\in\N\cup\{0\}$.
Here, 
$A\otimes^\mu B$ 
means 
$A\otimes B...\otimes B$, which includes  $\mu$ copies of $B$, 
where $\mu\in\N\cup\{0\}$. 
\end{thm}

\noindent 
{\bf Note.} 
(1) The fact that the two knots differ 
means 
that the two knots are not identical.  
There are two cases 
that the two knots are  not equivalent 
and that the two knots are  equivalent.

\smallbreak   \noindent (2)   
We review the empty knots $[n]$, where $n\in\N$, in \S\ref{knotproduct}.  
It holds that $[2]\otimes[2]$ is the negative Hopf link in $S^3$ 
(see \cite{Kauffman, KauffmanNeumann} and Theorem \ref{Hopf}.) 
The above  Theorem \ref{summerc} follows from 
the following  Theorem \ref{Illinois} 
by Theorem \ref{Hopf}.
 
\smallbreak   \noindent (3)  
We will show an example of the phenomenon which Theorem \ref{summerc}, \ref{Illinois}, 
and \ref{Tokyo}
assert in \S\ref{thecr}.

\bigbreak\noindent{\bf Theorem \ref{Illinois}.} {\it 
Take the same $J, K$ in Theorem \ref{summerc}.
Then
the $(2\nu+1)$-submanifolds $\subset S^{2\nu+3}$, 
$J\otimes^\nu($the empty knot $[2])$    
and   
$K\otimes^\nu($the empty knot $[2])$,   
differ by one twist-move, 
where $\nu\in\N\cup\{0\}$.
}

\bigbreak\noindent{\bf Theorem \ref{Tokyo}.}  {\it          
Let $m\in\N\cup\{0\}$. 
Suppose that two $(2m+1)$-dimensional closed oriented submanifolds  $\subset S^{2m+3}$, 
$J$ and $K$,  differ by one twist-move.
Then
the $(2m+2\nu+1)$-submanifolds $\subset S^{2m+2\nu+3}$, 
$J\otimes^\nu($the empty knot $[2])$  
and 
$K\otimes^\nu($the empty knot $[2])$, 
differ by one twist-move.
}

\bigbreak\noindent{\bf Note.}  
Of course the $m=0$ case is  Theorem \ref{summerc}.

\bigbreak\noindent{\bf Theorem \ref{aoiro}.}  {\it    
Let $k\in\N$.    
Let $K$ $($resp. $J)$ be $(4k+5)$-dimensional smooth submanifold $\subset S^{4k+7}$. 
Suppose that $K$ and $J$ differ by a single twist-move
and  are nonequivalent. 
Suppose that $K$ is equivalent to 
$A\otimes^{k+1}($the Hopf link$)$ for a 1-knot $A$. 

Then  
there is a unique equivalence class of simple $(4k+1)$-knots for $K$ $($resp. $J)$   
with the following properties. 

\smallbreak
\noindent $\mathrm{(i)}$
There is a representative element $K'$ of the above equivalence class for $K$ such that 
$K$ is equivalent to $K'\otimes($the Hopf link$)$.   

\noindent $\mathrm{(ii)}$
There is a representative element $J'$ of the above equivalence class for $J$ such that 
$J$ is equivalent to $J'\otimes($the Hopf link$)$.   

\noindent $\mathrm{(iii)}$
$K'$ and $J'$ differ 
by a single twist-move 
and are nonequivalent. 
}  

\bigbreak
\noindent{\bf Note. } 
Let $K$ be an $n$-dimensional spherical knot $\subset S^{n+2}$. 
If $\pi_1(S^{n+2}-K)=\Z$ and  
if $\pi_i(S^{n+2}-K)=0  \quad(2\leqq i<\frac{n}{2}, i\in\N)$, 
then we call $K$ a {\it simple knot}.  
See \cite{Levinesimp}. 
\smallbreak

\noindent{\bf Note.}  
There are countably infinitely many $(2p+5)$-dimensional spherical knots ($p\in\N\cup\{0\}$)
which are not the product of any $(2p+1)$-knot and the Hopf link 
by \cite{Kauffman, KauffmanNeumann}.

\smallbreak
\noindent{\bf Note.}   
If $k=0$, we have a different situation:   
There are nonequivalent 1-knots $K'$ and $J'$ with the following properties. 

\smallbreak
\noindent (1)  
$K'$ and $J'$ differ by two crossing-changes not by a crossing-change.  

\noindent (2) 
$K'\otimes^\mu$(the Hopf link) and 
$J'\otimes^\mu$(the Hopf link) differ by a twist-move 
and are nonequivalent for $\mu\in\N$.  
(As we stated before, the twist-move on 1-links is the crossing-change on 1-links.) 
\smallbreak

\noindent 
See Theorem \ref{Sunday}.

\bigbreak
We also prove relations between the pass-move on 1-links and 
the $(p, q)$-pass-move on high dimensional knots.

 \bigbreak\noindent{\bf Theorem \ref{Waltham}.} {\it      
Suppose that two 1-knots $J$ and $K$ differ by one pass-move.
Then 
the $(4\mu+1)$-submanifolds $\subset S^{4\mu+3}$, 
$J\otimes^\mu($the Hopf link$)$ 
and 
$K\otimes^\mu($the Hopf link$)$, 
differ by one $(2\mu+1,2\mu+1)$-pass-move $(\mu\in\N\cup\{0\})$.
}

\bigbreak\noindent{\bf Theorem \ref{mountain}.}  {\it 
Let $J, K$ be simple $(2l+1)$-knots, where $l\in\N$.  
Suppose that $J$ and $K$ differ by one $(l+1, l+1)$-pass-move.
Then
the $(2l+4\mu+1)$-submanifolds \newline 
$\subset S^{2l+4\mu+3}$, 
$J\otimes^\mu($the Hopf link$)$    
and
$K\otimes^\mu($the Hopf link$)$,   
differ by one \newline
$(l+2\mu+1, l+2\mu+1)$-pass-move.
}



\bigbreak\noindent{\bf Theorem \ref{ki}.} {\it    
Let $\mu\in\N$. 
Let $K$ $($resp. $J)$ be a $(4\mu+1)$-submanifold $\subset S^{4\mu+3}$. 
Let $K$ and $J$ be $(2\mu+1, 2\mu+1)$-pass-move-equivalent. 
Suppose that $K$ is equivalent to $K'\otimes^\mu($the Hopf link$)$ for a 1-knot $K'$.
Then  there is a 1-knot $J'$ with the following properties. 

\smallbreak
\noindent $\mathrm{(i)}$
$J$  is equivalent to $J'\otimes^\mu($the Hopf link$)$.   

\noindent $\mathrm{(ii)}$
$K'$ and $J'$ are pass-move-equivalent. 
}

\bigbreak\noindent{\bf Theorem \ref{aka}.}   {\it   
Let $p\in\N$.    
Let $K$ and $J$ be $(2p+5)$-dimensional smooth submanifolds $\subset S^{2p+7}$. 
 Suppose that $K$ and $J$ differ 
by a single $(p+3, p+3)$-pass-move 
and 
are nonequivalent. 
Suppose that $K$ is equivalent to \newline 
$
\begin{cases}
\text{$A\otimes^{\frac{p}{2}+1}($the Hopf link$)$ for a 1-knot $A$}&\text{if $p$ is even}\\
\text{$A\otimes($the Hopf link$)$ for a simple 3-knot $A$}&\text{if $p=1$ $($and hence $2p+5=7)$}\\
\text{$A\otimes^{\frac{p-1}{2}}($the Hopf link$)$ for a simple 7-knot $A$}&\text{if $p$ is odd and $p\neq1$.}
\end{cases}
$

Then  
there is a unique equivalence class of simple $(2p+1)$-knots for $K$ $($resp. $J)$  
with the following properties. 

\smallbreak
\noindent $\mathrm{(i)}$
There is a representative element $K'$ of the above equivalence class for $K$ such that 
$K$ is equivalent to $K'\otimes($the Hopf link$)$.   

\noindent $\mathrm{(ii)}$   
There is a representative element $J'$ of the above equivalence class for $J$ such that 
$J$  is equivalent to $J'\otimes($the Hopf link$)$.   

\noindent $\mathrm{(iii)}$
$K'$ and $J'$ differ by a single $(p+1,p+1)$-pass-move 
and are nonequivalent. 
}

\bigbreak
Next we discuss a relation between  
polynomial invariants of 1-links and those of high dimensional knots 
related by knot products. 

\smallbreak
Suppose that 1-links $K_+, K_-, K_0$ differ only in a 3-ball $B$ as shown below.

\bigbreak
\hskip1cm 
\unitlength 0.1in
\begin{picture}(46.22,12.10)(8.10,-20.20)
%
\special{pn 8}%
\special{ar 1356 1356 546 546  0.8902751 6.2831853}%
\special{ar 1356 1356 546 546  0.0000000 0.8502422}%
%
\special{pn 8}%
\special{pa 946 1716}%
\special{pa 1736 986}%
\special{fp}%
\special{sh 1}%
\special{pa 1736 986}%
\special{pa 1673 1017}%
\special{pa 1697 1022}%
\special{pa 1701 1046}%
\special{pa 1736 986}%
\special{fp}%
%
\special{pn 8}%
\special{pa 1266 1276}%
\special{pa 976 976}%
\special{fp}%
\special{sh 1}%
\special{pa 976 976}%
\special{pa 1008 1038}%
\special{pa 1013 1014}%
\special{pa 1037 1010}%
\special{pa 976 976}%
\special{fp}%
%
\special{pn 8}%
\special{pa 1416 1416}%
\special{pa 1726 1746}%
\special{fp}%
%
\special{pn 8}%
\special{ar 3086 1396 546 546  0.8902751 6.2831853}%
\special{ar 3086 1396 546 546  0.0000000 0.8502422}%
%
\special{pn 8}%
\special{pa 2996 1316}%
\special{pa 2706 1016}%
\special{fp}%
\special{sh 1}%
\special{pa 2706 1016}%
\special{pa 2738 1078}%
\special{pa 2743 1054}%
\special{pa 2767 1050}%
\special{pa 2706 1016}%
\special{fp}%
%
\special{pn 8}%
\special{pa 3146 1456}%
\special{pa 3456 1786}%
\special{fp}%
%
\special{pn 8}%
\special{pa 3210 1520}%
\special{pa 2990 1320}%
\special{fp}%
\special{pa 3240 1550}%
\special{pa 3100 1420}%
\special{fp}%
%
\special{pn 8}%
\special{pa 3140 1330}%
\special{pa 3480 1010}%
\special{fp}%
\special{sh 1}%
\special{pa 3480 1010}%
\special{pa 3418 1041}%
\special{pa 3441 1047}%
\special{pa 3445 1070}%
\special{pa 3480 1010}%
\special{fp}%
%
\special{pn 8}%
\special{pa 3010 1470}%
\special{pa 2690 1770}%
\special{fp}%
%
\special{pn 8}%
\special{ar 4886 1396 546 546  0.8902751 6.2831853}%
\special{ar 4886 1396 546 546  0.0000000 0.8502422}%
%
\special{pn 8}%
\special{pa 4570 1830}%
\special{pa 4540 970}%
\special{fp}%
\special{sh 1}%
\special{pa 4540 970}%
\special{pa 4522 1037}%
\special{pa 4542 1023}%
\special{pa 4562 1036}%
\special{pa 4540 970}%
\special{fp}%
%
\special{pn 8}%
\special{pa 5240 1820}%
\special{pa 5230 970}%
\special{fp}%
\special{sh 1}%
\special{pa 5230 970}%
\special{pa 5211 1037}%
\special{pa 5231 1023}%
\special{pa 5251 1036}%
\special{pa 5230 970}%
\special{fp}%
\put(11.2000,-21.5000){\makebox(0,0)[lb]{$K_+$}}%
\put(28.7000,-21.9000){\makebox(0,0)[lb]{$K_-$}}%
\put(47.1000,-21.8000){\makebox(0,0)[lb]{$K_0$}}%
\end{picture}%
\bigbreak

\noindent
Then the ordered set $(K_+, K_-, K_0)$ is called a {\it crossing-change-triple}. 
We also say that the ordered set $(K_+, K_-, K_0)$ is  {\it related by a crossing-change} in $B$. 
 
Let $A(K)$ be the Alexander-Conway polynomial of 1-links $K$. 
It is well-known that 
$$A(K_+)-A(K_-)=(t-1)\cdot A(K_0).$$
Note that there is another kind of setting of the variable. 
Here, we have the following.

\bigbreak\noindent{\bf Theorem \ref{skein1}.}  {\it    
Let $K_+, K_-, K_0$ be as above. 
There is a polynomial
$
\Delta_{2\mu+1}(\>K_*\otimes^\mu(\text{the Hopf link})\>) \in \Q[t, t^{-1}]$ 
whose $\Q[t, t^{-1}]$-balanced class is 
the $(2\mu+1)$-$\Q[t, t^{-1}]$-Alexander polynomial 
$A_{2\mu+1}(\>K_*\otimes^\mu
(\text{the Hopf link})\>) 
$  
$(\mu\in\N\cup\{0\}. \quad  *=+,-,0.)$   
such that 
\smallbreak

\hskip21mm
$
\Delta_{2\mu+1}(\>K_+\otimes^\mu
(\text{the Hopf link})\>) 
-
\Delta_{2\mu+1}(\>K_-\otimes^\mu
(\text{the Hopf link})\>) $

\hskip17mm
$=
(t-1)\cdot\Delta_{2\mu+1}(\>K_0\otimes^\mu
(\text{the Hopf link})\>).$
}

\bigbreak
\noindent{\bf Note.} 
We review the $p$-$\Q[t, t^{-1}]$-Alexander polynomial $A_p$ for 
$n$-dimensional closed oriented submanifolds $\subset S^{n+2}$ in \S\ref{Alex}.

\vskip3mm
The above  Theorem \ref{skein1}  
follows from 
the following  Theorem \ref{skein2}  
by Theorem \ref{Hopf}.  

\bigbreak\noindent{\bf Theorem \ref{skein2}.}  {\it      
Let $K_+, K_-, K_0$ be as above. 
There is a polynomial 
$\Delta_{\nu+1}(\>K_*
\otimes^\nu[2]\>)
\in \Q[t, t^{-1}]$ 
whose $\Q[t, t^{-1}]$-balanced class is 
the $(\nu+1)$-$\Q[t, t^{-1}]$-Alexander polynomial 
$A_{\nu+1}(\>
K_*
\otimes^\nu[2]\>)  \quad
 (\nu\in\N\cup\{0\}, *=+,-,0)$ 
such  that 
$$
\Delta_{\nu+1}(\>K_+\otimes^\nu[2]\>)
-
\Delta_{\nu+1}(\>K_-\otimes^\nu[2]\>)
=
(t+(-1)^{\nu+1})\cdot \Delta_{\nu+1}(\>K_0\otimes^\nu[2]\>),$$

\noindent 
where  
$[2]$ denotes the empty knot $[2]$. 
}

\bigbreak
\noindent{\bf Note.}  
 We will show an example of Theorem \ref{skein2} in \S\ref{thepol}.

\bigbreak  
The above  Theorem \ref{skein1},  
\ref{skein2} 
are related to the following  Theorem \ref{pepper}.  
The `$l=$even' case is proved in \cite{Ogasa09}. 
In this paper we prove the `$l=$odd' case.   

\bigbreak\noindent{\bf Theorem \ref{pepper}.}   {\it  
Let $K_+$ be a $(2l+1)$-dimensional spherical knot $\subset S^{2l+3} (l\in\N\cup\{0\})$. 
Let $K_-$, $K_0$ be $(2l+1)$-dimensional submanifolds $\subset S^{2l+3}$. 
Let $(K_+, K_-, K_0)$ be a twist-move-triple. 

Then there is a polynomial  
$\Delta_{l+1}(K_*)
\in \Q[t, t^{-1}]$ 
whose $\Q[t, t^{-1}]$-balanced class is 
the $(l+1)$-$\Q[t, t^{-1}]$-Alexander polynomial 
$A_{l+1}(K_*)$  
$(*=+,-,0)$  
and that  
$$
\Delta_{l+1}(K_+)
-
\Delta_{l+1}(K_-)
=(t+(-1)^{l+1})\cdot\Delta_{l+1}(K_0).$$
}
\bigbreak

\noindent
{\bf  Note.}  (1)  
We define the twist-move-triple in \S\ref{localn}. 
If $(K_+, K_-, K_0)$ is twist-move-triple, 
then $K_+$ is obtained from $K_-$,  and 
$K_0$ is $[K_+-(\partial h^{p+1} - \partial B^{n+2})]\cup(h^{p+1}\cap\partial B^{n+2})$. 
See the following diagram. 
In the right $B$, we move $K_0\cap B$ by isotopy.

\bigbreak
\includegraphics[width=12cm]{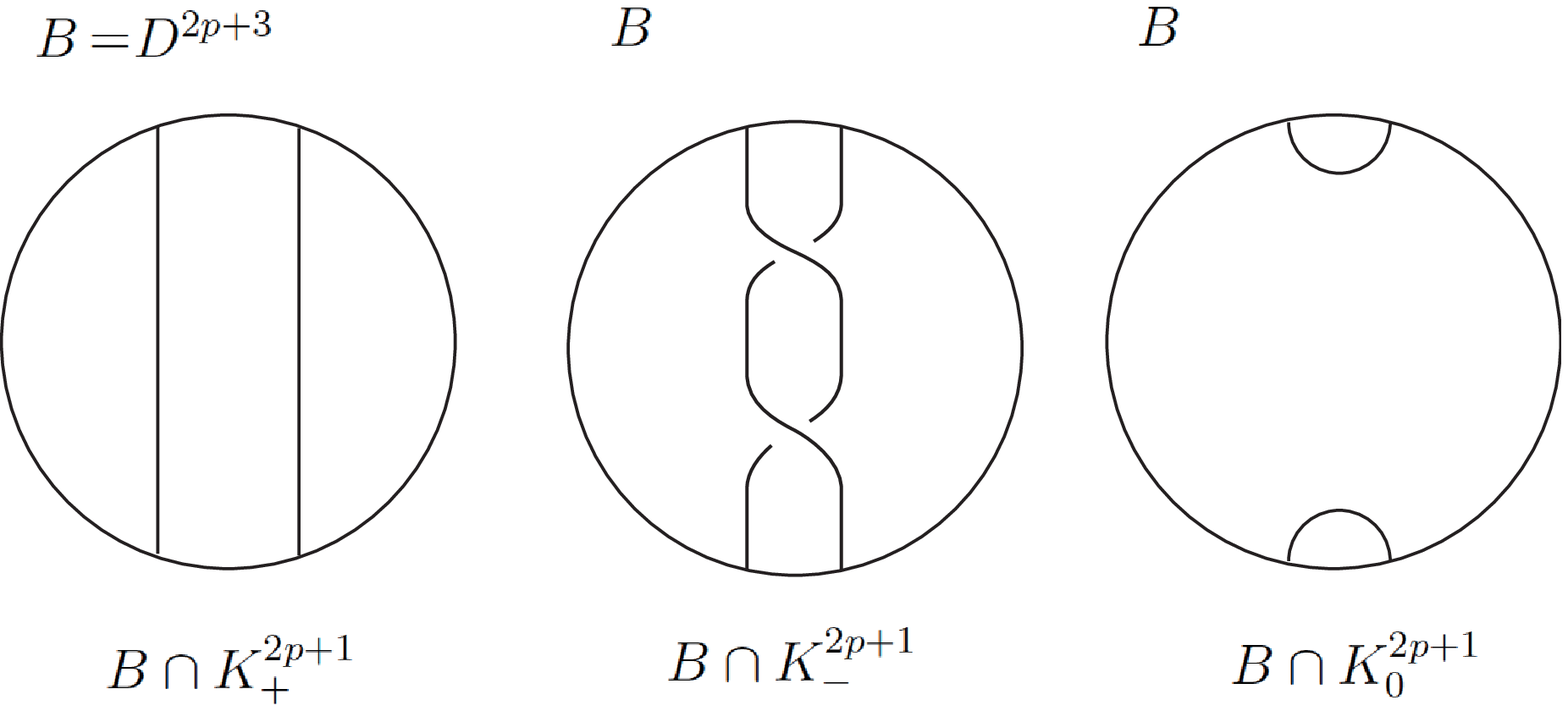}

\bigbreak

\noindent
(2)   We will show an example of Theorem \ref{pepper} in \S\ref{thepol}. 

\noindent (3)   
The identity in Theorem \ref{pepper}  
(resp. Theorem \ref{skein2})   
has a periodicity in dimensions. 
The identity in the `$l=$odd' case of Theorem \ref{pepper}  
has a different form from the identities 
in the `$l=$even' case of Theorem \ref{pepper},  
in Theorem \ref{knots},  in Theorem \ref{middle} and in the well-known case of classical links that is quoted above.

\section{Local moves on  classical links}\label{local1} 

\begin{defn}\label{pass} 
   (\cite{Kauffmanon}.)    
Two 1-links are {\it pass-move-equivalent}
 if one is obtained from the other 
by a sequence of pass-moves. 
See 
the following figure 
 for an illustration of the pass-move. 
Each of four arcs in the 3-ball 
may belong to different components of the 1-link.

\smallbreak
\unitlength 0.1in
\begin{picture}(45.49,19.00)(4.01,-26.71)
%
\special{pn 8}%
\special{ar 1322 1747 921 924  0.0000000 6.2831853}%
%
\special{pn 8}%
\special{pa 607 2328}%
\special{pa 1773 958}%
\special{fp}%
\special{sh 1}%
\special{pa 1773 958}%
\special{pa 1715 996}%
\special{pa 1738 999}%
\special{pa 1745 1022}%
\special{pa 1773 958}%
\special{fp}%
%
\special{pn 8}%
\special{pa 2031 1164}%
\special{pa 864 2541}%
\special{fp}%
\special{sh 1}%
\special{pa 864 2541}%
\special{pa 922 2503}%
\special{pa 898 2500}%
\special{pa 892 2477}%
\special{pa 864 2541}%
\special{fp}%
%
\special{pn 8}%
\special{pa 1515 1695}%
\special{pa 1328 1559}%
\special{fp}%
%
\special{pn 8}%
\special{pa 1315 1941}%
\special{pa 1135 1785}%
\special{fp}%
%
\special{pn 8}%
\special{pa 2868 1882}%
\special{pa 2514 1883}%
\special{fp}%
\special{sh 1}%
\special{pa 2514 1883}%
\special{pa 2581 1903}%
\special{pa 2567 1883}%
\special{pa 2581 1863}%
\special{pa 2514 1883}%
\special{fp}%
%
\special{pn 8}%
\special{pa 2527 1733}%
\special{pa 2901 1733}%
\special{fp}%
\special{sh 1}%
\special{pa 2901 1733}%
\special{pa 2834 1713}%
\special{pa 2848 1733}%
\special{pa 2834 1753}%
\special{pa 2901 1733}%
\special{fp}%
\put(23.2000,-22.3000){\makebox(0,0)[lb]{pass-move}}%
%
\special{pn 8}%
\special{ar 4028 1695 922 924  0.0000000 6.2831853}%
%
\special{pn 8}%
\special{pa 3835 1687}%
\special{pa 3990 1507}%
\special{fp}%
%
\special{pn 8}%
\special{pa 4022 1973}%
\special{pa 3610 2516}%
\special{fp}%
\special{sh 1}%
\special{pa 3610 2516}%
\special{pa 3666 2475}%
\special{pa 3642 2474}%
\special{pa 3634 2451}%
\special{pa 3610 2516}%
\special{fp}%
%
\special{pn 8}%
\special{pa 3777 1785}%
\special{pa 3339 2296}%
\special{fp}%
%
\special{pn 8}%
\special{pa 4074 1417}%
\special{pa 4480 900}%
\special{fp}%
\special{sh 1}%
\special{pa 4480 900}%
\special{pa 4423 940}%
\special{pa 4447 942}%
\special{pa 4455 965}%
\special{pa 4480 900}%
\special{fp}%
%
\special{pn 8}%
\special{pa 4080 1908}%
\special{pa 4261 1688}%
\special{fp}%
%
\special{pn 8}%
\special{pa 4319 1637}%
\special{pa 4731 1113}%
\special{fp}%
%
\special{pn 8}%
\special{pa 1606 1772}%
\special{pa 2134 2179}%
\special{fp}%
\special{sh 1}%
\special{pa 2134 2179}%
\special{pa 2093 2122}%
\special{pa 2092 2146}%
\special{pa 2069 2154}%
\special{pa 2134 2179}%
\special{fp}%
%
\special{pn 8}%
\special{pa 1258 1481}%
\special{pa 736 1054}%
\special{fp}%
%
\special{pn 8}%
\special{pa 1940 2425}%
\special{pa 1380 1999}%
\special{fp}%
%
\special{pn 8}%
\special{pa 1077 1727}%
\special{pa 529 1274}%
\special{fp}%
\special{sh 1}%
\special{pa 529 1274}%
\special{pa 568 1332}%
\special{pa 570 1308}%
\special{pa 593 1301}%
\special{pa 529 1274}%
\special{fp}%
%
\special{pn 8}%
\special{pa 3487 970}%
\special{pa 4847 2135}%
\special{fp}%
\special{sh 1}%
\special{pa 4847 2135}%
\special{pa 4809 2076}%
\special{pa 4806 2100}%
\special{pa 4783 2107}%
\special{pa 4847 2135}%
\special{fp}%
%
\special{pn 8}%
\special{pa 4615 2399}%
\special{pa 3230 1248}%
\special{fp}%
\special{sh 1}%
\special{pa 3230 1248}%
\special{pa 3268 1306}%
\special{pa 3271 1282}%
\special{pa 3294 1275}%
\special{pa 3230 1248}%
\special{fp}%
\end{picture}%
\bigbreak

If $K$ and $J$ are pass-move-equivalent and if $K$ and $K'$ is equivalent, 
then we also say that $K'$ and $J$ are pass-move-equivalent. 
\end{defn}
The following proposition is proved in \cite{Kauffmanon}.

\begin{thm}\label{passth}   
 {\rm (\cite{Kauffmanon}.)}   
Let $L_1$ and $L_2$ be 1-links. 
Then  $L_1$ and $L_2$ are pass-move-equivalent if and only if 
$L_1$ and $L_2$   satisfy 
one of the  following conditions (1) and (2). 

\smallbreak
\noindent 
(1)
 Both $L_1$ and $L_2$ are proper links, and 
$$\mathrm{Arf}(L_1) = \mathrm{Arf}(L_2).$$

\noindent 
(2) 
   Neither $L_1$ nor $L_2$ is a proper link, and 

$$\mathrm{lk}(K_{1j}, L_1-K_{1j})\equiv   
  \mathrm{lk}(K_{2j}, L_2-K_{2j}) 
  \hskip2mm\mathrm{mod 2}\hskip2mm for\hskip2mm  all\hskip2mm j.$$
\end{thm}

\section{Local moves on $n$-knots}\label{localn} 
We review $(p, q)$-pass-moves on $n$-knots ($p, q\in\N$, \quad $p+q=n+1$) 
and twist moves on high dimensional knots. 
\cite{Ogasa98n,  Ogasa04, Ogasa09} defined them.  See also \cite{Ogasa07, OgasaT3, OgasaIH}.  
Confirm that, if $(p, q)=(1, 1)$, 
$(p, q)$-pass-moves are pass-moves on 1-links. 
\bigbreak

We first define $(p, q)$-pass-moves on $n$-knots ($p, q\in\N$, \quad $p+q=n+1$).  
Let $K_+$,  $K_-$, $K_0$ be 
$n$-dimensional closed oriented submanifolds $\subset S^{n+2}$ ($n\in\N$). 
Let $B$ be an $(n+2)$-ball trivially embedded in $S^{n+2}$. 
Suppose that $K_+$ coincides with  $K_-$, $K_0$ in 
$\overline{S^{n+2}-B}$. 

\bigbreak
Take 
an $(n+1)$-dimensional $p$-handle $h^{p}_* (*=+, -)$  
and 
an $(n+1)$-dimensional \newline$(n+1-p)$-handle $h^{n+1-p}$  
in $B$ 
with the following properties. 

\smallbreak
\noindent(1) 
$h^p_* \cap \partial B$  is the attaching part of $h^p_* $.  
$h^{n+1-p}\cap \partial B$  is the attaching part of $h^{n+1-p}$.  

\noindent(2) 
$h^p_* $ (resp. $h^{n+1-p}$) is embedded trivially in $B$. 

\smallbreak

\noindent(3) \hskip4cm$h^{p}_* \cap h^{n+1-p}=\phi$. 

\smallbreak

\noindent(4)
The attaching part of $h^p_+$ coincides with that of $h^p_-$.
The linking number (in $B$) of 

\smallbreak\hskip2cm
 [$h^{p}_+\cup (-h^{p}_-)$]
and 
[$h^{n+1-p}$ whose attaching part is fixed in $\partial B$]
\smallbreak

\hskip2mm is one 
if an orientation is given.  
\smallbreak

Let $K_* (*=+,-)$ satisfy that  
$$K_*\cap \mathrm{Int}B=
(\partial h^{p}_*-\partial B)\cup(\partial h^{n+1-p}-\partial B).$$  

\noindent 
Let 
$$P=K_+ \cap (S^{n+2}-\mathrm{Int}B)$$ 
$$Q=h^{p}_+\cap \partial B$$ 
$$R=h^{n+1-p}\cap \partial B$$
$$T=P\cup Q\cup R.$$  
Then $T$ is 
an $n$-dimensional oriented closed submanifold 
$\subset (S^{n+2}-\mathrm{Int}B) \subset S^{n+2}$.
Let $K_0$ be 
$T$  $\subset S^{n+2}$.
Then we say that  $(K_+$, $K_-$, $K_0)$ is related by a {\it $(p,n+1-p)$-pass-move} in $B$.   
We also say that $(K_+$, $K_-$, $K_0)$ is a {\it $(p,n+1-p)$-pass-move-triple}.   
We say that $K_+$ and $K_-$ differ by one {\it $(p,n+1-p)$-pass-move} in $B$.   
We showed examples of pass-moves on high dimensional knots 
in \S\ref{LMH} and in \S\ref{overview}.

\bigbreak
If  $(K_+$, $K_-$, $K_0)$ is a $(p,n+1-p)$-pass-move-triple, 
then we also say that $(K_-$, $K_+$, $K_0)$ is a {\it $(p,n+1-p)$-pass-move-triple}.   
If $K_+$ and $K_-$ differ by one  $(p,n+1-p)$-pass-move in $B$, 
then we also say that $K_-$ and $K_+$ differ by one  {\it $(p,n+1-p)$-pass-move} in $B$.

\bigbreak
Let $(K_+, K_-, K_0)$ be related by a $(p,n+1-p)$-pass-move in $B$.   
Then there is a Seifert hypersurface $V_*$ for $K_*$ ($*=+.-.0$)  with the following properties. 

\smallbreak
\noindent 
(1) \hskip5cm$V_\sharp=V_0\cup h_\sharp^p\cup h^{n+1-p}  (\sharp=+,-).$
 $$V_\sharp\cap B=h_\sharp^p\cup h^{n+1-p}.$$
 
\noindent 
(2) 
\hskip5cm$V_0\cap \text{Int} B=\phi$. \vskip3mm \hskip3mm
$V_0\cap\partial B$ is the attaching part of $h_\sharp^p\cup h^{n+1-p}$. 

\noindent 
(The idea of the proof is the Thom-Pontrjagin construction.)
\bigbreak

\noindent 
Then the ordered set ($V_+, V_-, V_0$) is called a {\it $(p,n+1-p)$-pass-move-triple of 
Seifert hypersurfaces} for  $(K_+, K_-, K_0)$. 
 We say that 
an ordered set ($V_+, V_-, V_0$) is related by a {\it $(p,n+1-p)$-pass-move} in $B$. 
 We say that  
$V_-$ (resp. $V_+$) is obtained from $V_+$ (resp. $V_-$) 
by a {\it $(p,n+1-p)$-pass-move} in $B$.  

\bigbreak
\noindent 
{\bf Note.} When we construct $K_-$ and $K_0$ from $K_+$, 
we make a change only in $B$ and 
we do not impose any requirement on  
diffeomorphism type or homeomorphism type 
of $K_-$, $K_0$ other than the change only in $B$. 
In this sense, we use the word `local' in the above definition. 

\bigbreak

Figure \ref{localn}.1 is a diagram 
of a $(p, q)$-pass-move. 
Figure \ref{localn}.2, which consists of the three figures (1) (2) (3),       
is a diagram 
of a $(p, q)$-pass-move-triple.

\vskip5mm
\unitlength 0.1in
\begin{picture}(57.35,30.77)(1.70,-38.40)
%
\special{pn 8}%
\special{ar 4968 1942 937 937  0.8932296 6.2831853}%
\special{ar 4968 1942 937 937  0.0000000 0.8617069}%
%
\special{pn 8}%
\special{pa 4236 1373}%
\special{pa 5415 2755}%
\special{fp}%
%
\special{pn 8}%
\special{pa 5709 2500}%
\special{pa 4541 1118}%
\special{fp}%
%
\special{pn 8}%
\special{pa 5415 1118}%
\special{pa 4998 1515}%
\special{fp}%
\special{pa 4612 1922}%
\special{pa 4134 2369}%
\special{fp}%
\special{pa 5700 1382}%
\special{pa 5242 1789}%
\special{fp}%
\special{pa 4835 2216}%
\special{pa 4388 2643}%
\special{fp}%
%
\special{pn 8}%
\special{ar 1817 1980 937 937  5.5898109 6.2831853}%
\special{ar 1817 1980 937 937  0.0000000 5.5597350}%
%
\special{pn 8}%
\special{pa 1258 2722}%
\special{pa 2622 1521}%
\special{fp}%
%
\special{pn 8}%
\special{pa 2365 1230}%
\special{pa 1000 2420}%
\special{fp}%
%
\special{pn 8}%
\special{pa 986 1546}%
\special{pa 1390 1957}%
\special{fp}%
\special{pa 1802 2337}%
\special{pa 2255 2808}%
\special{fp}%
\special{pa 1246 1258}%
\special{pa 1661 1709}%
\special{fp}%
\special{pa 2093 2109}%
\special{pa 2527 2550}%
\special{fp}%
%
\special{pn 8}%
\special{pa 963 3251}%
\special{pa 1135 2336}%
\special{dt 0.045}%
\special{sh 1}%
\special{pa 1135 2336}%
\special{pa 1103 2398}%
\special{pa 1125 2388}%
\special{pa 1142 2405}%
\special{pa 1135 2336}%
\special{fp}%
\special{pa 1004 3230}%
\special{pa 1379 2611}%
\special{dt 0.045}%
\special{sh 1}%
\special{pa 1379 2611}%
\special{pa 1327 2658}%
\special{pa 1351 2657}%
\special{pa 1362 2678}%
\special{pa 1379 2611}%
\special{fp}%
%
\special{pn 8}%
\special{pa 2650 3027}%
\special{pa 2457 2509}%
\special{dt 0.045}%
\special{sh 1}%
\special{pa 2457 2509}%
\special{pa 2462 2578}%
\special{pa 2476 2559}%
\special{pa 2499 2564}%
\special{pa 2457 2509}%
\special{fp}%
\special{pa 2640 3007}%
\special{pa 2040 2539}%
\special{dt 0.045}%
\special{sh 1}%
\special{pa 2040 2539}%
\special{pa 2080 2596}%
\special{pa 2082 2572}%
\special{pa 2105 2564}%
\special{pa 2040 2539}%
\special{fp}%
\put(24.4700,-33.2200){\makebox(0,0)[lb]{$S^{p-1}\x D^{n+1-p}$}}%
\put(3.2200,-34.3300){\makebox(0,0)[lb]{$S^{n-p}\x D^p$}}%
\put(24.0600,-35.8700){\makebox(0,0)[lb]{$=\overline{\partial h^{n+1-p}-\partial B}$}}%
\put(13.0800,-9.7400){\makebox(0,0)[lb]{$B\cap K_+$}}%
\put(45.0000,-9.3300){\makebox(0,0)[lb]{$B\cap K_-$}}%
\put(1.7000,-36.3700){\makebox(0,0)[lb]{$=\overline{\partial h^p_+-\partial B}$}}%
\put(28.1000,-40.1000){\makebox(0,0)[lb]{Figure \ref{localn}.1}}%
\end{picture}%
\vskip3mm

\noindent 
{\bf A $(p, n+1-p)$-pass-move on an $n$-dimensional submanifold  
$\subset S^{n+2}$. 
Note $B=B^{n+2}=D^{n+2}\subset S^{n+2}$.   
This figure includes $h^p_+$ and $h^{n+1-p}$.}

\np
\hskip1mm\vskip1cm

\hskip4cm\includegraphics[width=10cm]{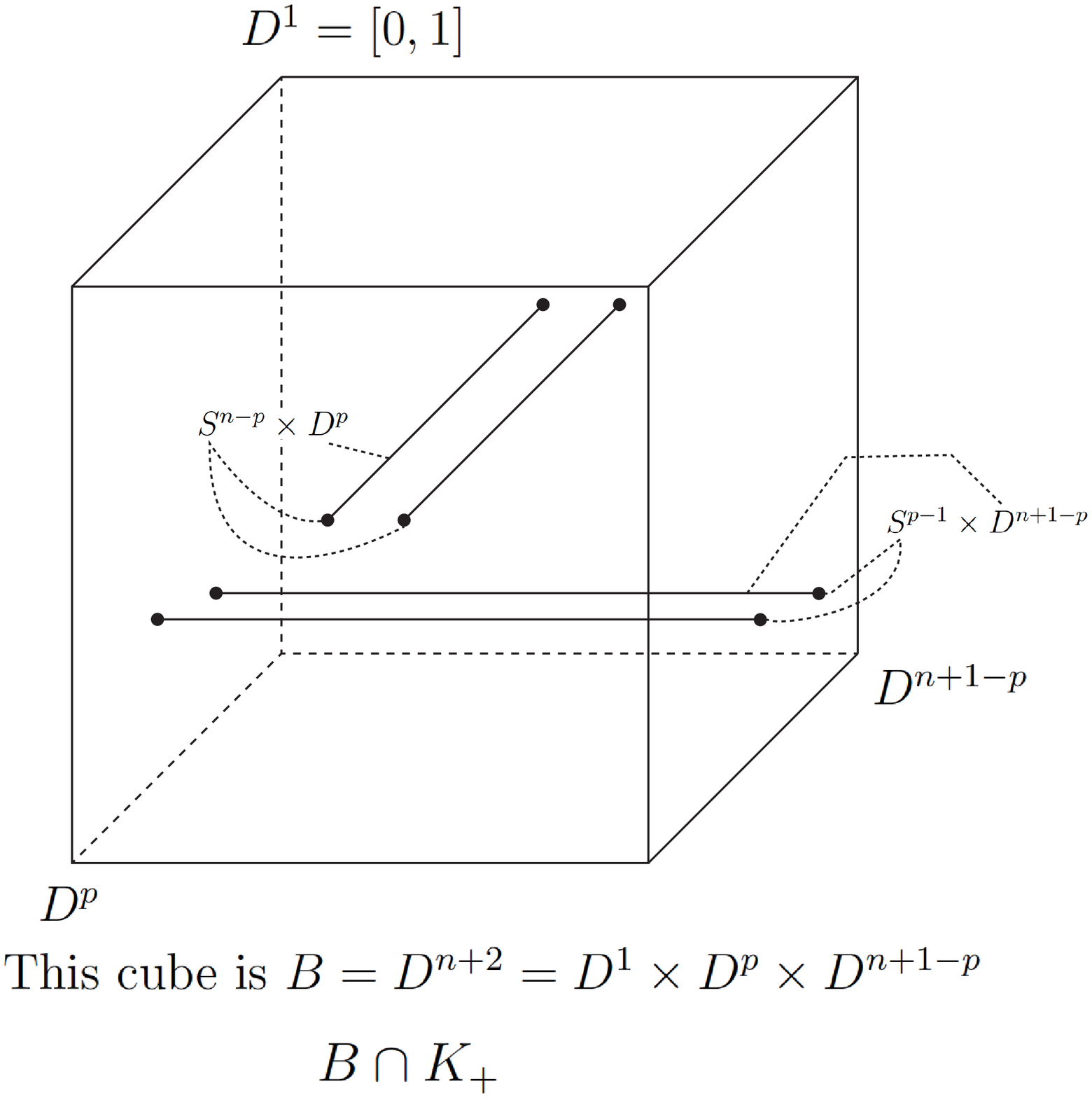}
 
\vskip10mm
\hskip40mm Figure \ref{localn}.2.(1): {\bf A $(p,n+1-p)$-pass-move-triple}

\np
\hskip2cm
\includegraphics[width=10cm]{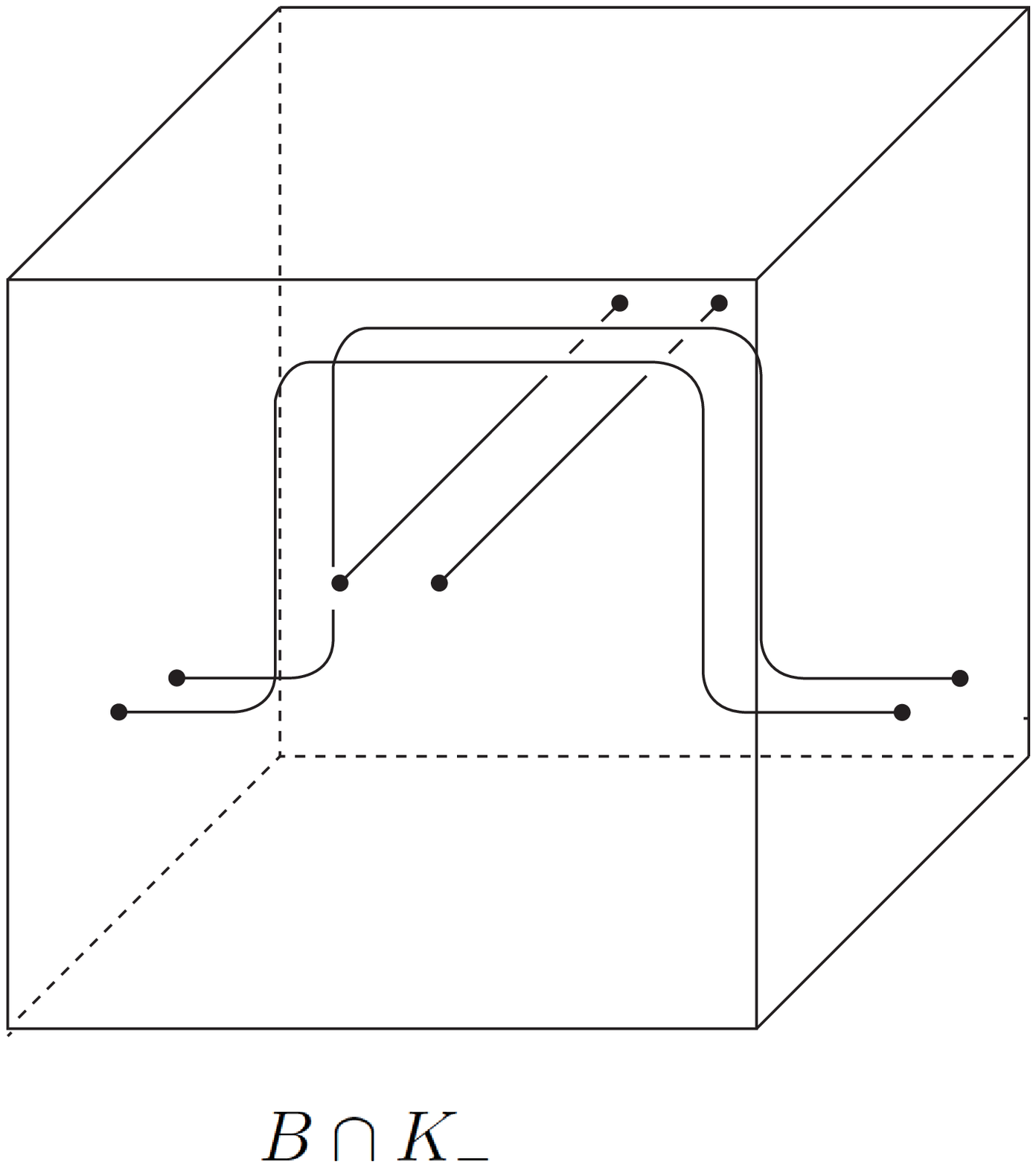}

\vskip-20mm
\hskip17mm Figure \ref{localn}.2.(2): {\bf A $(p,n+1-p)$-pass-move-triple}

\np
\vskip-25mm
\hskip2cm
\includegraphics[width=10cm]{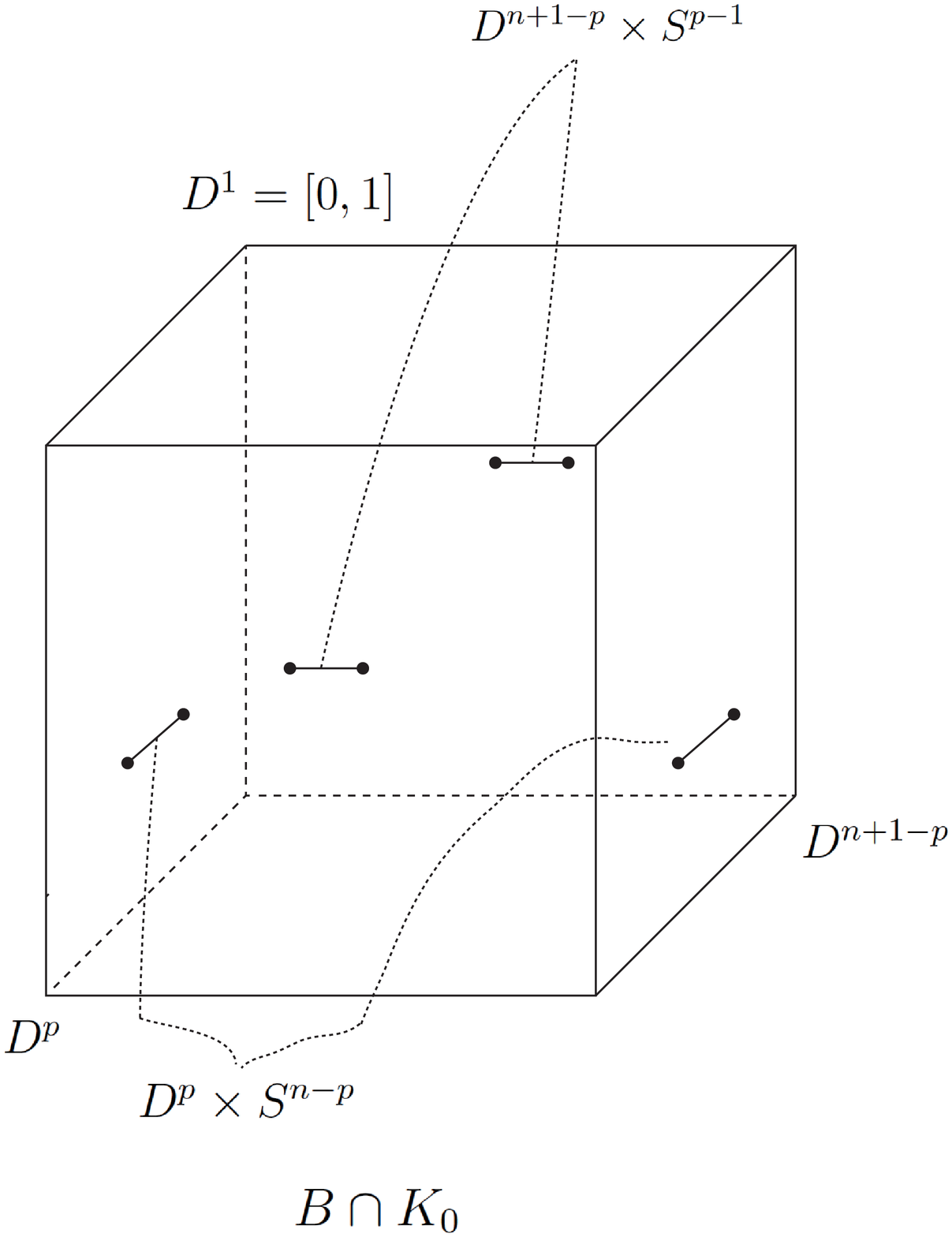}

\vskip-10mm
\hskip17mm Figure \ref{localn}.2.(3): {\bf A $(p,n+1-p)$-pass-move-triple}

\np
\hskip1mm\vskip-20mm
\hskip3cm \input 3+.tex  

\vskip15mm

\hskip3cm \input 3-.tex  



\np

\hskip1mm\vskip-20mm
\hskip3cm\input 30.tex  




\vskip1cm
In Figure \ref{localn}.3, 
which consists of the three figures (1) (2) (3), 
 we draw a (1,2)-pass-move-triple. 
(the $p=1$ and $n=2$ case).  
%
%
Since $(K_+, K_-, K_0)$ is related by a $(1,2)$-pass-move in $B$, 
$B$ has the following properties. 
We regard $B$ as (2-disc)$\x[0,1]\x\{t\vert -1\leqq t\leqq1\}$.  

\smallbreak
\noindent 
(i) $K_+-B$, $K_--B$, and $K_0-B$ coincide each other.   

\noindent 
(ii) $B\cap K_+$,  $B\cap K_-$, $B\cap K_0$ 
are shown as above. 

\bigbreak
In the above figures 
we draw 
$B_{-0.5}\cap K_*$, $B_0\cap K_*$, $B_{0.5}\cap K_*$, 
where $B_{t_0}$=(2-disc)$\x[0,1]\x\{t\vert t=t_0\}$.  
We suppose that 
each vector $\overrightarrow{x}$, $\overrightarrow{y}$ 
in the above figures
is a tangent vector of each disc at a point. 
(Note that we use $\overrightarrow{x}$ (resp. $\overrightarrow{y}$)
for different vectors.)
The orientation of each disc 
in
the above figures 
is determined by the each set 
$\{\overrightarrow{x},\overrightarrow{y}\}$. 
Around Figure 4.1 and 4.2 in \cite{Ogasa04}, 
we wrote more explanation of the figure of $B\cap K_+$ and that of $B\cap K_-$.

In \cite{Ogasa04} one more local move was discussed, 
called the `ribbon-move'.  
In  \cite{Ogasa04} are proved the following results. 
Let $K$ and $K'$ be 
2-dimensional closed oriented submanifolds 
$\subset S^4$.  
The following conditions (1) and (2) are equivalent.

\smallbreak
\noindent 
(1)  
$K$ is  (1,2)-pass-move-equivalent to $K'$.

\noindent 
(2)   
$K$ is  ribbon-move-equivalent to $K'$. 

Furthermore if $K$ is obtained from $K'$ by one ribbon-move, 
then $K$ is obtained from $K'$ by one (1,2)-pass-move.
\smallbreak

\np
We next define twist-moves on high dimensional knots.  
Let $K_+$,  $K_-$, $K_0$ be 
$(2p+1)$-dimensional closed oriented submanifold $\subset S^{2p+3}$ ($p\in\N\cup\{0\}$). 
Let $B$ be a $(2p+3)$-ball trivially embedded in $S^{2p+3}$. 
Suppose that $K_+$ coincides with  $K_-$, $K_0$ in $\overline{S^{2p+3}-B}$. 
Take 
a single $(2p+2)$-dimensional $(p+1)$-handle $h_+$ (resp. $h_-$) embedded in $B$
such that \newline
 [the handle]$\cap\partial B$ 
is the attaching part of the handle. 
Note. \cite{Haefligerunknot, Haefligerknot, Whitney, Whitneytrick} etc. 
imply that the core of $h_+$ (resp. $h_-$) is trivially embedded in $B$ under the above condition.     
%
%
Suppose that $(h_+-$its attaching part)$\cap(h_--$its attaching part)$=\phi$.  
Suppose that their attaching parts coincide.  
Thus we can suppose that  
we regard $h_+\cup h_-$ as an oriented $(2p+2)$-submanifold $\subset S^{2p+1}$
if we give the opposite orientation to $h_-$.  
Then we can define a $(p+1)$-Seifert matrix for the $(2p+2)$-submanifold $h_+\cup h_-$.  
We can suppose that the Seifert matrix is the matrix $(1)$.


\bigbreak
Let $K_* (*=+,-)$ satisfy that  
$K_*\cap \mathrm{Int}B=$ 
$(\partial h_*-\partial B)$.  
Note the following. When we define $K_+$, $h_+$exists in $B$ and $h_-$ does not exist in $B$. 
When we define $K_-$, $h_-$exists in $B$ and $h_+$ does not exist in $B$. 
%
%
Let $P=K_+ \cap (S^{2p+3}-\mathrm{Int}B)$. 
Let $Q=h_+\cap \partial B$. 
Let  
$T=P\cup Q$. 
Then $T$ is 
an $(2p+1)$-dimensional oriented closed submanifold in $S^{2p+3}-\mathrm{Int}B$.
Let $K_0$ be 
$T$ in $S^{2p+3}$.
Then we say that  an ordered set $(K_+$, $K_-$, $K_0)$ is related by a {\it twist-move.}     
$(K_+$, $K_-$, $K_0)$ is called a {\it twist-move-triple}.
We say that $K_+$ and $K_-$ differ by one {\it twist-move} in $B$. 
If $(K_+$, $K_-$, $K_0)$ is a twist-move-triple, 
then we also say that $(K_-$, $K_+$, $K_0)$ is a {\it twist-move-triple}.
If $K_+$ and $K_-$ differ by one twist-move in $B$,  
we also say that $K_-$ and $K_+$ differ by one {\it twist-move} in $B$. 
We will show an example of twist-moves on high dimensional knots 
after Theorem \ref{Tokyo} in \S\ref{thecr}. 
We will show examples of twist-move-triple of high dimensional knots  
after Theorem \ref{skein2} in \S\ref{thepol}, and 
after Theorem \ref{pepper} in \S\ref{thepol}.

\bigbreak
\noindent{\bf{Note.}} 
 The $XXII$-move in \cite{Ogasa09} is the twist-move in the `$p$=even' case.

\bigbreak
\noindent{\bf{Note.}} 
Suppose that $p$ is an odd natural number, put $p=2k+1$. 
The twist-move for $(4k+3)$-submanifolds $\subset S^{4k+5}$ 
($4k+3\in\N$, \quad $k\in\N\cup\{0\}$)   
has the following property: 
Suppose that $K_+$ is made into $K_-$  by the twist-move. 
Then $K_-$ is a nonspherical knot in general even if $K_+$ is a spherical knot. 
Furthermore the   $H_*(K_-;\Z)$ is not congruent to $H_*(K_+;\Z)$ in general.   
Example: A Seifert hypersurface $V_*$ for  a 3-knot $K_*$ ($*=+,-$); 
Framed link representation of $V_+$  is the Hopf link 
such that the framing of one component is zero 
and that that of the other is two. 
Framed link representation of $V_-$  is the Hopf link 
such that the framing of each component is two.

\bigbreak
Let $(K_+, K_-, K_0)$ be related by a twist-move in $B$.   
Then there is a Seifert hypersurface $V_*$ for $K_*$ ($*=+.-.0$)  with the following properties. 

\smallbreak
\noindent 
(1) $V_\sharp=V_0\cup h_\sharp$  ($\sharp=+,-$).
 $V_\sharp\cap B=h_\sharp$.  

\smallbreak
 \noindent 
(2) \hskip5cm$V_0\cap$ Int $B=\phi$. 
\smallbreak

\hskip5mm$V_0\cap\partial B$ is the attaching part of $h_\sharp^p$. 

\noindent 
(The idea of the proof is the Thom-Pontrjagin construction.)

\bigbreak
 The ordered set ($V_+, V_-, V_0$) is called a {\it twist-move-triple of 
Seifert hypersurfaces} for  $(K_+, K_-, K_0)$. 
We say that 
$V_-$ (resp. $V_+$) is obtained from $V_+$ (resp. $V_-$)   
by a {\it twist-move} in $B$.

\bigbreak
Figure \ref{localn}.4, which consists of the three figures (1) (2) (3),  
is a diagram 
of a twist-move-triple. 
The upper half 
of Figure \ref{localn}.5 is another diagram 
of a twist-move triple.
Compare 
the upper half 
of Figure \ref{localn}.5 
and 
the lower half. 
If $p=0$ (hence $n=2p+1=1$ ), 
the left figure in the upper half
and 
that in the lower half are same. 
That is, if $p=0$ (hence $n=2p+1=1$ ), 
a twist-move-triple is a crossing-change-triple of 1-links.  
Note that we move $B\cap K_0$ by isotopy in the right $B$ in the upper half 
of Figure \ref{localn}.5.


\vskip5mm
\includegraphics[width=10cm]{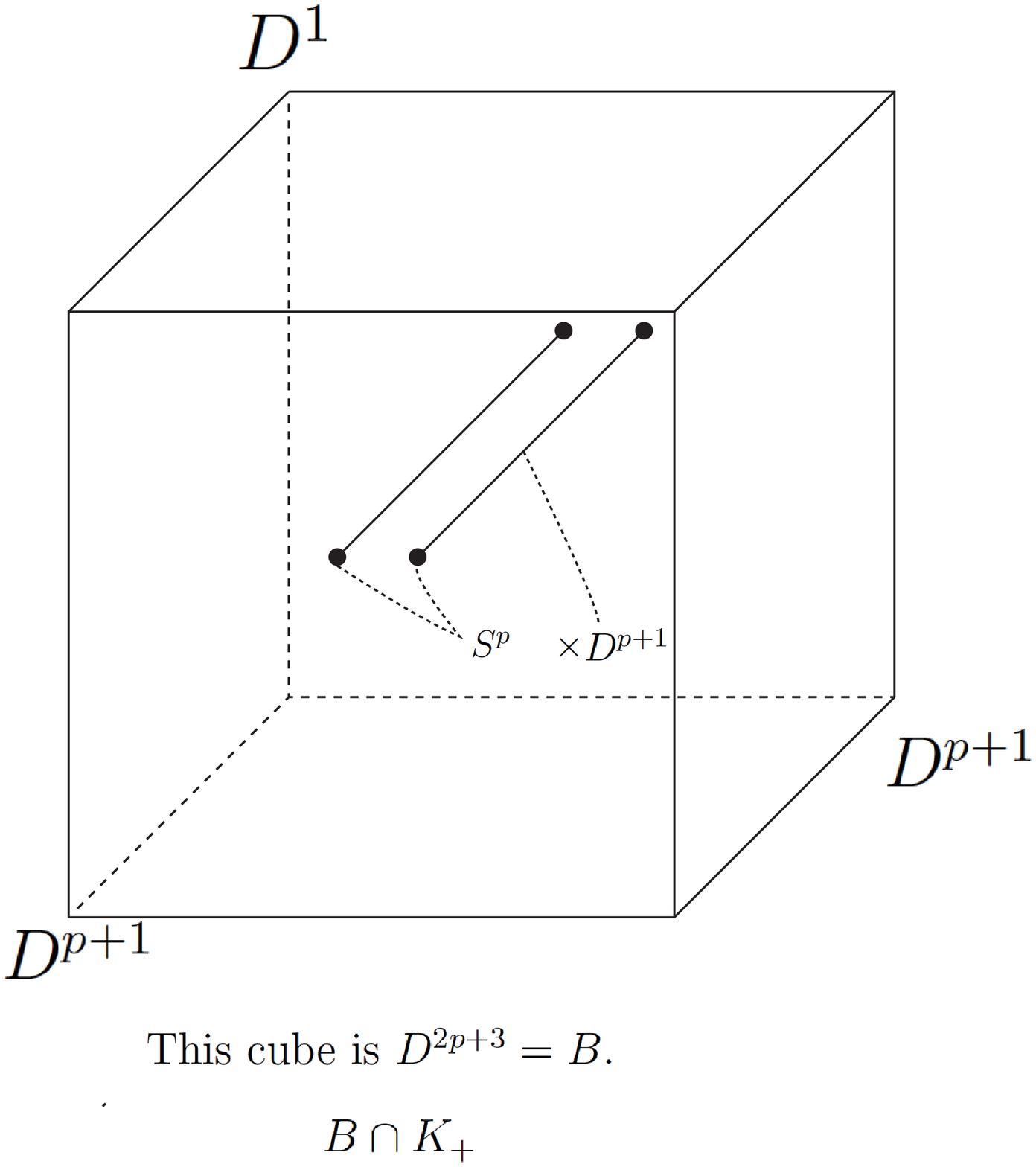}

\vskip10mm
\hskip-5mm Figure \ref{localn}.4.(1): {\bf A twist-move-triple}

\np
\includegraphics[width=10cm]{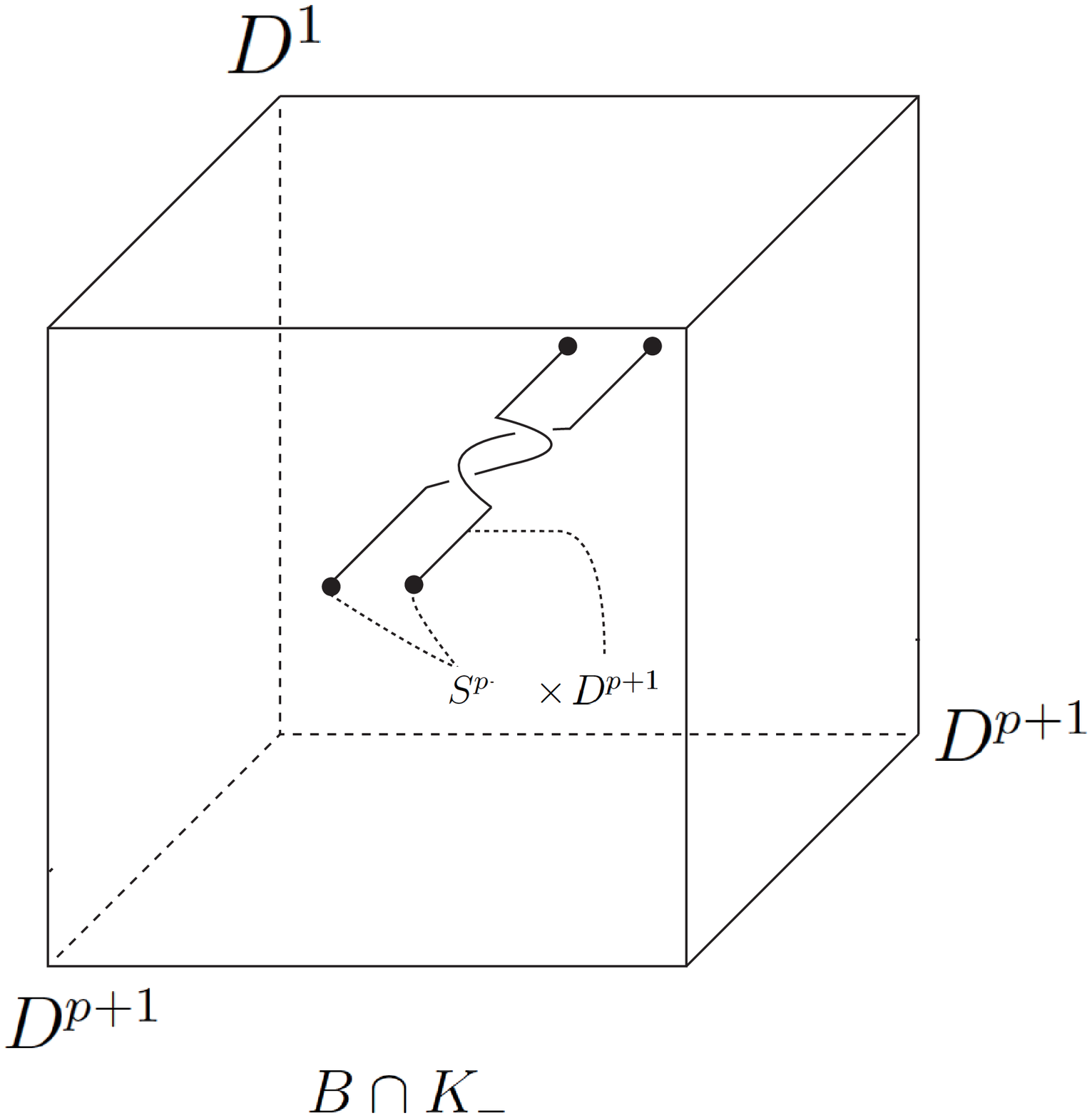}

\hskip-5mm Figure \ref{localn}.4.(2): {\bf A twist-move-triple}

\np
\includegraphics[width=10cm]{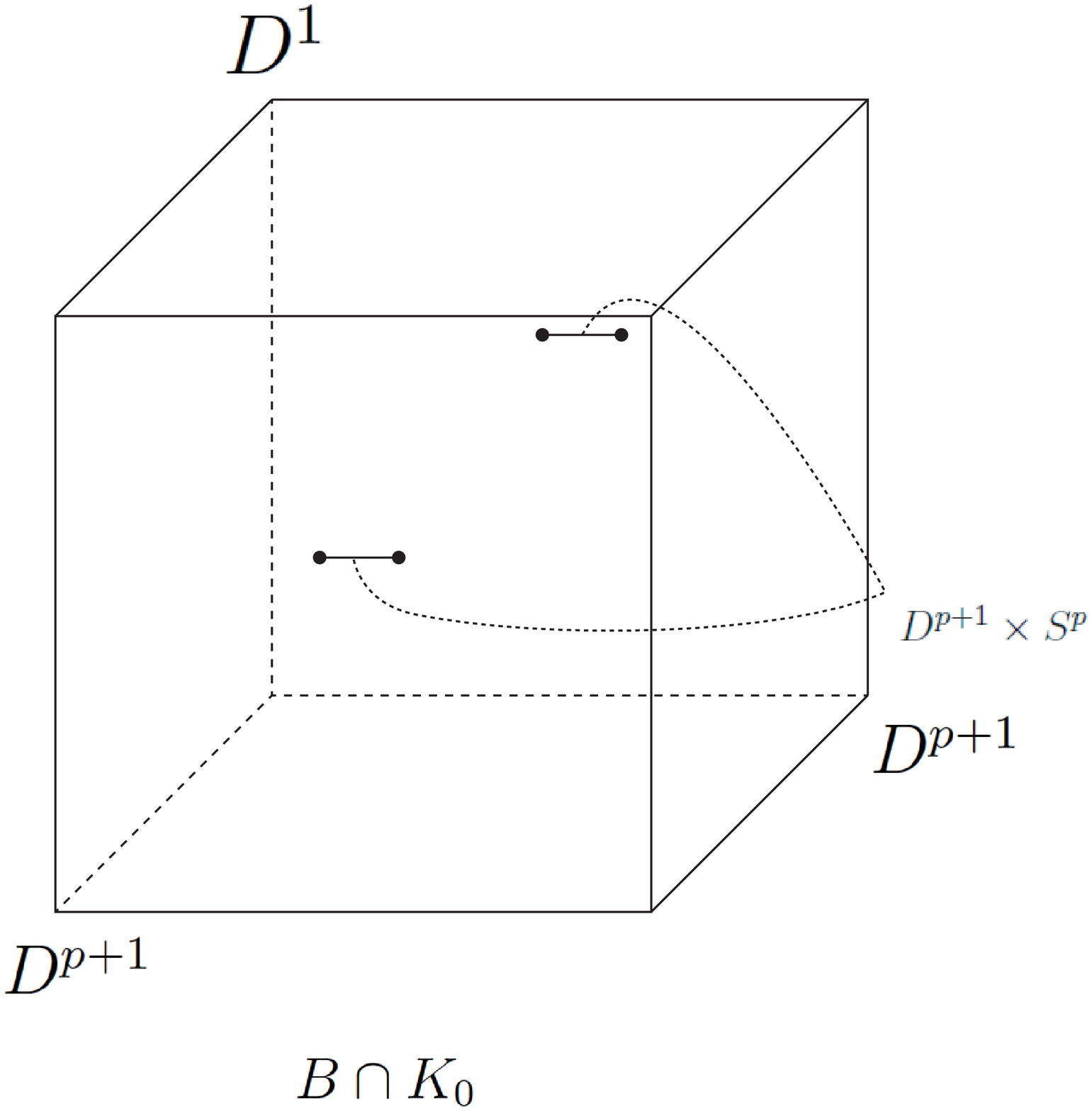}

\vskip10mm
Figure \ref{localn}.4.(3): {\bf A twist-move-triple}

\np

\includegraphics[width=12cm]{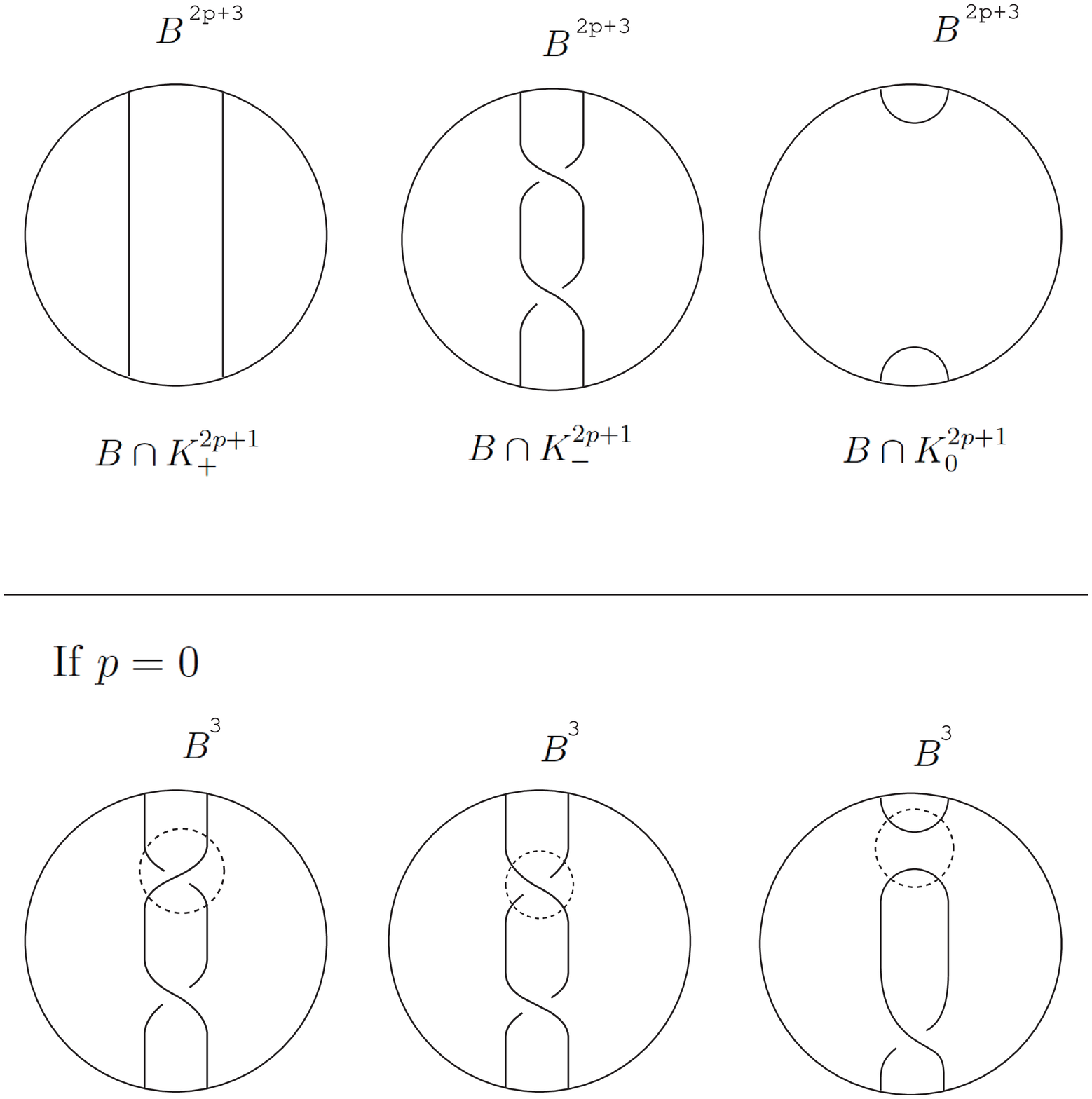}

\vskip2mm
 The triple of three  
\quad
\unitlength 0.1in
\begin{picture}(2.40,2.40)(0.60,-2.60)
%
\special{pn 8}%
\special{ar 180 140 120 120  0.0000000 0.1000000}%
\special{ar 180 140 120 120  0.4000000 0.5000000}%
\special{ar 180 140 120 120  0.8000000 0.9000000}%
\special{ar 180 140 120 120  1.2000000 1.3000000}%
\special{ar 180 140 120 120  1.6000000 1.7000000}%
\special{ar 180 140 120 120  2.0000000 2.1000000}%
\special{ar 180 140 120 120  2.4000000 2.5000000}%
\special{ar 180 140 120 120  2.8000000 2.9000000}%
\special{ar 180 140 120 120  3.2000000 3.3000000}%
\special{ar 180 140 120 120  3.6000000 3.7000000}%
\special{ar 180 140 120 120  4.0000000 4.1000000}%
\special{ar 180 140 120 120  4.4000000 4.5000000}%
\special{ar 180 140 120 120  4.8000000 4.9000000}%
\special{ar 180 140 120 120  5.2000000 5.3000000}%
\special{ar 180 140 120 120  5.6000000 5.7000000}%
\special{ar 180 140 120 120  6.0000000 6.1000000}%
\end{picture}%
 \quad 
makes a crossing-change-triple of a 1-dimensional link.

\vskip1cm
\hskip50mm Figure \ref{localn}.5 

{\bf A twist-move-triple of 1-links is a crossing-change-triple of 1-links.}

\np
\section{Products of knots}\label{knotproduct}   
We reviewed the knot product construction in \S\ref{Introduction}. 
%
%
In this section we state some remarks. 


Let 
$A\otimes^\mu B$ 
mean  
$A\otimes B,...,\otimes B$, which is composed of 
one copy of $A$ and $\mu$ copies of $B$, 
where $\mu\in\N\cup\{0\}$. 
Let $\otimes^\mu B$ 
mean  
$B\otimes,...,\otimes B$, which is composed of $\mu$ copies of $B$, 
where $\mu\in\N\cup\{0\}$.

\begin{defn}\label{empty}  
    {\rm (\cite{Kauffman,  KauffmanNeumann}.)}    
Let $n\in\N$. 
The {\it empty knot $[n]$} is a smooth map 
$S^1\to S^1$ such that 
$\theta\mapsto n\theta$, 
where $S^1=\{e^{2\pi i \theta}|\theta\in\R\}$
\end{defn}

We regard a Seifert hypersurface of the empty knot $[n]$
as a set of $n$ points $\subset S^1$. 
We can regard the empty knot $[n]$ as a fibred knot. 
In \cite{Kauffman,  KauffmanNeumann}  is defined a knot product of 
the empty knot and an $n$-dimensional closed oriented submanifold $\subset S^{n+2}$.

The {\it $($positive$)$ Hopf link} or 
the {\it $($linking number$)$ $(+1)$ Hopf link} is shown as follows.

\bigbreak
\unitlength 0.1in
\begin{picture}(11.06,11.06)(6.21,-15.52)
%
\special{pn 8}%
\special{ar 1174 999 553 553  0.2449787 6.2831853}%
\special{ar 1174 999 553 553  0.0000000 0.2204330}%
%
\special{pn 8}%
\special{ar 1870 978 552 552  0.2454152 6.2831853}%
\special{ar 1870 978 552 552  0.0000000 0.2208288}%
%
\special{pn 8}%
\special{sh 0}%
\special{pa 1458 531}%
\special{pa 1572 531}%
\special{pa 1572 637}%
\special{pa 1458 637}%
\special{pa 1458 531}%
\special{ip}%
%
\special{pn 8}%
\special{sh 0}%
\special{pa 1479 1368}%
\special{pa 1593 1368}%
\special{pa 1593 1475}%
\special{pa 1479 1475}%
\special{pa 1479 1368}%
\special{ip}%
%
\special{pn 8}%
\special{pa 1451 616}%
\special{pa 1470 590}%
\special{pa 1494 569}%
\special{pa 1517 549}%
\special{pa 1539 532}%
\special{pa 1543 524}%
\special{sp}%
%
\special{pn 8}%
\special{pa 1593 1354}%
\special{pa 1585 1384}%
\special{pa 1557 1394}%
\special{pa 1531 1409}%
\special{pa 1503 1422}%
\special{pa 1486 1449}%
\special{pa 1472 1461}%
\special{sp}%
%
\special{pn 8}%
\special{pa 1664 751}%
\special{pa 1650 857}%
\special{fp}%
%
\special{pn 8}%
\special{pa 1671 751}%
\special{pa 1763 808}%
\special{fp}%
%
\special{pn 8}%
\special{pa 1785 431}%
\special{pa 1678 410}%
\special{fp}%
\special{pa 1770 424}%
\special{pa 1721 516}%
\special{fp}%
\end{picture}%
\bigbreak

The {\it negative Hopf link} or 
the {\it $($linking number$)$ $(-1)$ Hopf link} is shown as follows.

\bigbreak
\unitlength 0.1in
\begin{picture}(11.06,11.06)(6.21,-15.29)
%
\special{pn 8}%
\special{ar 1174 976 553 553  0.2449787 6.2831853}%
\special{ar 1174 976 553 553  0.0000000 0.2204330}%
%
\special{pn 8}%
\special{ar 1870 955 552 552  0.2454152 6.2831853}%
\special{ar 1870 955 552 552  0.0000000 0.2190618}%
%
\special{pn 8}%
\special{sh 0}%
\special{pa 1458 507}%
\special{pa 1572 507}%
\special{pa 1572 614}%
\special{pa 1458 614}%
\special{pa 1458 507}%
\special{ip}%
%
\special{pn 8}%
\special{sh 0}%
\special{pa 1479 1345}%
\special{pa 1593 1345}%
\special{pa 1593 1452}%
\special{pa 1479 1452}%
\special{pa 1479 1345}%
\special{ip}%
%
\special{pn 8}%
\special{pa 1451 593}%
\special{pa 1470 567}%
\special{pa 1494 546}%
\special{pa 1517 526}%
\special{pa 1538 507}%
\special{pa 1543 500}%
\special{sp}%
%
\special{pn 8}%
\special{pa 1593 1331}%
\special{pa 1585 1361}%
\special{pa 1557 1372}%
\special{pa 1530 1386}%
\special{pa 1503 1400}%
\special{pa 1487 1427}%
\special{pa 1472 1438}%
\special{sp}%
%
\special{pn 8}%
\special{pa 1664 728}%
\special{pa 1650 834}%
\special{fp}%
%
\special{pn 8}%
\special{pa 1671 728}%
\special{pa 1763 784}%
\special{fp}%
%
\special{pn 8}%
\special{pa 1756 401}%
\special{pa 1778 380}%
\special{fp}%
%
\special{pn 8}%
\special{pa 1770 415}%
\special{pa 1877 465}%
\special{fp}%
%
\special{pn 8}%
\special{pa 1763 394}%
\special{pa 1834 330}%
\special{fp}%
%
\special{pn 8}%
\special{pa 1763 408}%
\special{pa 1813 344}%
\special{fp}%
\end{picture}%
\bigbreak

\begin{thm}\label{Hopf}  
  {\rm (\cite{Kauffman, KauffmanNeumann}.)}   
Let $[n]$ denote the empty knot of degree $n$. 
Then we have 
$$[2]\otimes[2]=\text{the negative Hopf link.} $$

\noindent 
For $\mu\in\N$, we have 

\smallbreak\hskip4cm
$\begin{matrix}
\text{{\hskip5mm\tiny $2\mu$}}\\
\text{\large $\otimes$}\\
\text{  }\\
\end{matrix}$
\hskip-2mm$[2]=$  
\hskip-2mm
$\begin{matrix}
\text{{\hskip5mm\tiny $\mu$}}\\
\text{\large $\otimes$}\\
\text{  }\\
\end{matrix}$
\hskip-1mm$($the negative Hopf link$).$  
\smallbreak

\noindent 
For any $n$-dimensional closed oriented submanifold $K\subset S^{n+2}$, 

\smallbreak\hskip4cm
$K$ 
\hskip-3mm
$\begin{matrix}
\text{{\hskip5mm\tiny $2\mu$}}\\
\text{\large $\otimes$}\\
\text{  }\\
\end{matrix}$
\hskip-2mm$[2]=$  
$K$ 
\hskip-2mm
$\begin{matrix}
\text{{\hskip5mm\tiny $\mu$}}\\
\text{\large $\otimes$}\\
\text{  }\\
\end{matrix}$
\hskip-1mm
$($the negative Hopf link$)$.   
\end{thm}

\noindent
{\bf Note.} See line($-12$) of page 389 
and line 18 of page 391 
of \cite{KauffmanNeumann}

\bigbreak
In this paper the Hopf link means the negative Hopf link.


\section{Review of the $\Q[t, t^{-1}]$-Alexander polynomials for $n$-knots and 
$n$-dimensional closed oriented submanifolds
} 
\label{Alex}     
  
We review 
the $\Q[t, t^{-1}]$-Alexander polynomials for $n$-knots and $n$-links and $n$-dimensional closed oriented submanifolds, 
Seifert matrices, 
Alexander matrices, etc.   
See 
\cite{Alexander, Levinepol, Levinecob, Levinesimp}. 

\bigbreak
Let $K=(K_1,...,K_\xi)$ be an $n$-dimensional closed oriented submanifold of 
$S^{n+2}$ ($n$\newline$
\in\N$). 
Let each $K_i$ be connected. 
It is known 
that 
any tubular neighborhood of $K$ is 
diffeomorphic to 
$K\x D^2$ (see pages 49, 50 of \cite{Kirby}). 
Let $X=\overline{S^{n+2}-K\x D^2}$. 
By using the orientation of $S^{n+2}$ and that of $K$,  
we can determine 
an orientation of $\partial D^2$. Take a homomorphism
$\alpha:H_1(X;\Z)\rightarrow\Z$ to carry 
%
all $[\partial D^2]$ with the orientations to $+1$.
Take the infinite cyclic covering 
$\pi:$ $\widetilde X\rightarrow X$ 
associated with $\alpha$. 
$\widetilde X$ is called 
the {\it canonical cyclic covering space} of $K$. 
We can regard $H_p(\widetilde X;\Z)$ as a $\Z[t, t^{-1}]$-module 
by using 
the covering translation 
$\widetilde X\rightarrow\widetilde X$ 
defined by $\alpha$. 
It is called the {\it $\Z[t, t^{-1}]$-$p$-Alexander module}.  
We can also regard $H_p(\widetilde X;\Q)$ as a $\Q[t, t^{-1}]$-module.  
It is called the {\it $\Q[t, t^{-1}]$-$p$-Alexander module}.

\bigbreak
According to module theory, 
it holds that any $\Q[t, t^{-1}]$-module  is congruent to 
$$(\Q[t, t^{-1}]/{\lambda_1})\oplus\cdot\cdot\cdot\oplus
(\Q[t, t^{-1}]/{\lambda_l})
\oplus(\oplus^k\Q[t, t^{-1}]),$$ 
where 
we have the following:

\smallbreak
\noindent
(1) 
$\lambda_*\in$ $\Q[t, t^{-1}]$ is not zero,

\noindent
(2)
$\lambda_*$ is not the $\Q[t, t^{-1}]$-balanced class of 1, 

\noindent
(3) 
$k$ is the rank of the free part.     

\smallbreak
Two polynomials $f(t), g(t)\in\Q[t, t^{-1}]$ are said to be 
{\it $\Q[t, t^{-1}]$-balanced } 
if there is an integer $n$ 
and a nonzero rational number $r$ 
such that 
$f(t)=r\cdot t^n\cdot g(t)$.

\smallbreak
Let $H_p(\widetilde X;\Q)$ be 
as above. 
Then the {\it $\Bbb Q[t, t^{-1}]$-$p$-Alexander polynomial} is 

\noindent
$
\left\{
\begin{array}{ll}
\mbox{the $\Bbb Q[t, t^{-1}]$-balanced class of 
the product $\lambda_1\cdot...\cdot\lambda_l$} & 
\mbox{if $k=0$ and  $H_p(\widetilde X;\Q)$ is nontrivial} \\
\mbox{0} & \mbox{if $k\neq0$}  \\
\mbox{1} & \mbox{if $H_p(\widetilde X;\Q)\cong0$.}      
\end{array}
\right.
$





\bigbreak
A {\it Seifert hypersurface} for 
an $n$-dimensional oriented closed submanifold  $K$ in $S^{n+2}$ 
is 
an $(n+1)$-dimensional oriented connected compact submanifold in $S^{n+2}$ 
whose boundary is $K$ ($n\in\N$).
Note that there are two cases that 
$K$ is not connected and that $K$ is connected.

\bigbreak
Let $V$ be a Seifert hypersurface for the above $n$-submanifold $K$.    
Note that the orientation of $V$ is compatible with that of $K$. 
Recall that Seifert hypersurfaces are connected by the definition 
(see \S\ref{knotproduct}). 
Let $x_1,..., x_\mu$ be $p$-cycles in $V$ 
which are basis of $H_p(V;\Z)$/Tor. 
Let $y_1,..., y_\nu$ be $(n+1-p)$-cycles in $V$ 
which are basis of $H_p(V;\Z)$/Tor. 
Push $y_i$ to the positive direction of the normal bundle of $V$. 
Call it $y_i^{+}$.  
Push $y_i$ to the negative direction of the normal bundle of $V$. 
Call it $y_i^{-}$.  
A $(p,n+1-p)$-{\it(positive) Seifert matrix} for the above submanifold $K$ associated with $V$ represented by 
an ordered basis,  
$\{x_1,..., x_\mu\}$,   and 
an ordered basis,  
$\{y_1,..., y_\nu\}$,  
is a $(\mu\x\nu)$-matrix 
$$S=(s_{ij})=({\mathrm{lk}}(x_i, y_j^{+})).$$   
We sometimes abbreviate 
$(p,n+1-p)$-Seifert matrix  
to 
$p$-Seifert matrix  
if we know what $n$ is. 
We sometimes let $S_p(K)$ denote a positive $p$-Seifert matrix for 
a closed oriented submanifold $K$ 
and $V$ and $\{x_i\}$ and $\{y_j\}$
if we know what $V$ and $\{x_i\}$ and $\{y_j\}$ are.

\bigbreak
A $(p,n+1-p)$-{\it negative Seifert matrix} 
for the above submanifold $K$ associated with $V$ represented by 
an ordered basis,  
$\{x_1,..., x_\mu\}$,   and 
an ordered basis,  
$\{y_1,..., y_\nu\}$,  
is a matrix 
$$N=(n_{ij})=({\mathrm{lk}}(x_i, y_j^{-})).$$   
%
%
%
We sometimes let $N_p(V)$ denote a negative $p$-Seifert matrix for 
a closed oriented submanifold $K$ 
and $V$ and $\{x_i\}$ and $\{y_j\}$
if we know what $K$ and $\{x_i\}$ and $\{y_j\}$ are. 
Let $S_p$ and $N_p$ be as above. Then we have the following.  
$S_p-N_p$ represents the map 
 $\{H_{p}(V;\Z)$/Tor\} $\x \{H_{n+1-p}(V;\Z)$/Tor\} $\rightarrow\Z$,  
which is defined by the intersection product. 
We call  $t\cdot S_p-N_p$ the $p$-{\it Alexander matrix}   
for $K$ associated with $V$ represented by 
an ordered basis,  
$\{x_1,..., x_\mu\}$,   and 
an ordered basis,  
$\{y_1,..., y_\nu\}$. 

\noindent 
Note that we sometimes define it to be $S_p-t\cdot N_p$. 
The difference of both is only setting the variables 
because we mainly discuss $\Q[t, t^{-1}]$-balanced-classes as follows. 
All we have to do is to change $t$ with $t^{-1}$.

\begin{prop}\label{square}  
Let $K$ be an $n$-dimensional oriented closed submanifold $\subset S^{n+2}$. 

\noindent
Let $S_p$ be 
a $(p,n+1-p)$-{positive Seifert matrix} for $K$ associated with $V$ represented by 
an ordered basis,  
$\{x_1,..., x_\mu\}$,   and 
an ordered basis,  
$\{y_1,..., y_\nu\}$.

\noindent
Let $N_p$ be 
a $(p,n+1-p)$-{negative Seifert matrix} for $K$ associated with $V$ represented by 
an ordered basis,  
$\{x_1,..., x_\mu\}$,   and 
an ordered basis,  
$\{y_1,..., y_\nu\}$.

\noindent
Suppose $\mu=\nu.$
Suppose that 
the linear map defined by a $(p-1)$-Alexander matrix is injective.

Then 
the $p$-$\Q[t, t^{-1}]$-Alexander polynomial is 
the $\Q[t, t^{-1}]$-balanced class of 
`the determinant of $p$-Alexander matrix,'  

$${\rm{det}} (t\cdot S_p-N_p).$$ 
\end{prop}

\noindent
{\bf Note.}  
Of course  $\mu\neq\nu$ in general. 

\smallbreak
\noindent
{\bf Proof.} 
Take the above $X=\overline{S^{n+2}-K\x D^2}$, $\widetilde X$, $V$.   
Let $V\x[-1,1]$ be the tubular neighborhood of $V$ in $X$. 
Let $Y=X-V$.  Consider the Meyer-Vietoris exact sequence: 
$$H_\natural(\amalg_{-\infty}^{\infty}V\x[-1,1];\Q)
\to
H_\natural(\amalg_{-\infty}^{\infty}Y;\Q)
\to  H_\natural(\widetilde X;\Q),$$  

\noindent 
where 
$\amalg_{-\infty}^{\infty}V\x[-1,1]$ is the lift of  $V\x[-1,1]$,  
and where  
$\amalg_{-\infty}^{\infty}Y$ is the lift of $Y$. 
This completes the proof. \qed

\bigbreak

Let $N_p$ be 
a $(p,n+1-p)$-{negative Seifert matrix} for $K$ associated with $V$ represented by 
an ordered basis,  
$\{x_1,..., x_\mu\}$,   and 
an ordered basis,  
$\{y_1,..., y_\nu\}$.  
Let $S_{n+1-p}$ be 
a $(n+1-p,p)$-{positive Seifert matrix} for $K$ associated with $V$ represented by 
an ordered basis,  
$\{y_1,..., y_\nu\}$,   and 
an ordered basis,  
$\{x_1,..., x_\mu\}$. 
By the definition of 
$x_i^{+}$ and  $y_i^{-}$, 
${\mathrm{lk}}(y_i, x_j^{+})$ 
$={\mathrm{lk}}(y_i^{-}, x_j)$. 
By page 541 of \cite{Levinepol}, 
${\mathrm{lk}}(y_i^{-}, x_j)$
$=(-1)^{p(n+1-p)+1}{\mathrm{lk}}(x_j, y_i^{-})$. 
Hence 
$$N_p=(-1)^{p\cdot n+1}S_{n+1-p}$$ 
(note that $p(1-p)$ is an even number).

\begin{prop}\label{Mars}   
Let $K$ be a $(2m+1)$-dimensional closed oriented submanifold  $\subset S^{2m+3}$. 
Let $S$ be an $(m+1,m+1)$-Seifert matrix. 
We have 
$$S=(-1)^{m}\cdot ^t \hskip-1mm S.$$

\end{prop}

The {\it signature} $\sigma(K)$ for 
a $(2p+1)$-dimensional closed oriented submanifold 
$K$\newline$
 \subset S^{2p+3} (p\in\N\cup\{0\}$)
is the signature of 
the matrix $S_{p+1}(K)+^t \hskip-1mm S_{p+1}(K)$.  Therefore, we have the following.

\begin{cla}\label{spring}    
Let $K$ be 
a $(4k+3)$-dimensional closed oriented submanifold 
$\subset S^{4k+5}$  $(k$\newline$
\in\N\cup\{0\})$.  
Then the signature of $K$ coincides with the signature of $\hat V$,  
where $\hat V$ is the closed oriented manifold 
which we obtain 
by attaching a $(4k+4)$-dimensional 0-handles to $\partial V$. 
\end{cla}

Let $K$ be a $(4k+1)$-dimensional spherical knot 
$(4k+1\geq1. \quad k\in\N\cup\{0\})$. 
We regard naturally 
$(H_{2k+1}(V;\Z)/\mathrm{Tor}) \otimes {\Z_2}$ 
as a subgroup of $H_{2k+1}(V;\Z_2)$.  
Then we can take basis  $x_1,..., x_\nu, y_1,..., y_\nu$ 
of 
$(H_{2k+1}(V;\Z)/\mathrm{Tor}) \otimes\Z_2$ 
such that 
$x_i\cdot x_j=0$, $y_i\cdot y_j=0$,  
$x_i\cdot y_j=\delta_{ij}$ 
for any pair $(i, j)$, 
where $\cdot$ denotes the $\Z_2$-intersection product. 
The {\it Arf  invariant} of $K$ is 

$$\text{mod 2 } \Sigma_{i=1}^\nu {\mathrm{lk}}(x_i, x_i^{+})\cdot {\mathrm{lk}}(y_i, y_i^{+}).$$



\bigbreak
Let $L=(L_1,..., L_\mu)$ be a $(4k+1)$-link 
$(4k+1\geq1. \quad k\in\N\cup\{0\}. \quad \mu\in\N-\{0\}.)$.  
We define the {Arf } invariant of $L$.  
There are two cases. 

\begin{itemize}
\item[(1)]
 Let $4k+1\geq5$. 
The Arf  invariant of $L$ is defined in the same manner
as the knot case. 

\item[(2)]
 Let $4k+1=1$. 
See Appendix of \cite{Kirby} and Note right above Note 1.2.1 of \cite{Ogasa02}.   
\end{itemize}

\section{Some results on invariants of $n$-knots and local moves on $n$-knots}
\label{RIL}   

\begin{thm}\label{knots}  
$($\cite{Ogasa09}.$)$ 
Let $K_+$, $K_-$ be spherical $n$-knots $\subset S^{n+2}$ $(n\in\N)$. 
Let $K_0$ be an $n$-submanifold $\subset S^{n+2}$. 
Suppose that  $(K_+$, $K_-$, $K_0)$ is related by a $(p,n+1-p)$-pass-move.
Let $p\neq n+1-p$.  
Then 
there is a polynomial 
$\Delta_p(K_*) 
\in \Q[t, t^{-1}]$ 
whose $\Q[t, t^{-1}]$-balanced class is 
the $p$-$\Q[t, t^{-1}]$-Alexander polynomial 
$A_p(K_*)$  
for the submanifold
$K_*$  
$(*=+,-,0)$  
such that 

$$
\Delta_p(K_+)-\Delta_p(K_-)=(t-1)\cdot \Delta_p(K_0)
.$$
\end{thm}

\begin{thm}\label{middle} 
 $($\cite{Ogasa09}.$)$  
Let $K_+$, $K_-$ be spherical $(4k+1)$-knots $(4k+1\in\N,  \quad k\in\N\cup\{0\})$. 
Let $K_0$ be a $(4k+1)$-submanifold $\subset S^{4k+3}$.
Suppose that  $(K_+$, $K_-$, $K_0)$ is related by a twist-move.  
Then 
there is a polynomial  
$\Delta_{2k+1}(K_*)
\in \Q[t, t^{-1}]$ 
whose $\Q[t, t^{-1}]$-balanced class is 
the $(2k+1)$-$\Q[t, t^{-1}]$-Alexander polynomial 
$A_{2k+1}(K_*)$  
$(*=+,-,0)$   
such that 

$$
\Delta_{2k+1}(K_+)-\Delta_{2k+1}(K_-)=(t-1)\cdot \Delta_{2k+1}(K_0)
.$$
\end{thm}


\bigbreak
We prove the following  Theorem \ref{pepper}.   
The `$l=$even' case is proved in \cite{Ogasa09}, which is the above Theorem \ref{middle}.   
In \S\ref{propol} of this paper we prove the `$l=$odd' case. 

\bigbreak
\noindent {\bf  Theorem \ref{pepper}.}   
{\it 
Let $K_+$ be a $(2l+1)$-dimensional spherical knot $\subset S^{2l+3} (l\in\N\cup\{0\})$. 
Let $K_-$, $K_0$ be $(2l+1)$-dimensional submanifolds $\subset S^{2l+3}$. 
Let $(K_+, K_-, K_0)$ be a twist-move-triple. 
Then there is a polynomial  
$\Delta_{l+1}(K_*)
\in \Q[t, t^{-1}]$ 
whose $\Q[t, t^{-1}]$-balanced class is 
the $(l+1)$-$\Q[t, t^{-1}]$-Alexander polynomial 
$A_{l+1}(K_*)$  
$(*$\newline$=+,-,0)$  
and that 
$$
\Delta_{l+1}(K_+)-\Delta_{l+1}(K_-)=(t+(-1)^{l+1})\cdot \Delta_{l+1}(K_0).  
$$
}

\begin{thm}\label{inertia}   
$($\cite{Ogasa09}.$)$ 
Let $K_+$, $K_-$, $K_0$ be as in Theorem \ref{middle}. 
Let $bP_{4k+2}$ be the $bP$-subgroup $\subset \Theta^{4k+1}$.
Suppose that $bP_{4k+2}$ is not congruent to the trivial group. 
Then we have 
$${\mathrm{Arf}}K_+-{\mathrm{Arf}}K_-=\{|bP_{4k+2}\cap I(K_0)|+1\} \hskip1mm {\rm mod 2},$$

\noindent 
where $I(K_0)$ is the inertia group of an oriented smooth manifold 
which is orientation preserving diffeomorphic to $K_0$ 
and the symbol $|\quad|$ denotes the order of a group.
\end{thm}

\smallbreak
\noindent{\bf{Note.}}  
See \S\ref{Alex} and 
\cite{KervaireMilnor, Kirby, Levinepol}
for the Arf invariant. 
See \cite{KervaireMilnor} for 
the homotopy sphere group $\Theta^{\star}$ and 
the $bP$-subgroup $\subset \Theta^{4k+1}$. 
See 
\cite{BS, Kk} 
for the inertia group.

\smallbreak

We state a problem. 

\begin{prob}\label{unique}   
Let $K_+$ be an $n$-dimensional spherical knot. 
Suppose that $(K_+, K_-, K_0)$ is a twist-move-triple 
(resp. $(p,n+1-p)$-pass-move-triple, where $p\neq n+1-p$). 
Let $\alpha(K)$ be an invariant of $K$ as a submanifold and 
 be a $\Q[t, t^{-1}]$-balanced class.  
Suppose that 
there are $f_+, f_-, f_0\in \Q[t, t^{-1}]$ such that 
the $\Q[t, t^{-1}]$-balanced class of $f_*$  is 
 $\alpha(K_*)$  
($*=+,-,0$)  
and that

\noindent
$
\begin{cases}
\text{$f_+-f_-=(t-1)\cdot f_0$}&\text{in the other cases than the following}\\
\text{$f_+-f_-=(t+1)\cdot f_0$}&\text{if $n=4k+3$, $k\in\N\cup\{0\}$ and if we consider 
the twist-move.}\\
\end{cases}
$

\noindent 
Then  
is $\alpha(K)$ the $\Q[t, t^{-1}]$-Alexander polynomial 
or what is determined by the $\Q[t, t^{-1}]$-Alexander polynomial? 
\end{prob}

\noindent{\bf{Note.}}  
It is well-known that 
the Alexander-Conway polynomial  
(resp. the Jones polynomial) of 1-links is essentially characterized by the well-known local move identity and the fact that it is trivial for the trivial knot.

\section{Theorems on relations between crossing-changes and knot products}\label{thecr} 

\begin{thm}\label{Chicago}   
Suppose that two 1-links $J$ and $K$ differ by one crossing-change.
Then
the 3-submanifolds $\subset S^5$,  
$J\otimes[2]$ 
and 
$K\otimes[2]$,   
differ by one twist-move, 
where $[2]$ denotes the empty knot $[2]$. 
\end{thm}

We will show an example of the phenomenon which Theorem  
\ref{Chicago} asserts after we state Theorem \ref{Illinois} and \ref{Tokyo}.

\begin{thm}\label{Illinois} 
Take the same $J, K$ in Theorem \ref{Chicago}.
Then
the $(2\nu+1)$-submanifolds $\subset S^{2\nu+3}$, 
$J\otimes^\nu[2]$  
and   
$K\otimes^\nu[2]$,   
differ by one twist-move, 
where $\nu\in\N\cup\{0\}$, and, 
where $[2]$ denotes the empty knot $[2]$. 
\end{thm}

\noindent{\bf Note.}  
It is true even in the $\nu=0$ case. 
Of course the $\nu=1$ case is Theorem \ref{Chicago}.

\begin{thm}\label{Tokyo}  
 Let $m\in\N\cup\{0\}$. 
Suppose that two $(2m+1)$-dimensional closed oriented submanifolds  $\subset S^{2m+3}$, 
$J$ and $K$,  differ by one twist-move.
Then
the $(2m+2\nu+1)$-submanifolds $\subset S^{2m+2\nu+3}$, 
$J\otimes^\nu[2]$  
and 
$K\otimes^\nu[2]$, 
differ by one twist-move,  
where $[2]$ denotes the empty knot $[2]$. 
\end{thm}

\smallbreak
\noindent{\bf Note.}  
Of course the $m=0$ case is Theorem \ref{Illinois}.

\bigbreak
We show an example of the phenomenon which Theorem  
\ref{summerc}, \ref{Chicago} \ref{Illinois}, and \ref{Tokyo}
assert. 

The following knot $T$ is the 1-dimensional trivial knot. 

\vskip5mm
\hskip4cm\includegraphics[width=4cm]{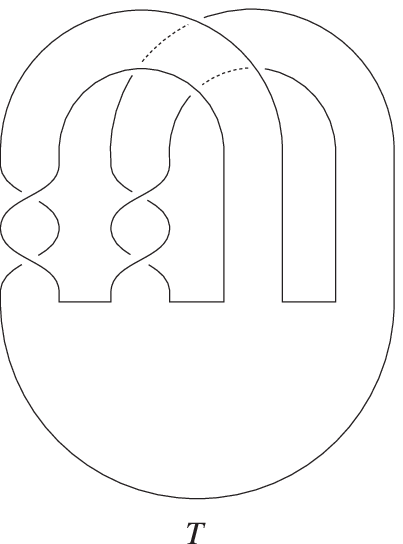}
\vskip5mm

\noindent 
Carry out a crossing-change
 in the 3-ball which is represented by the dotted circle in the following figure.  
We obtain a new knot $K$. Note that $K$ is  the trefoil knot. 

\vskip5mm
\hskip1cm\includegraphics[width=10cm]{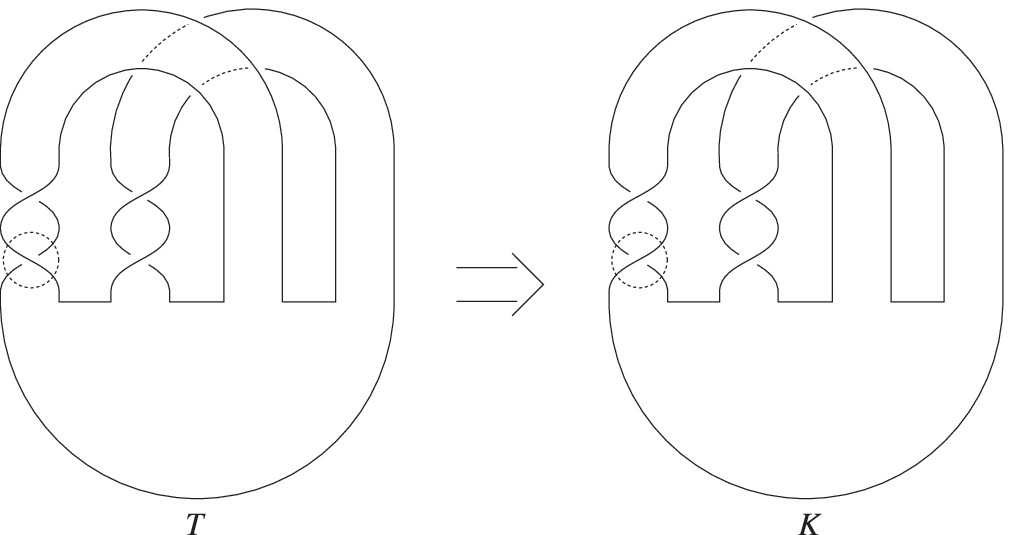}
\vskip5mm

\noindent 
Note that this crossing-change is same as 
the twist-move 
in the 3-ball which is represented by the dotted circle in the following figure.  

\vskip5mm
\hskip1cm\includegraphics[width=11cm]{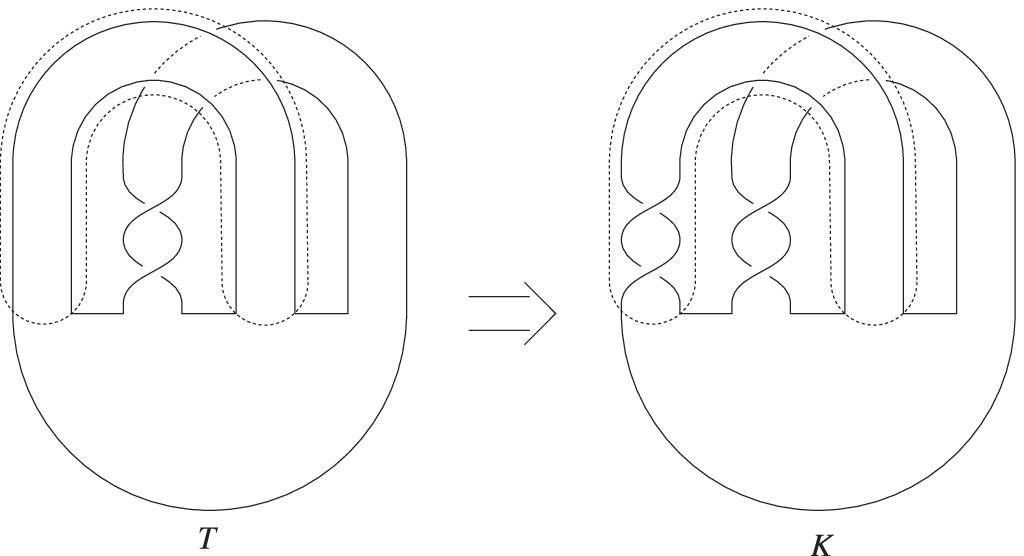}
\vskip5mm

\noindent 
Take $T\otimes[2]$ and  $K\otimes[2]$ in $S^5$. 
By \cite{Kauffman, KauffmanNeumann}   
we can determine the embedding type of $T\otimes[2]$ (resp. $K\otimes[2]$).  
A Seifert hypersurface for $T\otimes[2]$ is diffeomorphic to the following 4-manifold and 
its associated Seifert matrix is 
$\begin{pmatrix}
0&-1\\
0&-1\\
\end{pmatrix}$. 

\vskip5mm\hskip4cm
\unitlength 0.1in
\begin{picture}(18.22,11.30)(6.00,-15.56)
%
\special{pn 8}%
\special{ar 1174 999 553 553  0.2449787 6.2831853}%
\special{ar 1174 999 553 553  0.0000000 0.2204330}%
%
\special{pn 8}%
\special{ar 1870 978 552 552  0.2454152 6.2831853}%
\special{ar 1870 978 552 552  0.0000000 0.2208288}%
%
\special{pn 8}%
\special{sh 0}%
\special{pa 1458 531}%
\special{pa 1572 531}%
\special{pa 1572 637}%
\special{pa 1458 637}%
\special{pa 1458 531}%
\special{ip}%
%
\special{pn 8}%
\special{sh 0}%
\special{pa 1479 1368}%
\special{pa 1593 1368}%
\special{pa 1593 1475}%
\special{pa 1479 1475}%
\special{pa 1479 1368}%
\special{ip}%
%
\special{pn 8}%
\special{pa 1451 616}%
\special{pa 1470 590}%
\special{pa 1494 569}%
\special{pa 1517 549}%
\special{pa 1539 532}%
\special{pa 1543 524}%
\special{sp}%
%
\special{pn 8}%
\special{pa 1593 1354}%
\special{pa 1585 1384}%
\special{pa 1557 1394}%
\special{pa 1531 1409}%
\special{pa 1503 1422}%
\special{pa 1486 1449}%
\special{pa 1472 1461}%
\special{sp}%
%
\special{pn 8}%
\special{ar 1174 1003 553 553  0.2449787 6.2831853}%
\special{ar 1174 1003 553 553  0.0000000 0.2204330}%
%
\special{pn 8}%
\special{ar 1870 982 552 552  0.2454152 6.2831853}%
\special{ar 1870 982 552 552  0.0000000 0.2208288}%
%
\special{pn 8}%
\special{sh 0}%
\special{pa 1458 535}%
\special{pa 1572 535}%
\special{pa 1572 641}%
\special{pa 1458 641}%
\special{pa 1458 535}%
\special{ip}%
%
\special{pn 8}%
\special{sh 0}%
\special{pa 1479 1372}%
\special{pa 1593 1372}%
\special{pa 1593 1479}%
\special{pa 1479 1479}%
\special{pa 1479 1372}%
\special{ip}%
%
\special{pn 8}%
\special{pa 1451 620}%
\special{pa 1470 594}%
\special{pa 1494 573}%
\special{pa 1517 553}%
\special{pa 1539 536}%
\special{pa 1543 528}%
\special{sp}%
%
\special{pn 8}%
\special{pa 1593 1358}%
\special{pa 1585 1388}%
\special{pa 1557 1398}%
\special{pa 1531 1413}%
\special{pa 1503 1426}%
\special{pa 1486 1453}%
\special{pa 1472 1465}%
\special{sp}%
\put(6.0000,-16.0400){\makebox(0,0)[lb]{0}}%
\put(22.1000,-16.7000){\makebox(0,0)[lb]{$-2$}}%
\end{picture}%
\vskip5mm

\noindent 
A Seifert hypersurface for $K\otimes[2]$ is diffeomorphic to the following 4-manifold and 
its associated Seifert matrix is 
$\begin{pmatrix}
-1&-1\\
0&-1\\
\end{pmatrix}$. 

\vskip5mm\hskip4cm
\unitlength 0.1in
\begin{picture}(19.42,11.30)(4.80,-15.56)
%
\special{pn 8}%
\special{ar 1174 999 553 553  0.2449787 6.2831853}%
\special{ar 1174 999 553 553  0.0000000 0.2204330}%
%
\special{pn 8}%
\special{ar 1870 978 552 552  0.2454152 6.2831853}%
\special{ar 1870 978 552 552  0.0000000 0.2208288}%
%
\special{pn 8}%
\special{sh 0}%
\special{pa 1458 531}%
\special{pa 1572 531}%
\special{pa 1572 637}%
\special{pa 1458 637}%
\special{pa 1458 531}%
\special{ip}%
%
\special{pn 8}%
\special{sh 0}%
\special{pa 1479 1368}%
\special{pa 1593 1368}%
\special{pa 1593 1475}%
\special{pa 1479 1475}%
\special{pa 1479 1368}%
\special{ip}%
%
\special{pn 8}%
\special{pa 1451 616}%
\special{pa 1470 590}%
\special{pa 1494 569}%
\special{pa 1517 549}%
\special{pa 1539 532}%
\special{pa 1543 524}%
\special{sp}%
%
\special{pn 8}%
\special{pa 1593 1354}%
\special{pa 1585 1384}%
\special{pa 1557 1394}%
\special{pa 1531 1409}%
\special{pa 1503 1422}%
\special{pa 1486 1449}%
\special{pa 1472 1461}%
\special{sp}%
%
\special{pn 8}%
\special{ar 1174 1003 553 553  0.2449787 6.2831853}%
\special{ar 1174 1003 553 553  0.0000000 0.2204330}%
%
\special{pn 8}%
\special{ar 1870 982 552 552  0.2454152 6.2831853}%
\special{ar 1870 982 552 552  0.0000000 0.2208288}%
%
\special{pn 8}%
\special{sh 0}%
\special{pa 1458 535}%
\special{pa 1572 535}%
\special{pa 1572 641}%
\special{pa 1458 641}%
\special{pa 1458 535}%
\special{ip}%
%
\special{pn 8}%
\special{sh 0}%
\special{pa 1479 1372}%
\special{pa 1593 1372}%
\special{pa 1593 1479}%
\special{pa 1479 1479}%
\special{pa 1479 1372}%
\special{ip}%
%
\special{pn 8}%
\special{pa 1451 620}%
\special{pa 1470 594}%
\special{pa 1494 573}%
\special{pa 1517 553}%
\special{pa 1539 536}%
\special{pa 1543 528}%
\special{sp}%
%
\special{pn 8}%
\special{pa 1593 1358}%
\special{pa 1585 1388}%
\special{pa 1557 1398}%
\special{pa 1531 1413}%
\special{pa 1503 1426}%
\special{pa 1486 1453}%
\special{pa 1472 1465}%
\special{sp}%
\put(22.1000,-16.7000){\makebox(0,0)[lb]{$-2$}}%
\put(4.8000,-16.8000){\makebox(0,0)[lb]{$-2$}}%
\end{picture}%
\vskip5mm

\noindent 
$K\otimes[2]$ is obtained from $T\otimes[2]$ by a twist-move 
in the 5-ball 
which is represented by the dotted circle in the following figure.  

\vskip5mm
\includegraphics[width=12cm]{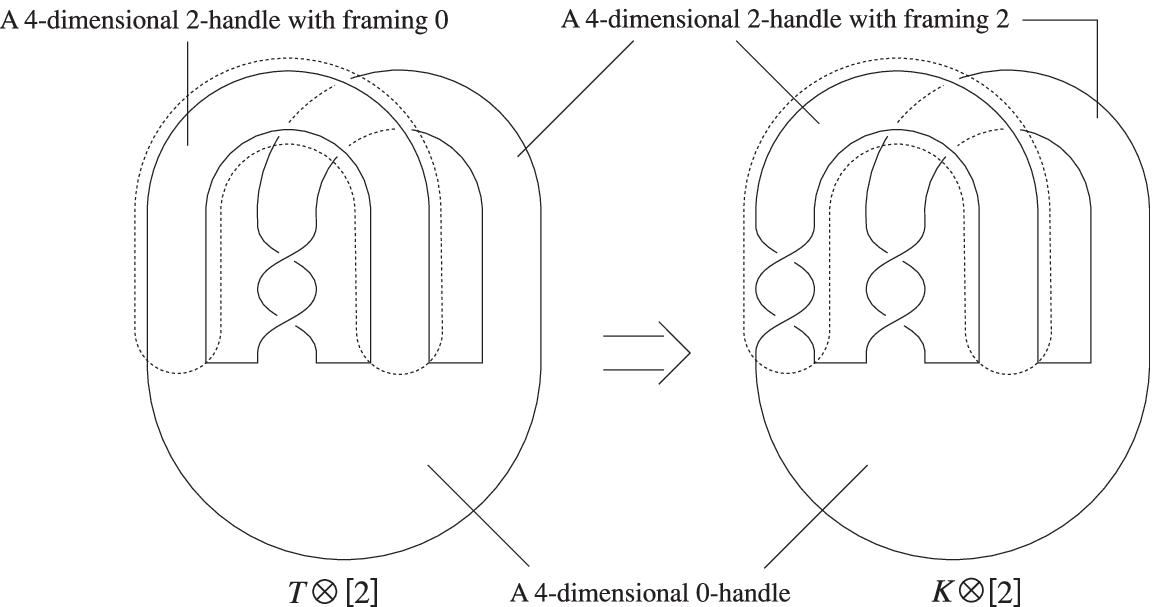}
\vskip5mm

\noindent 
Note that $K\otimes[2]$ is not homeomorphic to $T\otimes[2]$.  
Twist moves of $(4k+3)$-knots ($k\in\N\cup\{0\}$) 
change the homeomorphism types of submanifolds in general 
but we can determine the new embedding types which we obtain by twist-moves.

\vskip5mm

$K\otimes^{2\mu}[2]$= $K\otimes^\mu$(the Hopf link) in $S^{4\mu+3}$ 
is obtained from \newline
$T\otimes^{2\mu}[2]$= $T\otimes^\mu$(the Hopf link) in $S^{4\mu+3}$ 
by a twist-move 
in the $(4\mu+3)$-ball 
which is represented by the dotted circle in the following figure.  

\vskip5mm
\includegraphics[width=12cm]{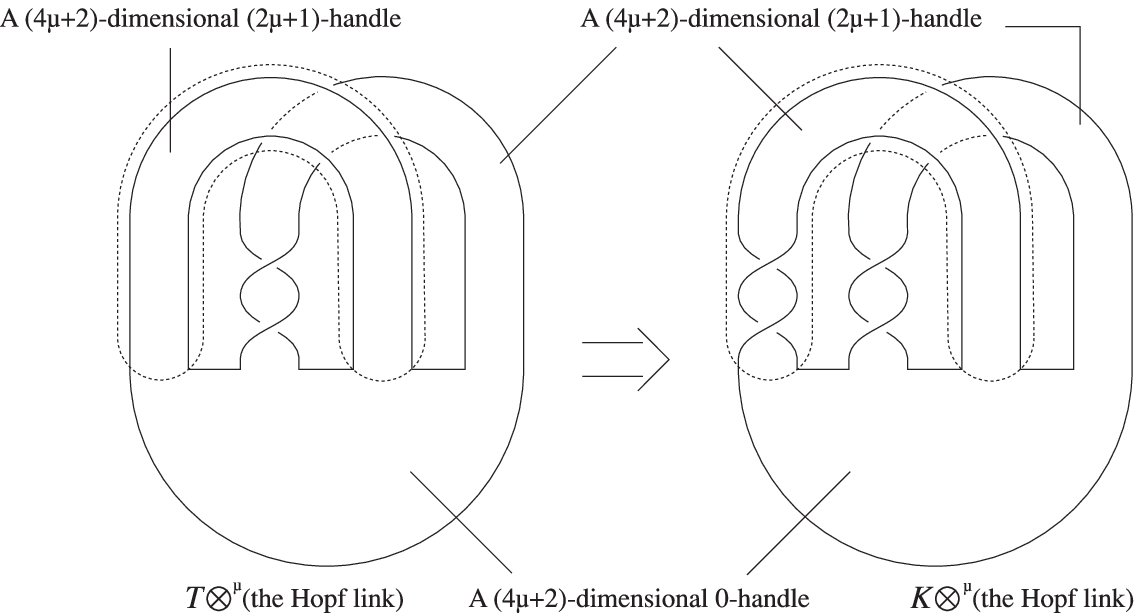}
\vskip5mm

\noindent
Note that a Seifert hypersurface for $T\otimes^\mu$(the Hopf link)
is diffeomorphic to 
the plumbing of the trivial $D^{2\mu+1}$-bundle over $S^{2\mu+1}$ and 
the $D^{2\mu+1}$-bundle over $S^{2\mu+1}$ 
associated with the tangent bundle. 

\noindent
Note that a Seifert hypersurface for $K\otimes^\mu$(the Hopf link)  
is diffeomorphic to 
the plumbing of two copies of 
the $D^{2\mu+1}$-bundle over $S^{2\mu+1}$ 
associated with the tangent bundle.

\noindent
For some $\mu$,  
$K\otimes^\mu$(the Hopf link)  is not diffeomorphic but homeomorphic to \newline
$T\otimes^\mu$(the Hopf link). Then $K\otimes^\mu$(the Hopf link) is an exotic sphere. 
For other $\mu$,  
$K\otimes^\mu$(the Hopf link)  is diffeomorphic to 
$T\otimes^\mu$(the Hopf link). 
Recall the discussion associated with the $bP$-subgroup in \cite{KervaireMilnor}.


\begin{thm}\label{Sunday}    
$(1)$  
There is a nontrivial 1-knot $K$ which is obtained  from the trivial knot 
by a single crossing-change 
with the following property. 
The $(2\nu+1)$-submanifold $\subset S^{2\nu+3}$, 
$K\otimes^\nu[2]$,   
is equivalent to the trivial $(2\nu+1)$-knot,  
where $\nu\in\N, \quad \nu\geqq2$, and, 
where $[2]$ denotes the empty knot $[2]$. 

\bigbreak
\noindent 
$(2)$  
There is a nontrivial 1-knot $P$ which is obtained  from the trivial knot 
by two crossing-changes 
not by a single crossing-change 
with the following property. 
The $(2\nu+1)$-submanifold $\subset S^{2\nu+3}$, 
$P\otimes^\nu[2]$,   
is equivalent to the trivial $(2\nu+1)$-knot,  
where $\nu\in\N, \quad \nu\geqq2$.

\bigbreak
\noindent 
$(3)$ 
There are nontrivial 1-knots $P$ and $Q$ with the following properties. 

\smallbreak
\noindent
$\rm{(i)}$ 
$P$ and $Q$ differ by a single crossing-change and are nonequivalent.

\noindent
$\rm{(ii)}$   Let $\nu\in\N$ and $\nu\geqq2$.   
The $(2\nu+1)$-submanifolds $\subset S^{2\nu+3}$, 
$P\otimes^\nu[2]$ 
and 
$Q\otimes^\nu[2]$ 
are equivalent spherical knots.

\bigbreak
\noindent 
$(4)$ 
There is a nontrivial 1-knot $P$ with the following properties. 

\noindent
$\rm{(i)}$ 
$P$ is obtained from the trivial 1-knot by two crossing-changes 
not by a single crossing-change.

\noindent
$\rm{(ii)}$ Let $\mu\in\N$. 
The $(4\mu+1)$-submanifold $\subset S^{4\mu+3}$, 
$P\otimes^\mu($the Hopf link$)$ 
is obtained from the trivial knot by a twist-move and is a nontrivial knot.

\bigbreak
\noindent 
$(5)$ 
There are nontrivial 1-knots $P$ and $Q$ with the following properties. 

\smallbreak
\noindent
$\rm{(i)}$ 
$P$ and $Q$ differ by a single crossing-change and are nonequivalent.

\noindent
$\rm{(ii)}$ Let $\mu\in\N$. 
The $(4\mu+1)$-submanifolds $\subset S^{4\mu+3}$, 
$P\otimes^\mu($the Hopf link$)$ 
and 
$Q\otimes^\mu($the Hopf link$)$ 
are  equivalent spherical knots and nontrivial knots.   

\end{thm}

\noindent{\bf  Note.} 
Recall that, 
if 
a 1-knot $K$ is obtained from a nonequivalent 1-knot $J$ 
by a single crossing-change, 
then  
the 1-knot $K$ is obtained from the nonequivalent knot $J$ by a single twist-move. 
Thus we can say that 
the `twist-move-unknotting-number' changes by a knot product.
A `non-twist-move-equivalent pair' 
is changed into 
a `twist-move-equivalent pair' by knot product.

\begin{note}\label{tori}    
Compare 
the above Theorem \ref{Sunday} 
with the following $(*)$. 

\noindent
$(*)$
There is a nontrivial 1-knot $K$ 
which is obtained from the trivial knot by one crossing-change  
with the following properties:    

\smallbreak
\noindent 
(1)
$K\otimes[2]$ is not a spherical knot. Hence it is a nontrivial knot.

\noindent 
(2) $K\otimes^\nu[2]$ ($\nu\geqq2, \nu\in\N$)  
is a spherical knot and the trivial knot. 
\smallbreak

\noindent We prove this in \S\ref{procr}. 
\end{note}


\begin{thm}\label{aoiro}    
Let $k\in\N$.    
Let $K$ $($resp. $J)$ be $(4k+5)$-dimensional smooth submanifold $\subset S^{4k+7}$. 
Suppose that $K$ and $J$ differ by a single twist-move
and  are nonequivalent. 
Suppose that $K$ is equivalent to 
$A\otimes^{k+1}($the Hopf link$)$ for a 1-knot $A$. 

Then  
there is a unique equivalence class of simple $(4k+1)$-knots for $K$ $($resp. $J)$   
with the following properties. 

\smallbreak
\noindent $\mathrm{(i)}$
There is a representative element $K'$ of the above equivalence class for $K$ such that 
$K$ is equivalent to $K'\otimes($the Hopf link$)$.   

\noindent $\mathrm{(ii)}$
There is a representative element $J'$ of the above equivalence class for $J$ such that 
$J$ is equivalent to $J'\otimes($the Hopf link$)$.   

\noindent $\mathrm{(iii)}$
$K'$ and $J'$ differ 
by a single twist-move 
and are nonequivalent. 
\end{thm}

\noindent
{\bf  Note.} 
If $k=0$ in Theorem \ref{aoiro}, we have a different result. 
See Theorem \ref{Sunday}.

\section{Theorems on relations between pass-moves and knot products}\label{thepa}  

\begin{thm}\label{Waltham}  
Suppose that two 1-knots $J$ and $K$ differ by one pass-move.
Then 
the $(4\mu+1)$-submanifolds $\subset S^{4\mu+3}$, 
$J\otimes^\mu($the Hopf link$)$ 
and \newline 
$K\otimes^\mu($the Hopf link$)$,
differ by one $(2\mu+1,2\mu+1)$-pass-move $(\mu\in\N\cup\{0\})$.
\end{thm}

\begin{note}\label{moon}   
Recall Theorem \ref{Hopf}. For any $n$-dimensional closed oriented submanifold $A$\newline 
$\subset S^{n+2}$, 
$A$
\hskip-1mm
$\begin{matrix}
\text{{\hskip5mm\tiny $\mu$}}\\
\text{\large $\otimes$}\\
\text{  }\\
\end{matrix}$
\hskip-1mm
$($the negative Hopf link$)=$
\hskip-1mm
$A$
\hskip-2mm
$\begin{matrix}
\text{{\hskip5mm\tiny $2\mu$}}\\
\text{\large $\otimes$}\\
\text{  }\\
\end{matrix}$
\hskip-1mm
$[2]$, 
where $n\in\N$ and $\mu\in\N\cup\{0\}$. 
\end{note}

Compare the above Theorem \ref{Waltham} and Note \ref{moon} 
with the following Theorem \ref{sky}.

\begin{thm}\label{sky}  
Take the same $J, K$ in Theorem \ref{Waltham}.
Then 
the $(4\mu+3)$-submanifolds $\subset S^{4\mu+5}$ $(\mu\in\N\cup\{0\})$, 
$J$ \hskip-8mm 
$\begin{matrix}
\text{{\hskip9mm\tiny $(2\mu+1)$}}\\
\text{\large $\otimes$}\\
\text{  }\\
\end{matrix}$\hskip-1mm
$[2]$  
and   
$K$ \hskip-8mm
$\begin{matrix}
\text{{\hskip9mm\tiny $(2\mu+1)$}}\\
\text{\large $\otimes$}\\
\text{  }\\
\end{matrix}$\hskip-1mm
$[2]$, 
are not homeomorphic in general 
and,  therefore, 
are NOT  $(2\mu+2,2\mu+2)$-pass-move-equivalent in general. 
\end{thm}

\begin{prob}\label{sea}   
In the above Theorem \ref{sky}, 
of course, 
if $J$ and $K$ are trivial knots,
then the above two $(4\mu+3)$-submanifolds are pass-move-equivalent.
What kind of pair, $J$ and $K$, in Theorem \ref{sky} satisfies
the condition that  
the above two $(4\mu+3)$-submanifolds 
are 
$(2\mu+2,2\mu+2)$-pass-move-equivalent?

\end{prob}

The `$\nu=$odd' case of Theorem 
\ref{grape} 
and Note \ref{bird} give partial solutions to Problem \ref{sea}.

\begin{thm}\label{mountain}    
Let $J, K$ be simple $(2l+1)$-knots, where $l\in\N$.   
Suppose that $J$ and $K$ differ by one $(l+1, l+1)$-pass-move.
Then
the $(2l+4\mu+1)$-submanifolds $\subset S^{2l+4\mu+3}$, 
$J\otimes^\mu($the Hopf link$)$    
and
$K\otimes^\mu($the Hopf link$)$,   
differ by one $(l+2\mu+1, l+2\mu+1)$-pass-move.
\end{thm}

\begin{prob}\label{lion}    
If we do NOT suppose that $J, K$ are simple knots in Theorem \ref{mountain},
do the above $(2l+4\mu+1)$-submanifolds 
differ by one pass-move?
Or, are they pass-move-equivalent?
\end{prob}

The above problem is really a problem of high dimensional knots.
The following one is also such a problem.

\vskip3mm
\begin{prob}\label{tiger}   
  (A generalization of Problem \ref{lion}.)
(1) Suppose that spherical $n$-knots 
(resp. $n$-dimensional closed oriented submanifolds $\subset S^{n+2}$) 
$J$ and $K$ differ by one $(p,n+1-p)$-pass-move, 
where $n\in\N$ and $p\in\N$.
Then
do
the $(n+4\mu+1)$-submanifolds $\subset S^{n+4\mu+3}$, 
$J\otimes^\mu($the Hopf link$)$    
and
$K\otimes^\mu($the Hopf link$)$,   
differ by one pass-move?  
Or, are they pass-move-equivalent?

\noindent
(2) How about the case where $J$ and $K$ are even dimensional simple $2m$-knots
and where $p=m$? Here, $m\in\N$. 

\noindent
(3) Of course there is a problem
in the case of the product with odd times copies of the empty knot $[2]$.
(Note Theorem \ref{sky} and its Proof.)
\end{prob}

\begin{thm}\label{grape}  
There is a nontrivial 1-knot $K$ 
which is obtained from the trivial knot by one pass-move 
with the following property:
$K\otimes^\nu[2]$
is the trivial $(2\nu+1)$-knot, where $\nu\geqq2, \quad\nu\in\N$.
\end{thm}

\smallbreak
\begin{note}\label{pine}    
By the above Theorem \ref{grape}, 
we have the following. 
Let $T$ be the trivial 1-knot.
The two 1-knots, $K$ and $T$, differ by a single pass-move 
and are nonequivalent.  
However the $(2\nu+1)$-dimensional spherical knot,   
$K\otimes^\nu[2]$,  
is equivalent to 
the $(2\nu+1)$-dimensional trivial knot, 
$T\otimes^\nu[2]$. 
Recall  $\nu\geqq2$ and $\nu\in\N$.
%
%
That is, 
they differ by ZERO times of pass-moves. 
Thus we can say that 
the `pass-move-unknotting-number' changes by knot products.
\end{note}

\smallbreak
\begin{note}\label{bird}    
Compare the above Theorem \ref{grape}   
with the following $(*)$. 

\noindent
$(*)$
There is a nontrivial 1-knot $K$ 
which is obtained from the trivial knot by one pass-move 
with the following properties:    

\smallbreak
\noindent 
(1)
$K\otimes[2]$,
is not a spherical knot. Hence it is a nontrivial knot.

\noindent 
(2) $K\otimes^\nu[2]$ ($\nu\geqq2, \nu\in\N$) 
is a spherical knot and the trivial knot. 
\smallbreak

\noindent We prove this in \S\ref{propa}. 
\end{note}


\begin{thm}\label{ki}   
Let $\mu\in\N$. 
Let $K$ $($resp. $J$$)$ be a $(4\mu+1)$-submanifold $\subset S^{4\mu+3}$. 
Let $K$ and $J$ be $(2\mu+1, 2\mu+1)$-pass-move-equivalent. 
Suppose that 
$K$ is equivalent to $K'\otimes^\mu($the Hopf link$)$ for a 1-knot $K'$.
Then  there is a 1-knot $J'$ with the following properties. 

\smallbreak
\noindent $\mathrm{(i)}$
$J$  is equivalent to $J'\otimes^\mu($the Hopf link$)$.   

\noindent $\mathrm{(ii)}$
$K'$ and $J'$ are pass-move-equivalent. 
\end{thm}

\noindent
{\bf Note.} 
If there is such a 1-knot $K'$, 
there are countably infinitely many nonequivalent 1-knots $P$ (resp. $Q$) 
such that 
$P\otimes^\mu($the Hopf link$)$ 
(resp. $Q\otimes^\mu($the Hopf link$)$)
is equivalent to $K$ (resp. $J$). 
We prove this in \S\ref{propa}.

\begin{thm}\label{aka}   
Let $p\in\N$.    
Let $K$ and $J$ be $(2p+5)$-dimensional smooth submanifolds $\subset S^{2p+7}$. 
Suppose that $K$ and $J$ differ 
by a single $(p+3, p+3)$-pass-move 
and 
are nonequivalent. 
Suppose that 
$K$ is equivalent to \newline
$
\begin{cases}
\text{$A\otimes^{\frac{p}{2}+1}($the Hopf link$)$ for a 1-knot $A$}&\text{if $p$ is even}\\
\text{$A\otimes($the Hopf link$)$ for a simple 3-knot $A$}&\text{if $p=1$ $($and $2p+5=7)$}\\
\text{$A\otimes^{\frac{p-1}{2}}($the Hopf link$)$ for a simple 7-knot $A$}&\text{if $p$ is odd and $p\neq1$.}
\end{cases}
$

Then  
there is a unique equivalence class of simple $(2p+1)$-knots for $K$ $($resp. $J)$   
with the following properties. 

\smallbreak
\noindent $\mathrm{(i)}$
There is a representative element $K'$ of the above equivalence class for $K$ such that 
$K$ is equivalent to $K'\otimes($the Hopf link$)$.   

\noindent $\mathrm{(ii)}$   
There is a representative element $J'$ of the above equivalence class for $J$ such that  
$J$  is equivalent to $J'\otimes($the Hopf link$)$.   

\noindent $\mathrm{(iii)}$
$K'$ and $J'$ differ by a single $(p+1,p+1)$-pass-move 
and are nonequivalent. 
\end{thm}

\bigbreak

Let $P$ be the 5-twist spun knot of the trefoil knot.  Note that $P$ is a 2-knot.    
Note that 
Proposition 4.3 of \cite{Ogasa04} 
and the last line of \S 7 in page 684 of \cite{Ogasa04} imply the following. 

\smallbreak
\noindent
(1)
$P$ is NOT ribbon-move-equivalent to $T$.  

\noindent
(2)
$P$ is NOT (1,2)-pass-move-equivalent to $T$.   


\begin{thm}\label{spun}  
Let $T$ be the trivial 2-knot.   
Let $P$ be as above. 
Although $P$ is NOT (1,2)-pass-move-equivalent to $T$,  we have the following.
Let $\nu\geqq2$ and $\nu\in\N$. 
The $(2\nu+2)$-submanifold, 
$P\otimes^\nu[2]$, 
is equivalent 
to the trivial  $(2\nu+2)$-knot, 
$T\otimes^\nu[2]$, 
and therefore,
is $(\nu+1, \nu+2)$-pass-move-equivalent   to the trivial knot. 
\end{thm}
\vskip3mm

\noindent{\bf Note.}  
Thus we can say that a knot product makes 
a `non-pass-move-equivalent pair' into
a `pass-move-equivalent pair'.

\section{Theorems on relations between local move identities of a knot polynomial
and knot products}\label{thepol}   
%
%
%
%
Let $(K_+, K_-, K_0)$ be a crossing-change-triple of 1-links. 
(See crossing-change-triple for \S\ref{Introduction}.)
Let $A(K)$ be the Alexander-Conway polynomial of 1-links $K$. 
It is well-known that  
$$A(K_+)-A(K_-)=(t-1)\cdot A(K_0).$$
Here, we have the following Theorems. 

\begin{thm}\label{skein1}    
Let $K_+, K_-, K_0$ be as above. 
There is a polynomial
$
\Delta_{2\mu+1}(\>K_*\otimes^\mu
(\text{the Hopf link})\>) 
\in \Q[t, t^{-1}]$ 
whose $\Q[t, t^{-1}]$-balanced class is 
the $(2\mu+1)$-$\Q[t, t^{-1}]$-Alexander polynomial 
$A_{2\mu+1}(\>K_*\otimes^\mu
(\text{the Hopf link})\>) 
$  
$(\mu\in\N\cup\{0\}. \quad  *=+,-,0.)$   
such that 

\hskip21mm
$
\Delta_{2\mu+1}(\>K_+\otimes^\mu(\text{the Hopf link})\>) -
\Delta_{2\mu+1}(\>K_-\otimes^\mu(\text{the Hopf link})\>) 
$

\hskip18mm
$
=
(t-1)\cdot \Delta_{2\mu+1}(\>K_0\otimes^\mu(\text{the Hopf link})\>).
$
\end{thm}

\noindent 
{\bf Note.}  
After taking knot product, $A(\quad)$ is changed into $A_{2\mu+1}(\quad)$.   
\smallbreak

\noindent
By  Theorem \ref{Hopf},  Theorem \ref{skein1} follows from  Theorem \ref{skein2}.

\begin{thm}\label{skein2} 
Let $K_+, K_-, K_0$ be as above. 
There is a polynomial 
$\Delta_{\nu+1}(\>K_*
\otimes^\nu[2]\>)
\in \Q[t, t^{-1}]$ 
whose $\Q[t, t^{-1}]$-balanced class is 
the $(\nu+1)$-$\Q[t, t^{-1}]$-Alexander polynomial 
$A_{\nu+1}(\>
K_*
\otimes^\nu[2]\>)  \quad
 (\nu\in\N\cup\{0\}, *=+,-,0)$ 
such  that 
$$
\Delta_{\nu+1}(\>K_+\otimes^\nu[2]\>)
-
\Delta_{\nu+1}(\>K_-\otimes^\nu[2]\>)
=
(t+(-1)^{\nu+1})\cdot \Delta_{\nu+1}(\>K_0\otimes^\nu[2]\>),$$

\noindent 
where  
$[2]$ denotes the empty knot $[2]$. 
\end{thm}

\noindent 
{\bf Note.}  
If $\nu$ is odd, then 
$K_+
\otimes^\nu[2]$ 
is not homeomorphic to 
$K_-
\otimes^\nu[2]$ 
in general. 
However the above theorem is  true.

\bigbreak
We show an example of Theorem  \ref{skein2}. 

Let $K_+$ be the trivial 1-knot,  
$K_-$ the trefoil knot, 
and $K_0$ the Hopf link as shown below.

\vskip5mm
\hskip1cm\includegraphics[width=14cm]{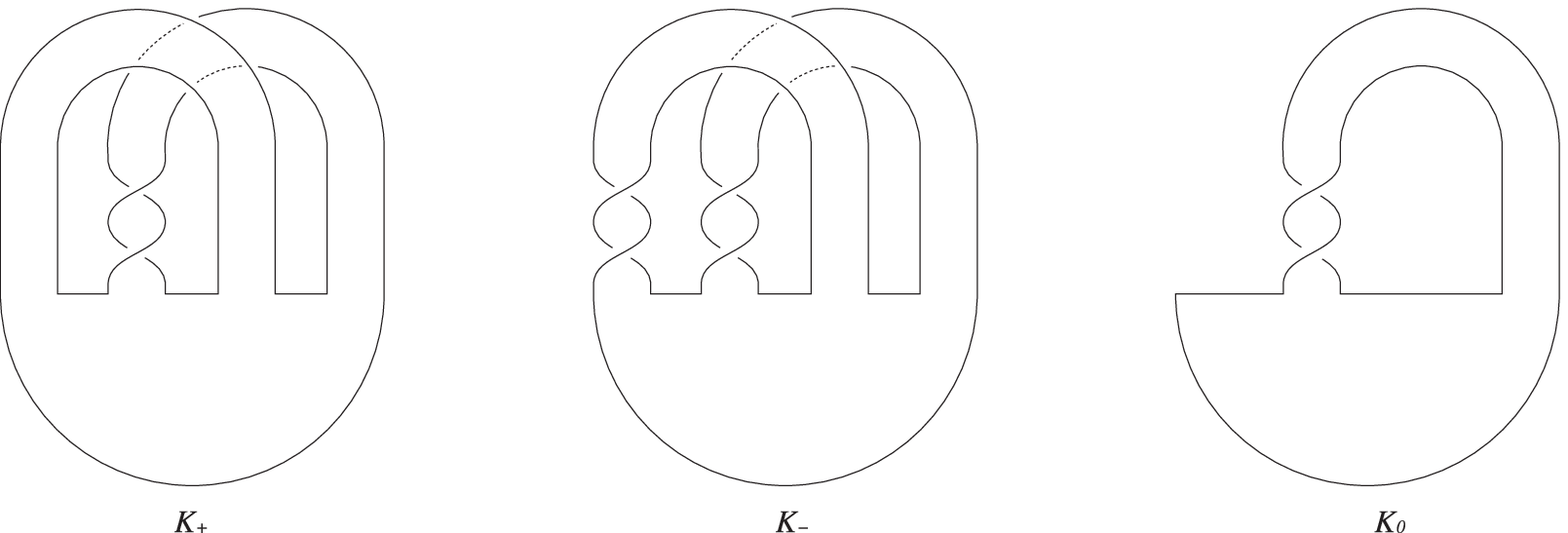}
\vskip5mm

Take $K_*\otimes[2]$ ($*=+,-,0$). 
By \cite{Kauffman, KauffmanNeumann}   
we can determine the embedding type of $K_*\otimes[2]$ ($*=+,-,0$).  
A Seifert matrix for 
$\begin{cases}
K_+\otimes[2]\\
K_-\otimes[2]\\
K_0\otimes[2]
\end{cases}$ 
is 
$\begin{cases}
\begin{pmatrix}
0&-1\\
0&-1\\
\end{pmatrix} \\

\begin{pmatrix}
-1&-1\\
0&-1\\
\end{pmatrix}\\

\quad(-1)
\end{cases}. $

\noindent 
Hence  
the 2-Alexander polynomial for 
$\begin{cases}
K_+\otimes[2]\\
K_-\otimes[2]\\
K_0\otimes[2]
\end{cases}$ 
is the $\Q[t,t^{-1}]$-balanced class of 

\noindent
$\begin{cases}
\mathrm{det} 
\Bigg\{t\begin{pmatrix}
0&-1\\
0&-1\\
\end{pmatrix}
+\begin{pmatrix}
0&0\\
-1&-1\\
\end{pmatrix}\Bigg\}
=-t  \\

\mathrm{det}
\Bigg\{t\begin{pmatrix}
-1&-1\\
0&-1\\
\end{pmatrix}
+\begin{pmatrix}
-1&0\\
-1&-1\\
\end{pmatrix}\Bigg\}=t^2+t+1\\

\mathrm{det}
\{t(-1)+(-1)\}=-t-1
\end{cases}. $
   \newline
Since $-t-(t^2+t+1)=(t+1)(-t-1)$, the identity in 
Theorem \ref{skein2} 
holds for the triple 
$(K_+\otimes[2], K_-\otimes[2], K_0\otimes[2])$.

\bigbreak




\smallbreak

The above Theorem \ref{skein1}, \ref{skein2} are related to the following Theorem \ref{pepper}. 
The `$l=$even' case is proved in \cite{Ogasa09}. 
In this paper we prove the `$l=$odd' case.

\begin{thm}\label{pepper}    
Let $K_+$ be a $(2l+1)$-dimensional spherical knot $\subset S^{2l+3} (l\in\N\cup\{0\})$. 
Let $K_-$, $K_0$ be $(2l+1)$-dimensional submanifolds $\subset S^{2l+3}$. 
Let $(K_+, K_-, K_0)$ be a twist-move-triple. 

Then there is a polynomial  
$\Delta_{l+1}(K_*)\in \Q[t, t^{-1}]$ 
whose $\Q[t, t^{-1}]$-balanced class is 
the $(l+1)$-$\Q[t, t^{-1}]$-Alexander polynomial 
$A_{l+1}(K_*)$  
$(*=+,-,0)$  
and that 
$$
\Delta_{l+1}(K_+)-\Delta_{l+1}(K_-)=(t+(-1)^{l+1})\cdot \Delta_{l+1}(K_0).$$
\end{thm}

\bigbreak
We show an example of Theorem  \ref{pepper}. 
Let $K_+$ be the trivial 1-knot,  
$K_-$ the trefoil knot, 
and $K_0$ the Hopf link as shown below.  
(These $K_*$ ($*=+,-,0$) are same as the examples $K_*$ of  Theorem \ref{skein2}.)  
Note that 
$(K_+, K_-, K_0)$ is a twist-move triple in the 3-ball which is represented by the dotted circle.

\vskip5mm
\includegraphics[width=15cm]{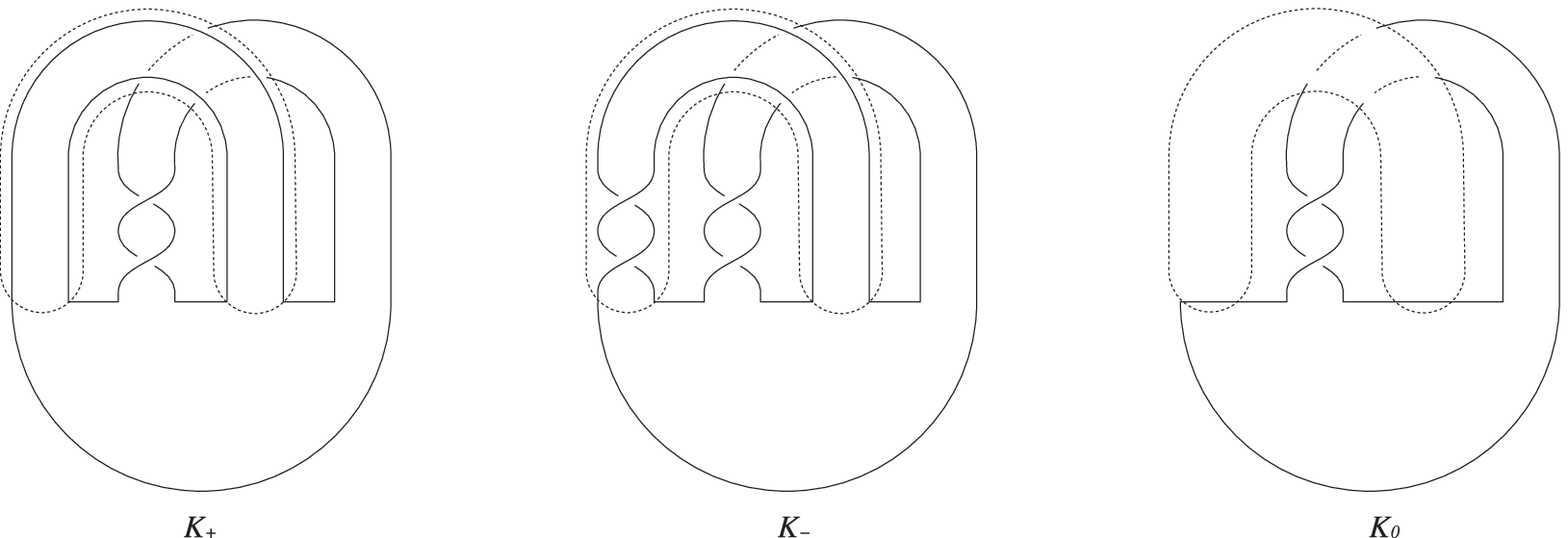}
\vskip5mm

\noindent 
Let $J_*$ be $K_*\otimes[2]$ ($*=+,-,0$). 
By \cite{Kauffman, KauffmanNeumann}   
a Seifert hypersurface for $J_*$ is diffeomorphic to the following 4-manifold.

\vskip5mm
\unitlength 0.1in
\begin{picture}(61.64,21.30)(3.70,-22.90)
%
\special{pn 8}%
\special{ar 1174 999 553 553  0.2449787 6.2831853}%
\special{ar 1174 999 553 553  0.0000000 0.2204330}%
%
\special{pn 8}%
\special{ar 1870 978 552 552  0.2454152 6.2831853}%
\special{ar 1870 978 552 552  0.0000000 0.2208288}%
%
\special{pn 8}%
\special{sh 0}%
\special{pa 1458 531}%
\special{pa 1572 531}%
\special{pa 1572 637}%
\special{pa 1458 637}%
\special{pa 1458 531}%
\special{ip}%
%
\special{pn 8}%
\special{sh 0}%
\special{pa 1479 1368}%
\special{pa 1593 1368}%
\special{pa 1593 1475}%
\special{pa 1479 1475}%
\special{pa 1479 1368}%
\special{ip}%
%
\special{pn 8}%
\special{pa 1451 616}%
\special{pa 1470 590}%
\special{pa 1494 569}%
\special{pa 1517 549}%
\special{pa 1539 532}%
\special{pa 1543 524}%
\special{sp}%
%
\special{pn 8}%
\special{pa 1593 1354}%
\special{pa 1585 1384}%
\special{pa 1557 1394}%
\special{pa 1531 1409}%
\special{pa 1503 1422}%
\special{pa 1486 1449}%
\special{pa 1472 1461}%
\special{sp}%
%
\special{pn 8}%
\special{ar 1174 1003 553 553  0.2449787 6.2831853}%
\special{ar 1174 1003 553 553  0.0000000 0.2204330}%
%
\special{pn 8}%
\special{ar 1870 982 552 552  0.2454152 6.2831853}%
\special{ar 1870 982 552 552  0.0000000 0.2208288}%
%
\special{pn 8}%
\special{sh 0}%
\special{pa 1458 535}%
\special{pa 1572 535}%
\special{pa 1572 641}%
\special{pa 1458 641}%
\special{pa 1458 535}%
\special{ip}%
%
\special{pn 8}%
\special{sh 0}%
\special{pa 1479 1372}%
\special{pa 1593 1372}%
\special{pa 1593 1479}%
\special{pa 1479 1479}%
\special{pa 1479 1372}%
\special{ip}%
%
\special{pn 8}%
\special{pa 1451 620}%
\special{pa 1470 594}%
\special{pa 1494 573}%
\special{pa 1517 553}%
\special{pa 1539 536}%
\special{pa 1543 528}%
\special{sp}%
%
\special{pn 8}%
\special{pa 1593 1358}%
\special{pa 1585 1388}%
\special{pa 1557 1398}%
\special{pa 1531 1413}%
\special{pa 1503 1426}%
\special{pa 1486 1453}%
\special{pa 1472 1465}%
\special{sp}%
\put(22.2000,-17.2000){\makebox(0,0)[lb]{$-2$}}%
\put(6.8000,-16.4000){\makebox(0,0)[lb]{0}}%
%
\special{pn 8}%
\special{pa 2730 160}%
\special{pa 2730 2060}%
\special{fp}%
%
\special{pn 8}%
\special{pa 5150 220}%
\special{pa 5150 2120}%
\special{fp}%
%
\special{pn 8}%
\special{ar 3464 1013 553 553  0.2449787 6.2831853}%
\special{ar 3464 1013 553 553  0.0000000 0.2204330}%
%
\special{pn 8}%
\special{ar 4160 992 552 552  0.2454152 6.2831853}%
\special{ar 4160 992 552 552  0.0000000 0.2208288}%
%
\special{pn 8}%
\special{sh 0}%
\special{pa 3748 545}%
\special{pa 3862 545}%
\special{pa 3862 651}%
\special{pa 3748 651}%
\special{pa 3748 545}%
\special{ip}%
%
\special{pn 8}%
\special{sh 0}%
\special{pa 3769 1382}%
\special{pa 3883 1382}%
\special{pa 3883 1489}%
\special{pa 3769 1489}%
\special{pa 3769 1382}%
\special{ip}%
%
\special{pn 8}%
\special{pa 3741 630}%
\special{pa 3760 604}%
\special{pa 3784 583}%
\special{pa 3807 563}%
\special{pa 3829 546}%
\special{pa 3833 538}%
\special{sp}%
%
\special{pn 8}%
\special{pa 3883 1368}%
\special{pa 3875 1398}%
\special{pa 3847 1408}%
\special{pa 3821 1423}%
\special{pa 3793 1436}%
\special{pa 3776 1463}%
\special{pa 3762 1475}%
\special{sp}%
%
\special{pn 8}%
\special{ar 5982 962 552 552  0.2454152 6.2831853}%
\special{ar 5982 962 552 552  0.0000000 0.2208288}%
\put(12.5000,-23.3000){\makebox(0,0)[lb]{$J_+$}}%
\put(37.9000,-24.6000){\makebox(0,0)[lb]{$J_-$}}%
\put(61.9000,-23.9000){\makebox(0,0)[lb]{$J_0$}}%
\put(29.1000,-17.7000){\makebox(0,0)[lb]{$-2$}}%
\put(44.3000,-17.2000){\makebox(0,0)[lb]{$-2$}}%
\put(64.2000,-16.4000){\makebox(0,0)[lb]{$-2$}}%
\put(3.7000,-20.7000){\makebox(0,0)[lb]{A Seifert hypersurface }}%
\put(8.2000,-23.1000){\makebox(0,0)[lb]{for}}%
\put(29.1000,-22.0000){\makebox(0,0)[lb]{A Seifert hypersurface }}%
\put(33.6000,-24.4000){\makebox(0,0)[lb]{for}}%
\put(53.1000,-21.6000){\makebox(0,0)[lb]{A Seifert hypersurface }}%
\put(57.6000,-24.0000){\makebox(0,0)[lb]{for}}%
\end{picture}%
\vskip5mm

\noindent 
$(J_+, J_-, J_0)$ is a twist-move triple in the 5-ball which is represented by the dotted circle in the following figure.

\vskip5mm
\includegraphics[width=16cm]{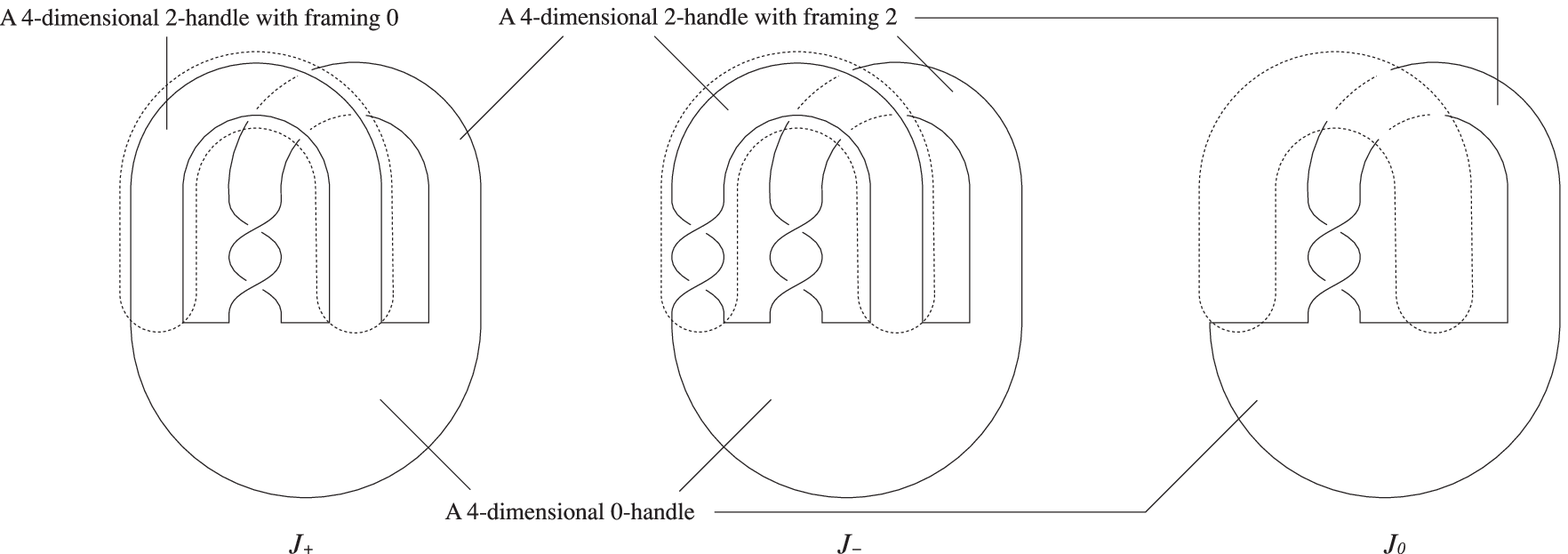}
\vskip5mm

\noindent 
By the calculus in  the example of  Theorem \ref{skein2}, 
the identity in Theorem \ref{pepper}   
holds for the triple  $(J_+, J_-, J_0)$.


\bigbreak


\begin{thm}\label{skeincro}    
Let $K_+$, $K_-$, $K_0$ be as in Theorem \ref{pepper}.  
There is a polynomial \newline
 $\Delta_{l+1+\nu} 
(\>K_* \otimes^\nu[2]\>)
\in \Q[t, t^{-1}]$ 
whose $\Q[t, t^{-1}]$-balanced class is 
the $(l+1+\nu)$-$\Q[t, t^{-1}]$-Alexander polynomial 
$A_{l+1+\nu} 
(\>K_*
\otimes^\nu[2]\>)$ 
$(\nu\in\N\cup\{0\}. \quad *=+,-,0. )$  
such that 
$$
\Delta_{l+1+\nu} 
(\>K_+\otimes^\nu[2]\>)
-
\Delta_{l+1+\nu} (\>K_-\otimes^\nu[2]\>)
= 
(t+(-1)^{l+1+\nu})\cdot \Delta_{l+1+\nu} (\>K_0\otimes^\nu[2]\>), 
$$

\noindent 
where 
$[2]$ denotes the empty knot  $[2]$.


\end{thm}

\noindent
\begin{thm}\label{proinertia}   
Let $K_+$ be a $(4l+1)$-dimensional spherical knot $(l\in\N\cup\{0\})$. 
Let $(K_+, K_-, K_0)$ be a twist-move-triple. 
Let $bP_{4l+2+4\mu}$ be the $bP$-subgroup $\subset \Theta^{4l+1+4\mu}$.
Suppose that $bP_{4l+2+4\mu}$ is not congruent to the trivial group. 
Then we have 

\bigbreak
\hskip1cm${\rm{Arf}} 
(\>K_+\otimes^\mu{\text{$($the Hopf link$)$}}\>) 
-
{\rm{Arf}}(\>K_-\otimes^\mu{\text{$($the Hopf link$)$}}\>)$

\noindent
\hskip1cm$=
\{|bP_{4k+2}\cap I(\>K_0\otimes^\mu{\text{$($the Hopf link$)$}}\>)|+1\}\hskip1mm {\rm mod 2},$
\bigbreak

\noindent 
where $
I(
\>K_0
\otimes^\mu
{\text{$($the Hopf link$)$}}\>)$ 
is the inertia group of an oriented smooth manifold 
which is orientation preserving diffeomorphic to 
$K_0\otimes^\mu{\text{$($the Hopf link$)$}}$    
and 
the symbol $|\quad|$ denotes the order of a group.
\end{thm}

Let $K_+$ and $K_-$ be  $n$-dimensional spherical knots $\subset S^{n+2}$. 
Let $K_0$ be an $n$-submanifold  $\subset S^{n+2}$. 
Let $(K_+, K_-, K_0)$ be related by a $(p, q)$-pass-move in $B^{n+2}$. 
Let $p+q=n+1$.
Let $p\neq q.$  
%
%
Recall that we have the following by \cite{Ogasa09}. (It is quoted in Theorem \ref{knots} in this paper.)

\smallbreak
There is a polynomial 
$\Delta_p(K_*) 
\in \Q[t, t^{-1}]$ 
whose $\Q[t, t^{-1}]$-balanced class is 
the $p$-$\Q[t, t^{-1}]$-Alexander polynomial 
$A_p(K_*)$  
for the submanifold
$K_*$  
$(*=+,-,0)$  
such that 
$$
\Delta_p(K_+)-\Delta_p(K_-)=(t-1)\cdot \Delta_p(K_0)
.$$


\noindent
\begin{thm}\label{skeinpass}   
Let $K_+$, $K_-$, $K_0$ be as above. 
There is a polynomial    
 $\Delta_{p+\nu} (\>K_*\otimes^\nu[2]\>)$\newline$
\in \Q[t, t^{-1}]$ 
whose $\Q[t, t^{-1}]$-balanced class  is 
the $(p+\nu)$-$\Q[t, t^{-1}]$-Alexander polynomial 
$A_{p+\nu} (\>K_*\otimes^\nu[2]\>)
\quad (*=+,-,0, \quad \nu\in\N\cup\{0\})$  
such that 
$$
\Delta_{p+\nu} 
(\>K_+
\otimes^\nu[2]\>)
-
\Delta_{p+\nu} 
(\>K_-
\otimes^\nu[2]\>)
=
(t+(-1)^{\nu+1})\cdot 
\Delta_{p+\nu} 
(\>K_0
\otimes^\nu[2]\>),$$

\noindent
where 
$[2]$ denotes the empty knot  $[2]$. 
\end{thm}

\section{A remark on the $\Z[t, t^{-1}]$-case  }\label{specialX}    
Some of our results on polynomial invariants could be extended to the case where 
the word, `$\Q[t, t^{-1}]$-balanced class,' is replaced with 
the word, `($\Z[t, t^{-1}]$-balanced class of) an element of $\Z[t, t^{-1}]$'. 
However we must take care of the following Theorem \ref{noX}.   

Two polynomials $f(t), g(t)\in\Z[t, t^{-1}]$ are said to be 
{\it $\Z[t, t^{-1}]$-balanced } 
if there is an integer $n$ 
such that 
$f(t)=\pm t^n\cdot g(t)$.

\begin{thm}\label{noX}   
There is a smooth oriented codimension two closed submanifold $K$\newline$
\subset S^{n+2}$ 
with the following property. 

\smallbreak
\noindent $(1)$ 
Any Seifert hypersurface of $K$  
satisfies the condition that 
the $p$-th betti number is equal to the $(n+1-p)$-th betti number 
$(1\leqq p\leqq n).$ 

\noindent $(2)$  
There are Seifert hypersurfaces $V$ and $W$ for $K$  such that 
the $\Z[t, t^{-1}]$-balanced class of 
the determinant of a $p$-Alexander matrix of $V$ and 
that of $W$ are different for a nonnegative integer $p$ 
  even if 
the linear map defined by 
a $(p-1)$-Alexander matrix  associated with $V$ (resp. $W$) is injective. 
\end{thm}

\noindent{\bf Note.} 
The $\Z[t, t^{-1}]$-balanced class of 
the determinant of a $p$-Alexander matrix of $V$ (resp. $W$) 
is same as that of any  $p$-Alexander matrix of $V$ (resp. $W$) . 
If the linear map defined by 
a $(p-1)$-Alexander matrix  associated with $V$ (resp. $W$) is injective, 
then 
that defined by 
any $(p-1)$-Alexander matrix  associated with $V$ (resp. $W$) is injective. 
This is because they do not depend on the choice of the basis of the homology groups 
of $V$ (resp. $W$).

\smallbreak
\noindent{\bf Proof.} 
Let $K$ be a closed oriented smooth 3-dimensional submanifold $\subset S^5$ 
such that the diffeomorphism type of $K$ is represented 
by the following framed link: 
Take the $(a, a)$-torus link $\subset S^3$ ($a\in\N-\{1\}$). 
The framing of each component is zero.  
Then $K$ is diffeomorphic to a homology sphere and, therefore, 
  any Seifert hypersurface for $K$ satisfies the condition (1) in Theorem \ref{noX}. 
%
%
The following figure is the $a=2$ case of the framed link.

\hskip30mm
\includegraphics[width=5cm]{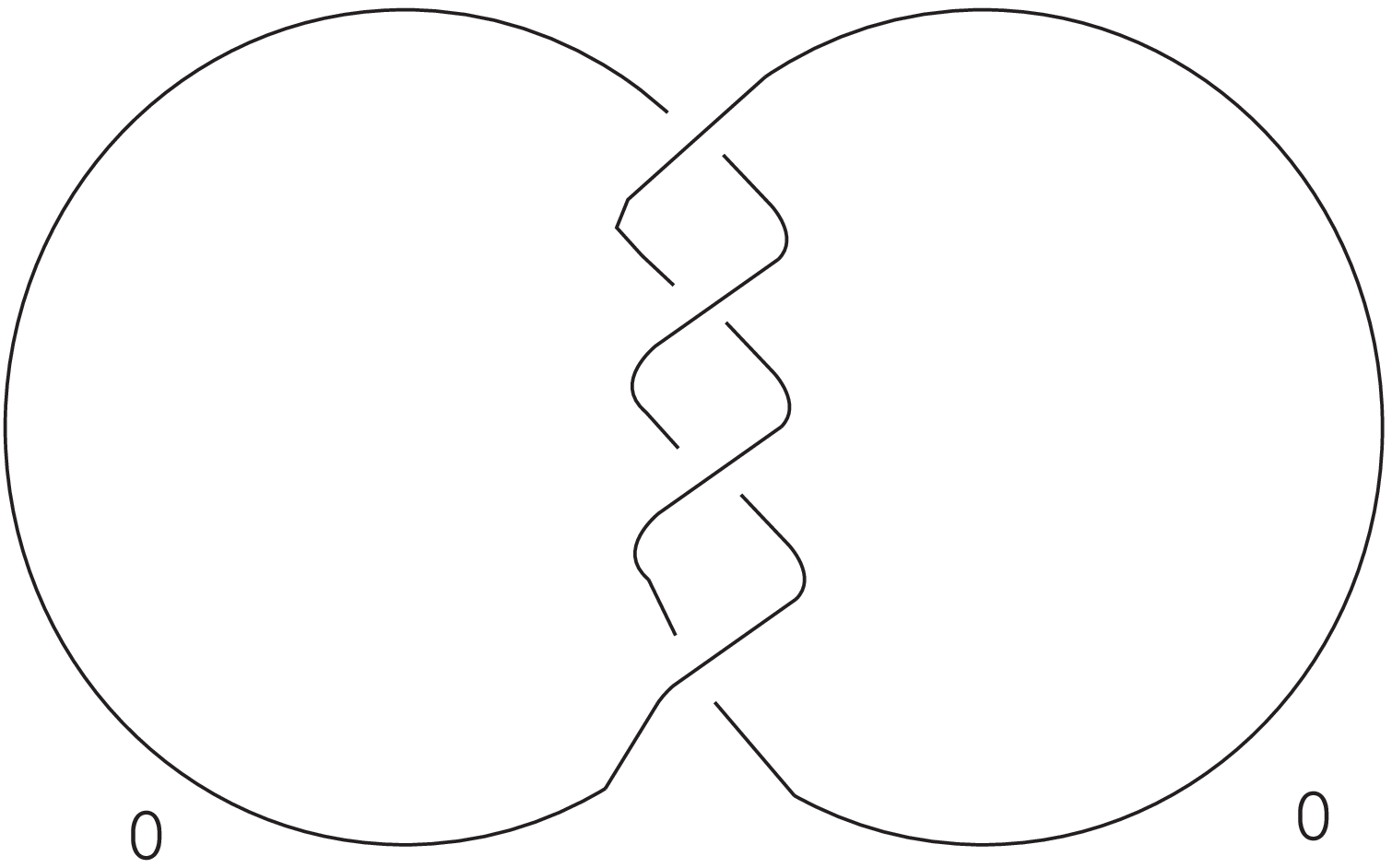}

\bigbreak
We make two kinds of Seifert hypersurfaces $V$, $W$ for $K$ as follows.

\bigbreak
{\bf The first.}  
Regard $\R^5=\R^4\x\{t\in\R\}$. 
Regard the above framed link 
which represents $K$ 
as a 4-manifold, too. 
Of course this 4-manifold has a handle decomposition 

\noindent 
(a 4-dimensional 0-handle)$\cup$(a 4-dimensional 2-handle)
$\cup$(a 4-dimensional 2-handle), 

\noindent
which is defined by the framed link representation.

Suppose that the diffeomorphism type of a Seifert hypersurface $V$ is 
this 4-manifold.  

Suppose that $V$ in $\R^5$ satisfies the following. 

\smallbreak
\noindent 
(1) The 4-dimensional 0-handle is embedded in $\R^4\x\{t=0\}$. 

\noindent 
(2) One of the 4-dimensional 2-handles is embedded in $\R^4\x\{t=0\}$, call it $h^2$. 

\noindent 
(3) The other of the 4-dimensional 2-handles is embedded in $\R^4\x\{t\leqq0\}$. 
Only the attached part is embedded in $\R^4\x\{t=0\}$. 
We can do this because the framing is zero. 
\smallbreak

Thus we obtain a Seifert hypersurface $V$ for $K$. 


Then a positive 2-Seifert matrix associated with $V$ is 
$\begin{pmatrix}
0&a\\
0&0\\
\end{pmatrix}$. 
The negative 2-Seifert matrix associated with the positive 2-Seifert matrix 
is 
$\begin{pmatrix}
0&0\\
-a&0\\
\end{pmatrix}$. 
Hence the 2-Alexander matrix associated with these two matrices is 
$\begin{pmatrix}
0&at\\
a&0\\
\end{pmatrix}$. 
Its determinant is $-a^2t$.

\noindent
Note that the linear map defined by a 1-Alexander matrix associated with $V$ 
is injective. 

\bigbreak 
{\bf The second.}    
Take the above Seifert hypersurface $V$. 
Suppose that a 5-dimensional 3-handle $k^3$ is embedded in $\R^4\x\{t\geqq0\}$.  
Attach  $k^3$ along the 2-sphere embedded in $V$ which includes the core of 
the above $h^2$. Suppose that only the attach part is embedded in  $\R^4\x\{t=0\}$. 
By this surgery by $k^3$, $V$ is changed into another Seifert  hypersurface $W$ for $K$. 
Then the framed link representation of $W$ is as follows:  
Take the $(a, a)$-torus link $\subset S^3$. 
The framing of one component is zero.  
The other component is the dot circle  
(see \cite{Kirby} for the dot circle). 
Then  $W$ is a rational homology ball.

The following figure is the $a=2$ case of $W$. 

\hskip20mm
\includegraphics[width=6cm]{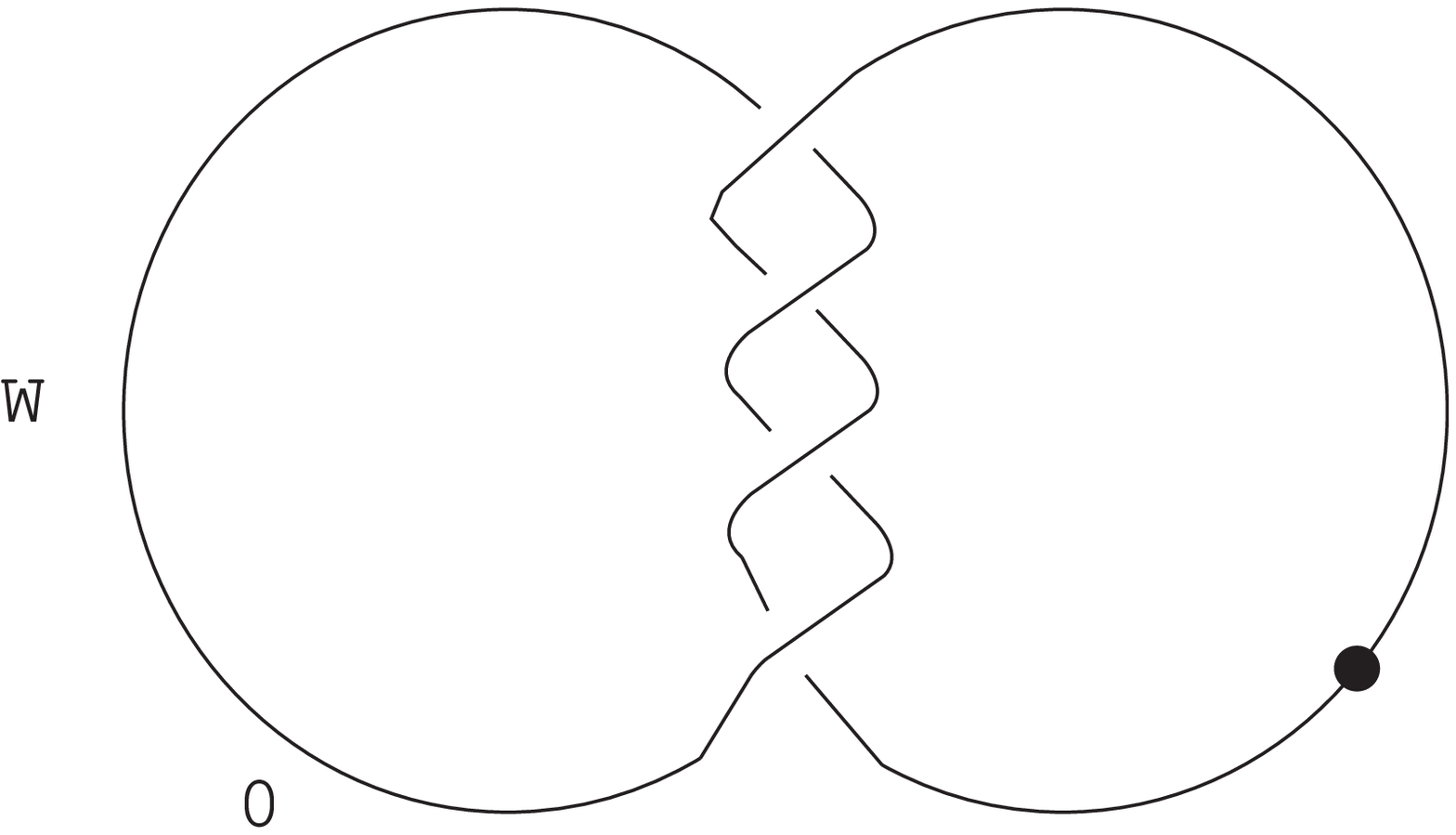}
\bigbreak

Hence the following holds: 
the positive 2-Seifert matrix associated with $W$ is `empty'. 
The negative 2-Seifert matrix associated with $W$ is `empty'.
Hence the 2-Alexander matrix associated with $W$ is `empty'.   
Note that 
the $\Z[t, t^{-1}]$-balanced class of the determinant of the 2-Alexander matrix `empty' is 
that of 1.
Note that the linear map defined by a 1-Alexander matrix associated with $W$ 
is injective. 
%
%
Note that 
the $\Z[t, t^{-1}]$-balanced class of the determinant of the 2-Alexander matrix `empty' is 
NOT that of $a^2t$.  (Recall that we define $a\in\N-\{1\}$.)

This completes the proof.  \qed

\section{Proof of Theorems in \S \ref{thecr}}\label{procr}    
We use the following proposition. 

\begin{prop}\label{KLS}     
Let $K$ be an $n$-dimensional closed oriented submanifold $\subset S^{n+2}$.  
Then we have the following. 
\smallbreak
\noindent
$(1)$ \hskip3cm$S_{p+1}(K\otimes[2])=(-1)^{(n-p)}S_p(K)\otimes S_0([2]).$

\smallbreak

\noindent
$(2)$ \hskip3cm$N_{p+1}(K\otimes[2])=(-1)^{(n-p)}N_p(K)\otimes N_0([2]).$
\smallbreak

\noindent
$(3)$ \hskip6cm$S_0([2])=(1).$

\smallbreak

\noindent
$(4)$  \hskip6cm$N_0([2])=(-1).$

\smallbreak
 
\noindent
$(5)$   
\hskip6cm $S_{\nu-1}(\otimes^\nu[2])=(-1)^{\frac{(\nu-1)\nu}{2}}.$  

\hskip8mm
Note that 
$\otimes^\nu[2]$
is a $(2\nu-3)$-submanifold $\subset S^{2\nu-1}$.

\smallbreak

\noindent
$(6)$  \hskip3cm$S_{p+\nu}(K\otimes^\nu[2])=
(-1)^{(n-p)\nu+\frac{\nu(\nu-1)}{2}}\quad S_p(K).$ 

\smallbreak

\noindent
$(7)$  
\hskip3cm$S_{p+2\mu}(K\otimes^\mu\text{$($the Hopf link$)$})=(-1)^\mu S_p(K).$  

\hskip8mm
$($Let $\nu=2\mu$ in $(6).$ We obtain $(7)$.$)$

\noindent
$(8)$ 
Let $\alpha$, $\beta$ be disjoint cycles of dimension $p$ and $q$ in $S^{p+q+1}$. 
$$\mathrm{lk}(\alpha, \beta)=(-1)^{pq+1}\mathrm{lk}(\beta,\alpha).$$ 

\end{prop}

By \S6 of \cite{KauffmanNeumann} we have (1)-(7). 
By  page 541 of \cite{Levinepol} we have (8).

\bigbreak
\noindent{\bf Proof of Theorem \ref{Chicago}. }     
Take 1-links $K_+, K_-\subset S^3=\partial B^4\subset B^4.$  
We can suppose that the 1-links, $K_+, K_-$,  differ by only one crossing-change in 
a 3-ball $A$ trivially embedded in $S^3$ 
as shown below.  

\bigbreak
\unitlength 0.1in
\begin{picture}(28.22,12.30)(8.10,-20.40)
%
\special{pn 8}%
\special{ar 1356 1356 546 546  0.8902751 6.2831853}%
\special{ar 1356 1356 546 546  0.0000000 0.8502422}%
%
\special{pn 8}%
\special{pa 946 1716}%
\special{pa 1736 986}%
\special{fp}%
\special{sh 1}%
\special{pa 1736 986}%
\special{pa 1673 1017}%
\special{pa 1697 1022}%
\special{pa 1701 1046}%
\special{pa 1736 986}%
\special{fp}%
%
\special{pn 8}%
\special{pa 1266 1276}%
\special{pa 976 976}%
\special{fp}%
\special{sh 1}%
\special{pa 976 976}%
\special{pa 1008 1038}%
\special{pa 1013 1014}%
\special{pa 1037 1010}%
\special{pa 976 976}%
\special{fp}%
%
\special{pn 8}%
\special{pa 1416 1416}%
\special{pa 1726 1746}%
\special{fp}%
%
\special{pn 8}%
\special{ar 3086 1396 546 546  0.8902751 6.2831853}%
\special{ar 3086 1396 546 546  0.0000000 0.8502422}%
%
\special{pn 8}%
\special{pa 2996 1316}%
\special{pa 2706 1016}%
\special{fp}%
\special{sh 1}%
\special{pa 2706 1016}%
\special{pa 2738 1078}%
\special{pa 2743 1054}%
\special{pa 2767 1050}%
\special{pa 2706 1016}%
\special{fp}%
%
\special{pn 8}%
\special{pa 3146 1456}%
\special{pa 3456 1786}%
\special{fp}%
%
\special{pn 8}%
\special{pa 3210 1520}%
\special{pa 2990 1320}%
\special{fp}%
\special{pa 3240 1550}%
\special{pa 3100 1420}%
\special{fp}%
%
\special{pn 8}%
\special{pa 3140 1330}%
\special{pa 3480 1010}%
\special{fp}%
\special{sh 1}%
\special{pa 3480 1010}%
\special{pa 3418 1041}%
\special{pa 3441 1047}%
\special{pa 3445 1070}%
\special{pa 3480 1010}%
\special{fp}%
%
\special{pn 8}%
\special{pa 3010 1470}%
\special{pa 2690 1770}%
\special{fp}%
\put(11.0000,-21.6000){\makebox(0,0)[lb]{$A\cap K_+$}}%
\put(28.7000,-22.1000){\makebox(0,0)[lb]{$A\cap K_-$}}%
\end{picture}%
\vskip10mm

Take a Seifert surface $V_+$ (resp. $V_-$) for the 1-link $K_+$ (resp. $K_-$) 
such that the submanifolds, $V_+, V_-$, differ only in the 3-ball $A\subset S^3$ as shown below.

\bigbreak
\unitlength 0.1in
\begin{picture}(28.22,12.20)(8.10,-20.30)
%
\special{pn 8}%
\special{ar 1356 1356 546 546  0.8902751 6.2831853}%
\special{ar 1356 1356 546 546  0.0000000 0.8502422}%
%
\special{pn 8}%
\special{pa 946 1716}%
\special{pa 1736 986}%
\special{fp}%
\special{sh 1}%
\special{pa 1736 986}%
\special{pa 1673 1017}%
\special{pa 1697 1022}%
\special{pa 1701 1046}%
\special{pa 1736 986}%
\special{fp}%
%
\special{pn 8}%
\special{pa 1266 1276}%
\special{pa 976 976}%
\special{fp}%
\special{sh 1}%
\special{pa 976 976}%
\special{pa 1008 1038}%
\special{pa 1013 1014}%
\special{pa 1037 1010}%
\special{pa 976 976}%
\special{fp}%
%
\special{pn 8}%
\special{pa 1416 1416}%
\special{pa 1726 1746}%
\special{fp}%
%
\special{pn 8}%
\special{ar 3086 1396 546 546  0.8902751 6.2831853}%
\special{ar 3086 1396 546 546  0.0000000 0.8502422}%
%
\special{pn 8}%
\special{pa 2996 1316}%
\special{pa 2706 1016}%
\special{fp}%
\special{sh 1}%
\special{pa 2706 1016}%
\special{pa 2738 1078}%
\special{pa 2743 1054}%
\special{pa 2767 1050}%
\special{pa 2706 1016}%
\special{fp}%
%
\special{pn 8}%
\special{pa 3146 1456}%
\special{pa 3456 1786}%
\special{fp}%
%
\special{pn 8}%
\special{pa 3210 1520}%
\special{pa 2990 1320}%
\special{fp}%
\special{pa 3240 1550}%
\special{pa 3100 1420}%
\special{fp}%
%
\special{pn 8}%
\special{pa 3140 1330}%
\special{pa 3480 1010}%
\special{fp}%
\special{sh 1}%
\special{pa 3480 1010}%
\special{pa 3418 1041}%
\special{pa 3441 1047}%
\special{pa 3445 1070}%
\special{pa 3480 1010}%
\special{fp}%
%
\special{pn 8}%
\special{pa 3010 1470}%
\special{pa 2690 1770}%
\special{fp}%
\put(10.6000,-21.5000){\makebox(0,0)[lb]{$A\cap V_+$}}%
\put(28.3000,-22.0000){\makebox(0,0)[lb]{$A\cap V_-$}}%
%
\special{pn 8}%
\special{pa 1160 870}%
\special{pa 1160 1100}%
\special{dt 0.045}%
\special{pa 1160 1100}%
\special{pa 1160 1099}%
\special{dt 0.045}%
\special{pa 1320 840}%
\special{pa 1310 1230}%
\special{dt 0.045}%
\special{pa 1310 1230}%
\special{pa 1310 1229}%
\special{dt 0.045}%
\special{pa 1510 870}%
\special{pa 1510 1090}%
\special{dt 0.045}%
\special{pa 1510 1090}%
\special{pa 1510 1089}%
\special{dt 0.045}%
\special{pa 1440 840}%
\special{pa 1420 1110}%
\special{dt 0.045}%
\special{pa 1420 1110}%
\special{pa 1420 1109}%
\special{dt 0.045}%
\special{pa 1620 920}%
\special{pa 1620 1000}%
\special{dt 0.045}%
\special{pa 1620 1000}%
\special{pa 1620 999}%
\special{dt 0.045}%
\special{pa 1110 1620}%
\special{pa 1110 1790}%
\special{dt 0.045}%
\special{pa 1110 1790}%
\special{pa 1110 1789}%
\special{dt 0.045}%
\special{pa 1200 1540}%
\special{pa 1190 1830}%
\special{dt 0.045}%
\special{pa 1190 1830}%
\special{pa 1190 1829}%
\special{dt 0.045}%
\special{pa 1340 1450}%
\special{pa 1350 1860}%
\special{dt 0.045}%
\special{pa 1350 1860}%
\special{pa 1350 1859}%
\special{dt 0.045}%
\special{pa 1470 1560}%
\special{pa 1480 1810}%
\special{dt 0.045}%
\special{pa 1480 1810}%
\special{pa 1480 1809}%
\special{dt 0.045}%
\special{pa 1600 1650}%
\special{pa 1610 1770}%
\special{dt 0.045}%
\special{pa 1610 1770}%
\special{pa 1610 1769}%
\special{dt 0.045}%
%
\special{pn 8}%
\special{pa 2900 930}%
\special{pa 2900 1160}%
\special{dt 0.045}%
\special{pa 2900 1160}%
\special{pa 2900 1159}%
\special{dt 0.045}%
\special{pa 3060 900}%
\special{pa 3050 1290}%
\special{dt 0.045}%
\special{pa 3050 1290}%
\special{pa 3050 1289}%
\special{dt 0.045}%
\special{pa 3250 930}%
\special{pa 3250 1150}%
\special{dt 0.045}%
\special{pa 3250 1150}%
\special{pa 3250 1149}%
\special{dt 0.045}%
\special{pa 3180 900}%
\special{pa 3160 1170}%
\special{dt 0.045}%
\special{pa 3160 1170}%
\special{pa 3160 1169}%
\special{dt 0.045}%
\special{pa 3360 980}%
\special{pa 3360 1060}%
\special{dt 0.045}%
\special{pa 3360 1060}%
\special{pa 3360 1059}%
\special{dt 0.045}%
\special{pa 2850 1680}%
\special{pa 2850 1850}%
\special{dt 0.045}%
\special{pa 2850 1850}%
\special{pa 2850 1849}%
\special{dt 0.045}%
\special{pa 2940 1600}%
\special{pa 2930 1890}%
\special{dt 0.045}%
\special{pa 2930 1890}%
\special{pa 2930 1889}%
\special{dt 0.045}%
\special{pa 3080 1510}%
\special{pa 3090 1920}%
\special{dt 0.045}%
\special{pa 3090 1920}%
\special{pa 3090 1919}%
\special{dt 0.045}%
\special{pa 3210 1620}%
\special{pa 3220 1870}%
\special{dt 0.045}%
\special{pa 3220 1870}%
\special{pa 3220 1869}%
\special{dt 0.045}%
\special{pa 3340 1710}%
\special{pa 3350 1830}%
\special{dt 0.045}%
\special{pa 3350 1830}%
\special{pa 3350 1829}%
\special{dt 0.045}%
\end{picture}%
\vskip10mm

Take the collar neighborhood  $A\x[0,1]$ of $B^4$. 
Note that $A\x[0,1]$ and \newline
$\overline{\hskip1mm B^4-(A\x[0,1])\hskip1mm }$ 
are diffeomorphic to the 4-ball, 
where $\overline{Q}$ denotes the closure of $Q$ in $B^4$ if $Q\subset B^4$.    
Note that the intersection of the two 4-balls is 
$(\partial A\x[0,1])\cup (A\x\{1\})$ and is diffeomorphic to the 3-ball.

\smallbreak

Push $V_+$ (resp. $V_-$) into $B^4$, fixing $\partial V_+=K_+$ (resp. $\partial V_-=K_-$), 
call the submanifold $\subset B^4$,   $V'_+$ (resp. $V'_-$).   
Suppose that the submanifolds, $V'_+, V'_-$, differ only in the 4-ball $A\x[0,1]\subset B^4$. 
Suppose that $V'_*\cap (A\x[0,1])$ is 
$(\hskip1mm((\partial V_*)\cap A)\x[0,\frac{1}{2}]\hskip1mm) 
\cup 
(\hskip1mm(V_*\cap A)\x\{\frac{1}{2}\}\hskip1mm)$
($*=+,-$).

\bigbreak
Let $B^2=\{(x, y)\in\R^2|x^2+y^2\leqq1\}$. 
Then there are smooth maps 
%
$f_+:B^4\to B^2$  and 
%
%
$f_-:B^4\to B^2$
with the following properties.  
(Reason: Use the Thom-Pontrjagin construction.)

\smallbreak
\noindent 
(i) The submanifold $f_+^{-1}((0,0))\subset B^4$ is $V'_+$. 

\noindent \hskip6mm 
The submanifold $f_-^{-1}((0,0))\subset B^4$ is $V'_-$. 

\noindent
(ii)  
$f_+$ and $f_-$ coincide on  
$\overline{\hskip1mm B^4-(A\x[0,1])\hskip1mm}$ 

\smallbreak


The knot product 
$K_*\otimes($the empty knot $[2])$   
is defined as follows ($*=+,-$).   
See \S\ref{Introduction}, \S\ref{knotproduct} and \cite{Kauffman, KauffmanNeumann} 
for knot products. 
Take a smooth map $g:B^2\to B^2$ such that \newline
$(r\text{cos}\theta, r\text{sin}\theta) 
\mapsto 
(r\text{cos}2\theta, r\text{sin}2\theta)$, 
where we use the standard polar coordinate. 
Recall that $g|_{\partial B^2=S^1}:S^1\to S^1$ is the empty knot $[2]$.  

\noindent
Let $M_*=
\{(x, y)\in B^4\times B^2|f_*(x)-g(y)=(0,0)\in B^2\}$.

\noindent 
The 3-dimensional closed oriented submanifold 
$\partial M_*\subset \partial(B^4\times B^2)$  
is the knot product 
$K_*\otimes [2]$ in the standard 5-sphere.   
Note that $\partial(B^4\times B^2)$ is the standard 5-sphere. 

\bigbreak
$\partial (\hskip1mm(A\x[0,1])\x B^2\hskip1mm)
\cap\partial (\hskip1mmB^4\x B^2\hskip1mm)$ 
is the 5-ball, call it $\check A$. 

\noindent
$\partial 
(\hskip1mm \overline{\hskip1mm B^4-(A\x[0,1])\hskip1mm}\x B^2\hskip1mm )
\cap\partial (B^4\x B^2)$  
is also the 5-ball. 
It is 
$(\partial (B^4\x B^2))- $Int$\check A$.  

\bigbreak
By \cite{Kauffman, KauffmanNeumann} we have the following. 

\noindent(1)  
$\partial M_*$ is the double branched covering space of $\partial B^4$ along $K_*$. 

\noindent(2) 
A Seifert hypersurface for $\partial M_*$ is the double branched covering space of $B^4$ 
along $V'_*$.

\bigbreak
By (ii) several lines above here, 
there is a diffeomorphism map \newline$\alpha:\partial(B^4\x B^2)\to\partial(B^4\x B^2)$ 
with the following properties: 

\noindent 
(1) $\alpha\vert_{(\partial (B^4\x B^2)) - \rm{Int}\check A}$ is the identity map.

\noindent 
(2) $\alpha\vert_{(\partial (B^4\x B^2)) - \rm{Int}\check A}
(\hskip1mm(\partial M_+)-\rm{Int}\check A\hskip1mm) 
= (\partial M_-)-\rm{Int}\check A.$

\noindent 
(3) $\alpha\vert_{(\partial (B^4\x B^2))- \rm{Int}\check A}
(\hskip1mm(\text{the Seifert hypersurface for } \partial M_+)-\rm{Int}\check A\hskip1mm)
\newline
 = 
(\text{the Seifert hypersurface for } \partial M_-)-\rm{Int}\check A.$

\bigbreak
By \cite{Kauffman, KauffmanNeumann} we have the following. 

\noindent (1)  
$\check A\cap \partial M_*$ is the double branched covering space of 
$A$ along $A\cap K_*$  
and is  $S^1\x D^2$. 
Note that $A\cap K_*$ is a set of two arcs.

\noindent (2) 
$\check A\cap $ (the Seifert hypersurface for $\partial M_*$) is the double branched covering space of 
$A\x[0,1]$ along $(A\x[0,1])\cap V'_*$ 
and is $D^2\x D^2$.  
Note that 
$(A\x[0,1])\cap V'_*$ is  a rectangle.

\smallbreak

Hence we have the following: 
The intersection of \newline
$\check A\cap $ (the Seifert hypersurface for $\partial M_*$),  
which is $D^2\x D^2$, 
 and the standard 4-sphere $\partial\check A$ 
is $S^1\x D^2$, which is $\check A\cap \partial M_*$.  
Thus we can regard this $D^2\x D^2$ as 
a 4-dimensional 2-handle embedded in the standard 5-ball $\check A$ 
which is attached to the standard 4-sphere $\partial\check A$.

\bigbreak
Since $K_+$ and $K_-$ differ by only one crossing-change, we can suppose that 
there is a Seifert matrix $X_*=(x^*_{i,j})$ for $K_*$ with the following property ($*=+.-$).

\noindent
$
\left\{
\begin{array}{ll}
x^+_{11}=1, &x^-_{11}=0 \\
x^+_{ij}=x^-_{ij}&\text{if $(i,j)\neq(1,1). \quad i,j\leqq\nu.  \quad i,j\in\N-\{0\}.$}\\
\end{array}
\right.
$

Here, 
$x^*_{11}$ is derived from $A\cap V_*$  $(*=+, -)$.  

By \cite{Kauffman, KauffmanNeumann} and Proposition \ref{KLS}, 
we can suppose that 
the $\nu\x\nu$-matrix $-X_*=(-x^*_{i,j})$ 
is a 2-Seifert matrix for $K_*\otimes[2]$.  
Here, 
$-x^*_{11}$ is derived from \newline 
$\check A\cap$ (the Seifert hypersurface for $\partial M_*$) $(*=+, -)$.

Recall the definition of twist-moves in \S\ref{localn}. 
By the above two paragraphs,  
$K_+\otimes[2]$ and $K_-\otimes[2]$ differ by one twist-move. 

This completes the proof. \qed

\bigbreak
\noindent{\bf Note.} The above  proof is the $m=0, \nu=1$ case of the Proof of Theorem \ref{Tokyo}.

\bigbreak
\noindent{\bf Proof of Theorem \ref{Illinois}.}   
This proof is the $m=0$ case of Proof of Theorem \ref{Tokyo}.  \qed



\bigbreak
\noindent{\bf Proof of Theorem \ref{Tokyo}.}   
We first prove the $\nu=1$ case.

Take $(2m+1)$-dimensional oriented closed submanifolds 
$K_+, K_-\subset S^{2m+3}=\partial B^{2m+4}$\newline$\subset B^{2m+4}.$  
The $(2m+1)$-dimensional oriented closed submanifolds, $K_+, K_-$,  differ by only one 
twist-move in 
a $(2m+3)$-ball $A$ trivially embedded in $S^{2m+3}$.  
See the left two figures 
in the upper half of Figure \ref{localn}.5.

Take a Seifert hypersurface $V_+$ (resp. $V_-$) for the $(2m+1)$-dimensional oriented closed submanifold $K_+$ (resp. $K_-$) 
such that the submanifolds, $V_+, V_-$, differ by only one 
twist-move in the $(2m+3)$-ball 
$A\subset S^{2m+3}$  
and that 
$V_*\cap A$ is the $(2m+2)$-dimensional $(m+1)$-handle 
$h_*(*=+,-)$. See the definition of twist-moves in \S\ref{localn}. 
%
%





\smallbreak
Take the collar neighborhood  $A\x[0,1]$ of $B^{2m+4}$. 
Note that $A\x[0,1]$ and \newline
$\overline{\hskip1mm B^{2m+4}-(A\x[0,1])\hskip1mm}$ 
are diffeomorphic to the $(2m+4)$-ball, 
where $\overline{Q}$ denotes the closure of $Q$ in $B^{2m+4}$ if $Q\subset B^{2m+4}$.    
Note that the intersection of the two $(2m+4)$-balls is 
$(\partial A\x[0,1])\cup (A\x\{1\})$ and is diffeomorphic to the $(2m+3)$-ball.

\smallbreak

Push $V_+$ (resp. $V_-$) into $B^{2m+4}$, fixing $\partial V_+=K_+$ (resp. $\partial V_-=K_-$), 
call the submanifold $\subset B^{2m+4}$,   $V'_+$ (resp. $V'_-$).   
Suppose that the submanifolds, $V'_+, V'_-$, differ only in 
the $(2m+4)$-ball $A\x[0,1]\subset B^{2m+4}$. 
Suppose that $V'_*\cap (A\x[0,1])$ is 
$(\hskip1mm((\partial V_*)\cap A)\x[0,\frac{1}{2}]\hskip1mm) 
\cup 
(\hskip1mm(V_*\cap A)\x\{\frac{1}{2}\}\hskip1mm)$ ($*=+,-$).

\bigbreak
Let $B^2=\{(x, y)\in\R^2|x^2+y^2\leqq1\}$. 
Then there are smooth maps 
%
$f_+:B^{2m+4}\to B^2$  and 
$f_-:B^{2m+4}\to B^2$
with the following properties.  
(Reason: Use the Thom-Pontrjagin construction.)

\smallbreak
\noindent 
(i) The submanifold $f_+^{-1}((0,0))\subset B^{2m+4}$ is $V'_+$. 

\noindent \hskip6mm 
The submanifold $f_-^{-1}((0,0))\subset B^{2m+4}$ is $V'_-$. 

\noindent
(ii)  
$f_+$ and $f_-$ coincide on  
$\overline{\hskip1mm B^{2m+4}-(A\x[0,1])\hskip1mm}$ 

\smallbreak


The knot product 
$K_*\otimes($the empty knot $[2])$   
is defined as follows ($*=+,-$).   
See \S\ref{Introduction}, \S\ref{knotproduct} and \cite{Kauffman, KauffmanNeumann} 
for knot products. 
Take a smooth map $g:B^2\to B^2$ such that \newline
$(r\text{cos}\theta, r\text{sin}\theta) 
\mapsto 
(r\text{cos}2\theta, r\text{sin}2\theta)$, 
where we use the standard polar coordinate. 
Recall that $g|_{\partial B^2=S^1}:S^1\to S^1$ is the empty knot $[2]$.  

\noindent
Let $M_*=
\{(x, y)\in B^{2m+4}\times B^2|f_*(x)-g(y)=(0,0)\in B^2\}$.

\noindent 
The $(2m+3)$-dimensional closed oriented submanifold 
$\partial M_*\subset \partial(B^{2m+4}\times B^2)$  
is the knot product 
$K_*\otimes [2]$ in the standard $(2m+5)$-sphere.   
Note that $\partial(B^{2m+4}\times B^2)$ is the standard $(2m+5)$-sphere.

\bigbreak
$\partial (\hskip1mm(A\x[0,1])\x B^2\hskip1mm)
\cap\partial (\hskip1mm B^{2m+4}\x B^2\hskip1mm)$ 
is the $(2m+5)$-ball, call it $\check A$. 

\noindent
$\partial (\hskip1mm \overline{\hskip1mm B^{2m+4}-(A\x[0,1])\hskip1mm }\x B^2\hskip1mm )
\cap\partial (B^{2m+4}\x B^2)$  
is also the $(2m+5)$-ball. 
It is \newline
$(\partial (B^{2m+4}\x B^2))- $Int$\check A$.  

\bigbreak
By \cite{Kauffman, KauffmanNeumann} we have the following. 

\noindent(1)  
$\partial M_*$ is the double branched covering space of $\partial B^{2m+4}$ along $K_*$. 

\noindent(2) 
A Seifert hypersurface for $\partial M_*$ is the double branched covering space 
of $B^{2m+4}$ along $V'_*$.

\bigbreak
By (ii) several lines above here, 
there is a diffeomorphism map \newline  
$\alpha:\partial(B^{2m+4}\x B^2)\to\partial(B^{2m+4}\x B^2)$ 
with the following properties: 

\noindent 
(1) $\alpha\vert_{(\partial (B^{2m+4}\x B^2)) - \rm{Int}\check A}$ is the identity map.

\noindent 
(2) $\alpha\vert_{(\partial (B^{2m+4}\x B^2)) - \rm{Int}\check A}
(\hskip1mm(\partial M_+)-\rm{Int}\check A\hskip1mm) 
= (\partial M_-)-\rm{Int}\check A.$

\noindent 
(3) $\alpha\vert_{(\partial (B^{2m+4}\x B^2))- \rm{Int}\check A}
(\hskip1mm(\text{the Seifert hypersurface for } \partial M_+)-\rm{Int}\check A\hskip1mm)
\newline
=
(\text{the Seifert hypersurface for } \partial M_-)-\rm{Int}\check A.$

\bigbreak
By \cite{Kauffman, KauffmanNeumann} we have the following. 

\noindent (1)  
$\check A\cap \partial M_*$ is the double branched covering space of 
$A$ along $A\cap K_*$  
and is  $S^{m+1}\x D^{m+2}$. 
Note that $A\cap K_*$ is $S^m\x D^{m+1}$.

\noindent (2) 
$\check A\cap $ (the Seifert hypersurface for $\partial M_*$) 
is the double branched covering space of 
$A\x[0,1]$ along $(A\x[0,1])\cap V'_*$ 
and is $D^{m+2}\x D^{m+2}$.  
Note that 
$(A\x[0,1])\cap V'_*$ is $D^{m+1}\x D^{m+1}$.

\smallbreak

Hence we have the following: 
The intersection of \newline
$\check A\cap $ (the Seifert hypersurface for $\partial M_*$),  
which is $D^{m+2}\x D^{m+2}$, 
 and the standard $(m+4)$-sphere $\partial\check A$ 
is $S^{m+1}\x D^{m+2}$, which is $\check A\cap \partial M_*$.  
Thus we can regard this $D^{m+2}\x D^{m+2}$ as 
a $(2m+4)$-dimensional $(m+2)$-handle embedded 
in the standard $(2m+5)$-ball $\check A$ 
which is attached to the standard $(2m+4)$-sphere $\partial\check A$.

\bigbreak
Since $K_+$ and $K_-$ differ by only one twist-move,   
we can suppose that there is a $(m+1)$-Seifert matrix 
$X_*=(x^*_{i,j})$ for $K_*$ with the following property ($*=+.-$).

\noindent
$
\left\{
\begin{array}{ll}
x^+_{11}=1, &x^-_{11}=0 \\
x^+_{ij}=x^-_{ij}&\text{if $(i,j)\neq(1,1). \quad i,j\leqq\nu.  \quad i,j\in\N-\{0\}.$}\\
\end{array}
\right.
$

Here, 
$x^*_{11}$ is derived from $A\cap V_*$  $(*=+, -)$.  

By \cite{Kauffman, KauffmanNeumann} and Proposition \ref{KLS}, 
we can suppose that 
the $\nu\x\nu$-matrix $-X_*=(-x^*_{i,j})$ 
is a $(m+2)$-Seifert matrix for $K_*\otimes[2]$.  
Here, 
$-x^*_{11}$ is derived from \newline 
$\check A\cap$ (the Seifert hypersurface for $\partial M_*$) $(*=+, -)$.

Recall the definition of twist-moves in \S\ref{localn}. 
By the above two paragraphs,  
$K_+\otimes[2]$ and $K_-\otimes[2]$ differ by one twist-move. 

This completes the proof of the $\nu=0$ case.

\smallbreak

Next we prove the $\nu>1$ case. 

If the `$2m=2\alpha$ and $\nu=1$' case is true, then 
the `$2m=2\alpha+2$ and $\nu=1$' case is true. 
Hence `$2m=2\alpha$ and $\nu=2$' case is true. 
By the induction, Theorem \ref{Tokyo} is true.   \qed

\bigbreak

Here, we prove the following Proposition \ref{zero}, which is used 
in Proof of Theorem \ref{proinertia} in \S\ref{propol}.

Suppose that two $(2m+1)$-dimensional oriented closed submanifolds $\subset S^{2m+3}$, 
$K_+$ and $K_-$,  differ by only one twist-move 
as in Proof of Theorem \ref{Tokyo}. 
By Theorem \ref{Tokyo}, 
the $(2m+2\nu+1)$-submanifolds 
$\subset S^{2m+2\nu+3}$,  
$K_+\otimes^\nu[2]$ 
and 
$K_-\otimes^\nu[2]$,   
differ by one twist-move in a $(2m+2\nu+3)$-ball $\check A$.
Then there is a unique 
closed oriented $(2m+2\nu+1)$-submanifold $K_0^\otimes\subset S^{2m+2\nu+3}$ 
such that a triple 
$(K_+\otimes^\nu[2]$,  $K_-\otimes^\nu[2]$, $K_0^\otimes)$ 
is related by a twist-move in $\check A$. 
Note that the equivalence class of the submanifold 
$K_0^\otimes\subset S^{2m+2\nu+3}$  
is determined uniquely. 
Note that we have the following. 
Take 
the Seifert hypersurface for $\partial M_*$ 
in the $(2m+2\nu+1)$ case of Proof of Theorem \ref{Tokyo}  ($*=+,-$).
Then \newline
$\partial 
(\hskip1mm(\text{the Seifert hypersurface for } \partial M_*)-\rm{Int}\check A\hskip1mm)$ 
in $S^{2m+2\nu+3}$ is $K_0^\otimes\subset S^{2m+2\nu+3}$  ($*=+,-$).
Take $K_0$ in Proof of Theorem \ref{Tokyo}. 
Then we have the following.

\begin{prop}\label{zero}   
The $(2m+2\nu+1)$-submanifolds $\subset S^{2m+2\nu+3}$, \newline
$K_0^\otimes$ and 
$K_0\otimes^\nu($the empty knot $[2])$,   
are equivalent. 
\end{prop}

\smallbreak
\noindent{\bf Proof of Proposition \ref{zero}.} 
It suffices to prove the $\nu=1$ case.

Take $V'_+$ and $V'_-$ in Proof of Theorem \ref{Tokyo}.  
The submanifolds \newline
$V'_+\cap \overline{\hskip1mm B^{2m+4}-(A\x[0,1])\hskip1mm}$  in  
$\overline{\hskip2mm B^{2m+4}-(A\x[0,1])\hskip2mm }$
and \newline
$V'_-\cap \overline{\hskip1mm B^{2m+4}-(A\x[0,1])\hskip1mm }$  in  
$\overline{\hskip1mm B^{2m+4}-(A\x[0,1])\hskip1mm }$
are equivalent.

\smallbreak
Furthermore, the submanifold \newline
$\partial(\hskip1mm V'_+\cap \overline{\hskip1mm B^{2m+4}-(A\x[0,1])\hskip1mm}\hskip1mm )$  
in  
$\partial(\hskip1mm \overline{\hskip1mm B^{2m+4}-(A\x[0,1])\hskip1mm }\hskip1mm )$ 
($*=+,-$) 
is equivalent to $K_0$ in the standard $(2m+3)$-sphere.

\smallbreak
Take 
$M_*\cap
(\hskip1mm \overline{\hskip1mm B^{2m+4}-(A\x[0,1])\hskip1mm }\times B^2\hskip1mm )$ ($*=+, -$). 
By the construction, 
the submanifold 
$\partial(M_*\cap
(\hskip1mm \overline{\hskip1mm B^{2m+4}-(A\x[0,1])\hskip1mm }\times B^2\hskip1mm )$ 
in 
$\partial
(\hskip1mm \overline{\hskip1mm B^{2m+4}-(A\x[0,1])\hskip1mm }\times B^2\hskip1mm )$ 
($*=+, -$)  
is 
$K_0^\otimes$ and 
$K_0\otimes^\nu[2]$    in the standard $(2m+2\nu+3)$-sphere. 
\qed

\bigbreak

\noindent{\bf Proof of Theorem \ref{Sunday}. } 

\noindent{\bf Proof of Theorem \ref{Sunday}.(1).}
Take the 1-knot $K$ in the following figure.

\includegraphics[width=10cm]{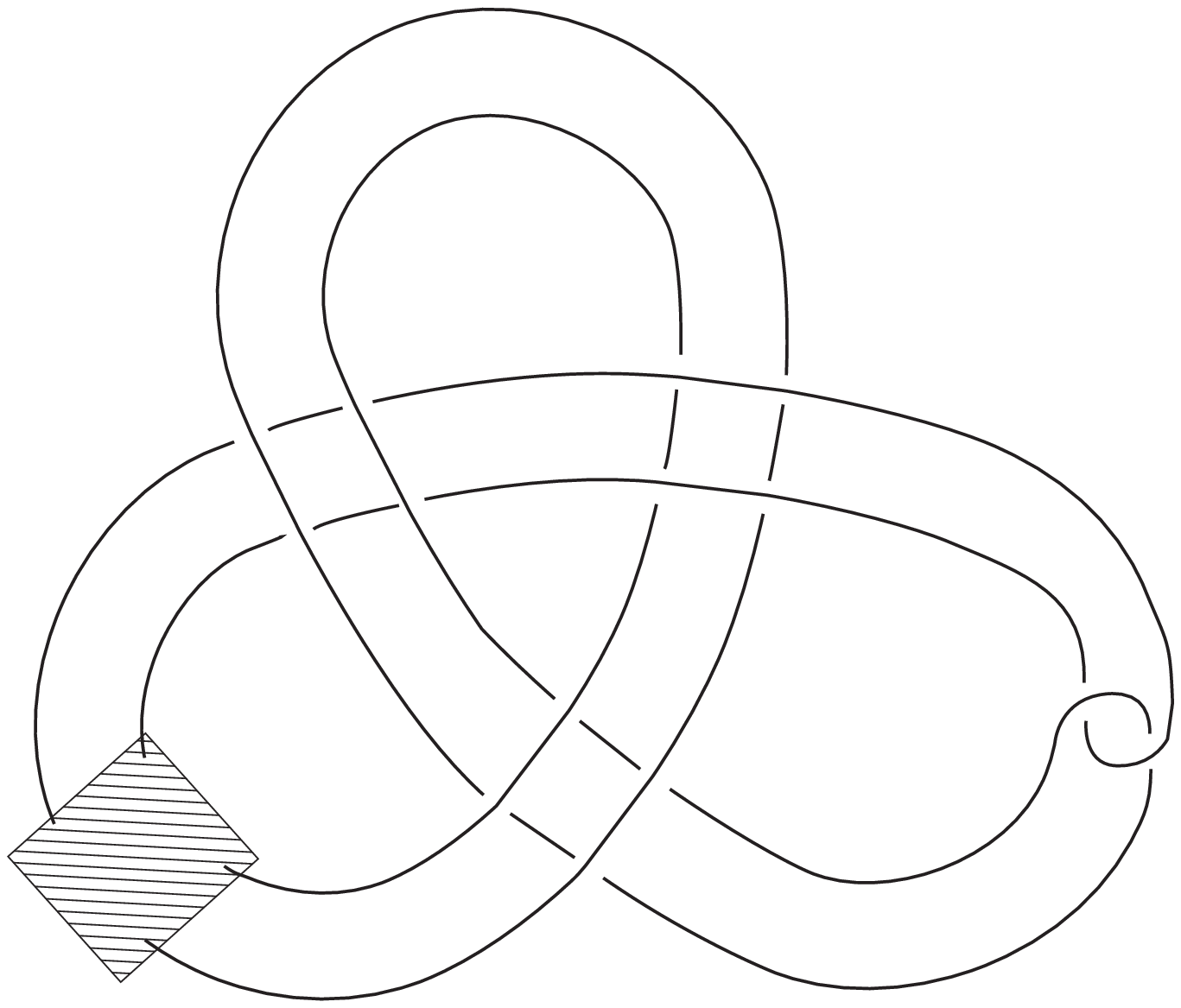}
\vskip-1cm

Twist in the shaded part so that its Seifert matrix is 
$\begin{pmatrix}
1&1\\
0&0\\
\end{pmatrix}$.

\bigbreak
By \cite{Whitten},  $K$ is a nontrivial knot. 
Note that the unknotting number of $K$ is one.

By \cite{Kauffman, KauffmanNeumann} and Proposition \ref{KLS},   
the $(2\nu+1)$-submanifold 
$K\otimes^\nu[2]$ 
$\subset S^{2\nu+3}$ has a Seifert hypersurface $V$ 
with the following conditions. 

\smallbreak
\noindent 
(1) 
$V$ has a handle decomposition 

\noindent 
(a $(2\nu+2)$-dimensional 0-handle)$\cup((2\nu+2)$-dimensional $(\nu+1)$-handles), where 
there may be no $(2\nu+2)$-dimensional $(\nu+1)$-handle. 

\noindent 
(2)  
A Seifert matrix $S$ associated with $V$ 
 is 
$\begin{pmatrix}
1&1\\
0&0\\
\end{pmatrix}$
or 
$\begin{pmatrix}
-1&-1\\
0&0\\
\end{pmatrix}$.

Note that 
$S+(-1)^{\nu+1}(^t\hskip-1mm  S)$ 
represents  
the intersection product of  
$H_{\nu+1}(V;\Z)/$(Tor).  
Recall $\nu\geqq2$.
By (1),  Tor$H_{\nu+1}(V;\Z)\cong0$. 
Since the determinant of 
$S+(-1)^{\nu+1}(^t\hskip-1mm  S)$ is  $+1$ or $-1$, 
$\partial V=K\otimes^\nu[2]$ 
is a homology sphere. 
By (1) right above, 
$\pi_1\partial V=1$.  
By \cite{Smale},  
$\partial V$ 
is homeomorphic 
to the standard sphere.  
Hence 
the $(2\nu+1)$-submanifold 
$K\otimes^\nu[2]$ 
is a spherical knot. 
By (1)  right above, 
the $(2\nu+1)$-knot  
$K\otimes^\nu[2]$ 
is a simple knot.

\noindent
Recall that the trivial $(2\nu+1)$-knot has Seifert matrices,  
$\begin{pmatrix}
1&1\\
0&0\\
\end{pmatrix}$
and 
$\begin{pmatrix}
-1&-1\\
0&0\\
\end{pmatrix}$
. 

\noindent   
By \cite{Levinesimp}, 
the $(2\nu+1)$-knot  
$K\otimes^\nu[2]$ 
is equivalent to the trivial knot. 
\qed

\smallbreak
\noindent{\bf Proof of Theorem \ref{Sunday}.(2).}
Take the above nontrivial 1-knot $K$. 
Take the knot-sum $K\sharp K$. 
By \cite{Homma, Papa} and the van Kampen's theorem,   
$K\sharp K$ is nontrivial.   
By \cite{Scharlemann},  the unknotting number of $K\sharp K$ is two. 

Note $(K\sharp K)\otimes^\nu[2]$ is the trivial knot. 

$K\sharp K$ is an example which we want. 
\qed

\smallbreak
\noindent{\bf Proof of Theorem \ref{Sunday}.(3).}
The pair of $K$ and $K\sharp K$ is an example which we want. 
Reason: 
$K$ is a prime knot and  $K\sharp K$ is not. 
Hence  $K$ is not equivalent to $K\sharp K$. 
$(K\sharp K)\otimes^\nu[2]$ 
and 
$K\otimes^\nu[2]$ are trivial knots and hence equivalent spherical knots.   
\qed

\smallbreak
\noindent{\bf Proof of Theorem \ref{Sunday}.(4).}
Let $A$ be the trefoil knot (and hence a nontrivial 1-knot). 
By \cite{Homma, Papa} and the van Kampen's theorem, 
  $A\sharp K$ is nontrivial. 
By \cite{Scharlemann},   the unknotting number of $A\sharp K$ is two. 
By \cite{Kauffman, KauffmanNeumann} and \cite{Levinesimp}, 
$(A\sharp K)\otimes^\mu$(the Hopf link) and  
$A\otimes^\mu$(the Hopf link) are equivalent and nontrivial . 
By Theorem \ref{Illinois}, 
$A\otimes^\mu$(the Hopf link) 
is obtained from the trivial knot by a twist-move. 
Hence 
$(A\sharp K)\otimes^\mu$(the Hopf link) is obtained from the trivial knot by a twist-move. 
Hence $A\sharp K$ is an example which we want.   \qed

\smallbreak\noindent{\bf Proof of Theorem \ref{Sunday}.(5).}
The pair of $A$ and $A\sharp K$ is an example which we want.
\qed


\bigbreak
\noindent {\bf Proof of $(*)$ in Note \ref{tori}.}   
Take the nontrivial 1-knot $K$ 
in Proof of Theorem \ref{Sunday}.   
The diffeomorphism type of 
$K\otimes[2]$  is the 3-manifold 
which is the double branched covering space of $S^3$ along $K$ 
(see \cite{Kauffman, KauffmanNeumann}.)
It is not homeomorphic to the standard 3-sphere 
(see \cite{MorganBass, Smith}.)  
$K$ is an example which we want. 
\qed

\bigbreak 
In order to prove Theorem \ref{aoiro}, we prove some propositions and a theorem.

\begin{prop}\label{KKN2da}   
Let $l\in\N$. 
Let $K$ be a simple $(2l+1)$-knot. 
Then \newline
$K\otimes($the Hopf link$)$ is a simple knot. 
\end{prop}

\noindent
{\bf Proof of Proposition \ref{KKN2da}.} 
By \cite{Kauffman, KauffmanNeumann} and Proposition \ref{KLS}, 
$K\otimes($the Hopf link$)$ 
satisfies the following. 

\smallbreak
\noindent
(1) 
A Seifert hypersurface  
has a handle decomposition of 
one $(2l+6)$-dimensional 0-handle and $(2l+6)$-dimensional $(l+3)$-handles, 
where there may not be $(2l+6)$-dimensional $(l+3)$-handles.

\noindent
(2) 
A $(l+3)$-Seifert matrix $Y$ associated with the above Seifert hypersurface 
is \newline
$(-1)\x$ (a $(l+1)$-Seifert matrix of $K$).

\smallbreak
By these (1) and (2) and \cite{Smale}, 
$K\otimes($the Hopf link)  is a spherical knot and a simple knot. \qed

\bigbreak

\begin{prop}\label{aiiro}  
Let $p\in\N$. 
Suppose that $K$ is a simple $(2p+5)$-knot and 
that, if $2p+5=7$, the signature of $K$ is a multiple of 16. 
Then there is a simple $(2p+1)$-knot $A$ with the following properties. 

\smallbreak
\noindent$\mathrm{(i)}$
$K$ is equivalent to $A\otimes($the Hopf link$)$. 

\noindent$\mathrm{(ii)}$
If $X$ is a $(p+3)$-Seifert matrix of $K$, 
then $-X$ is a $(p+1)$-Seifert matrix of $A$. 

\noindent$\mathrm{(iii)}$
The equivalence class of such a knot is unique. 
\end{prop}

\noindent{\bf Proof of Proposition \ref{aiiro}.}
Take a $(p+3)$-Seifert matrix $X$ of $K$. 
By \cite{Levinesimp} and Proposition \ref{KLS},  
there is a simple $(2p+1)$-knot $A$ 
such that  a $(p+1)$-Seifert matrix is the matrix $-X$.  
By Proposition \ref{KLS} and \ref{KKN2da},  
$A\otimes($the Hopf link$)$ is a simple $(2p+5)$-knot
such that  a $(p+2)$-Seifert matrix is the matrix $X$.  
By \cite{Levinesimp}, 
 $A\otimes($the Hopf link$)$ is equivalent to $K$.  
By \cite{Levinesimp}, 
(iii) holds. 
\qed

\bigbreak

\begin{prop}\label{KKN}  
Let $K$ be a 1-knot. 
Let $\mu\in\N$. 
Then $K\otimes^\mu($the Hopf link$)$
is a $(4\mu+1)$-dimensional simple knot $($and hence  a spherical knot$).$   
\end{prop}

\noindent{\bf Proof of Proposition \ref{KKN}.}   
By \cite{Kauffman, KauffmanNeumann} and Proposition \ref{KLS}, 
$K\otimes^\mu($the Hopf link$)$ satisfies the following. 

\smallbreak
\noindent 
(1)  
A Seifert hypersurface has a handle decomposition of 
one $(4\mu+2)$-dimensional 0-handle and 
$(4\mu+2)$-dimensional $(2\mu+1)$-handles, 
where there may not be $(2\mu+1)$-handles.

\noindent 
(2)  
Let $Y$ be a Seifert matrix of $K$. 
A Seifert matrix associated with the above Seifert hypersurface 
is $Y$ or $-Y$. 

\smallbreak
By (2) and 
the fact that $K$ is diffeomorphic to the single circle,  
det $(Y-^t\hskip-1mm Y)$ is $\pm1$. 
By this fact and the above  (1) (2) and \cite{Smale},  
$K\otimes^\mu($the Hopf link$)$ 
is a spherical knot and a simple knot. 
\qed

\bigbreak

\begin{prop}\label{shuiro}  
Let $K$ be a simple 5-knot. 
Then there is a 1-knot $A$ such that $K$ is equivalent to $A\otimes($the Hopf link$)$. 
\end{prop}

\noindent{\bf Proof of Proposition \ref{shuiro}.}
Take a 3-Seifert matrix $X$ of $K$. 
By using a Seifert surface,  there is a 1-knot $A$ 
such that  a Seifert matrix is $-X$.  
By Proposition \ref{KLS} and \ref{KKN},  
$A\otimes($the Hopf link$)$ is a simple 5-knot
such that  a 3-Seifert matrix is $X$.  
By \cite{Levinesimp}, 
 $A\otimes($the Hopf link$)$ is equivalent to $K$.  
\qed

\smallbreak
\noindent{\bf Note. } 
The equivalence class of $A$ is not unique. 
There are countably infinitely many equivalence types of 1-knots of this property.
Reason: Use the nontrivial knot $K$ in Proof of Theorem \ref{Sunday}. 
Take knot-sums.

\begin{prop}\label{nichiyo}  
Let $k\in\N$.   
Let $K$ be a simple $(4k+1)$-knot. 
Let $J$ be a $(4k+1)$-submanifold in $S^{4k+3}$ such that 
$J$ is obtained from $K$ by one twist-move. 
Then $J$ is a simple knot. 
\end{prop}

\noindent{\bf Proof of Proposition \ref{nichiyo}.}
By the definition of twist-moves, 
there is a $(4k+3)$-ball $B$ trivially embedded in $S^{4k+3}$ 
in which this twist-move is carried out.

By using the Thom-Pontrjagin construction, 
we can prove that 
there is a Seifert hypersurface $V_K$ for $K$ and $V_J$ for $J$
with the following properties. 

\smallbreak
\noindent(1)
$V_K\cap B$ (resp. $V_J\cap B$) is a $(4k+2)$-dimensional $(2k+1)$-handle $h$ 
that is attached to $\partial B$ as explained in 
the definition of the twist-moves in \S\ref{localn}.   

\noindent(2) 
$V_K\cap (S^{4k+3}-$Int $B)=V_J\cap (S^{4k+3}-$ Int $B)$.  
\smallbreak

The handle $h$ makes an order zero $(2k+1)$-cycle in  $V_K$ (resp. $V_J$).   
The intersection product between an order zero $(2k+1)$-cycle and itself in 
a compact oriented $(4k+2)$-manifold is zero.  
%
By using this fact, the Meyer-Vietoris exact sequence and the van Kampen's theorem, 
$J$ is 
homeomorphic to the standard sphere.

Let $N(J)$ (resp. $N(K)$) be the tubular neighborhood of $J$ (resp. $K$) in $S^{4k+3}$.  
Then we have

\noindent 
$S^{4k+3}-$Int$N(K)=$
$(S^{4k+3}-$Int$N(K)-$Int$B)\cup($a $(4k+3)$-dimensional $(2k+2)$-handle). 

\noindent 
$S^{4k+3}-$Int$N(J)=$
$(S^{4k+3}-$Int$N(J)-$Int$B)\cup($a $(4k+3)$-dimensional $(2k+2)$-handle).

By the definition of twist-moves, 
$S^{4k+3}-$Int$N(K)-$Int$B=S^{4k+3}-$Int$N(J)-$Int$B$.  

Since $K$ is a simple knot,  $\pi_1(S^{4k+3}-$Int$N(K))=\Z$.
Use the van Kampen's theorem for the above unions. 
Hence  $\pi_1(S^{4k+3}-$Int$N(J))=\Z$. 



Let $i\in\N$ and $i\leqq p$. 
there is an $i$-Seifert matrix $X_{iK}$ (resp. $X_{iJ}$) for $K$ (resp. $J$)
such that 
$X_{iK}=X_{iJ}$. (Reason: Use $V_K$ and $V_J$.) 
Consider the homology groups, the homotopy groups, and the fundamental group of 
the infinite cyclic covering space for $K$ (resp. $J$). 
Hence $J$ is a simple knot. 
\qed

\bigbreak
We prove the following theorem. 
\begin{thm}\label{apple}    
Let $n\in\N\cup\{0\}.$
Let $K$ be an $n$-dimensional closed oriented submanifold $\subset S^{n+2}.$
Take a map 

\smallbreak
\hskip3cm$K\mapsto K\otimes($the Hopf link$)$ 
\smallbreak

\noindent
from
the set of 
$n$-dimensional closed oriented submanifolds $\subset S^{n+2}$
to
the set of 
$(n+4)$-dimensional closed oriented submanifolds $\subset S^{n+6}$.

Then we have the following.

\smallbreak
\noindent$(1)$ 
Let $K$ be a simple $(2l+1)$-knot $(l\geqq2, \quad l\in\N).$   
That is, suppose that the domain of the map is the set of simple $(2l+1)$-knots.  
Then the image of the map is 
the set of simple $(2l+5)$-knots.
Furthermore  
the map 

\smallbreak
\hskip3cm$\{$simple $(2l+1)$-knots$\}\to\{$simple $(2l+5)$-knots$\}$.
\smallbreak
\hskip68mm$K\mapsto K\otimes($the Hopf link$)$ 
\smallbreak

\noindent 
is a one-to-one map.

\noindent$(2)$ 
Let $K$ be a simple 3-knot. 
That is, suppose that the domain of the map is the set of simple 3-knots.  
Then the image of the map is included in 
the set of simple 7-knots.
Furthermore the map 
\smallbreak
\hskip3cm$\{$simple 3-knots$\}\to\{$simple 7-knots$\}$.
\smallbreak
\hskip56mm$K\mapsto K\otimes($the Hopf link$)$ 
\smallbreak

\noindent is injective but not onto.  

\smallbreak

\noindent$(3)$ 
Let $K$ be a 1-knot. 
That is, suppose that the domain of the map is the set of 1-knots.  
Then the image of the map is  
the set of simple 5-knots.
Furthermore the map 
\smallbreak
\hskip3cm$\{1$-knots$\}\to\{$simple 5-knots$\}$.
\smallbreak
\hskip43mm$K\mapsto K\otimes($the Hopf link$)$ 
\smallbreak
\noindent is onto but not injective.  
The inverse image of any element by this map is an infinite set.

\end{thm}

\noindent
\begin{prob}\label{tora}   
What happens if we define the domain is another set in Theorem \ref{apple}? 
\end{prob}

\noindent{\bf Proof of Theorem \ref{apple}.}    
%
Proposition \ref{KKN2da}, \ref{aiiro}  imply  Theorem \ref{apple}.(1).

There is a simple 7-knot with the following property $(*)$:  

\noindent
$(*)$ The signature is a multiple of 8 but not a multiple of 16.

\noindent
There is not a simple 3-knot with the above property $(*)$.
  See \cite{Levinesimp}. 

By these facts and 
Proposition \ref{KKN2da}, \ref{aiiro},  
Theorem \ref{apple}.(2) holds. 

Theorem \ref{apple}.(3) follows from 
Proposition \ref{KKN}, \ref{shuiro}, 
Note to Proposition \ref{shuiro}, 
and
Theorem \ref{Sunday}. \qed

\bigbreak 
\noindent
{\bf Proof of  Theorem \ref{aoiro}. }  
By the definition of the twist-move in \S\ref{localn} 
 there is 
a $(2k+3)$-Seifert matrix $X$  (resp. $Y$) for $J$ (resp. $K$)   
with the following properties. 

\smallbreak
\noindent 
(1) $X$ and $Y$ are $c\times c$ matrices  for a natural number $c$.  

\noindent
(2) Let $x_{ij}$ denote $(i,j)$-element of $X$.  
Let $y_{ij}$ denote $(i,j)$-element of $Y$.  
There is a natural number $a$  $\leqq c$ such that 
$$
\left\{
\begin{array}{ll}
x_{ij}=y_{ij}-1 & \text{ if $(i,j)=(a, a)$}\\
x_{ij}=y_{ij}    & \text{ if $(i,j)\neq(a, a).$}\\
\end{array}
\right.
$$

\noindent(3)   
$X$ and $Y$ are not S-equivalent.  See \cite{Levinesimp} for S-equivalent. 
\smallbreak

Note that we can take Seifert matrices which satisfy the above conditions.
If necessarily, carry out surgery on a Seifert hypersurface 
by handles embedded in $S^{4k+7}$. 
\bigbreak

By Theorem \ref{apple}, $K$ is a simple knot. 
By Proposition \ref{nichiyo}, $J$ is a simple knot.

By Proposition \ref{KLS} and Theorem \ref{apple}, 
there are simple $(4k+1)$-knots $K'$ and $J'$ 
with the following properties. 

\smallbreak
\noindent $\mathrm{(i)}$
$J$ is equivalent to $J'\otimes($the Hopf link$)$.     
$-X$ is a $(2k+1)$-Seifert matrix for $J'$. 

\noindent $\mathrm{(ii)}$
$K$ is equivalent to $K'\otimes($the Hopf link$)$.   
$-Y$ is a $(2k+1)$-Seifert matrix for $K'$.

\noindent $\mathrm{(iii)}$
The equivalence class of such a knot is unique.

\bigbreak
Therefore we can make 
a $(4k+1)$-dimensional simple knot $\bar J$ (resp. $\bar K$) with the following properties. 

\bigbreak
\noindent 
(I) A handle decomposition of a Seifert hypersurface is a set of  
a $(4k+2)$-dimensional 0-handle
and 
$(4k+2)$-dimensional $(2k+1)$-handles, 
where there may not be a $(4k+2)$-dimensional $(2k+1)$-handle. 

\noindent
(II) A Seifert matrix associated with this Seifert hypersurface is $-X$ (resp. $-Y$). 

\noindent
(III) $\bar J$ and $\bar K$ differ by a single twist-move  
and are nonequivalent.

\noindent 
Reason: 
Since $J$ (resp. $K$) is  
 homeomorphic  to the standard sphere,   
we can realize (I) (II)(III) by using $(2k+1, 2k+1)$-pass-moves. 

\vskip3mm
By \cite{Levinesimp}, 
simple $(4k+1)$-dimensional spherical knot $\bar J$ (resp. $\bar K$) 
is equivalent to $J'$ (resp. $K'$).

This completes the proof. \qed
\bigbreak

\noindent{\bf {Note.}} 
Let $p$ be any natural number. 
There are countably infinitely many pair $(J, K)$ of $(2p+5)$-dimensional knots $J$ and $K$ with the following properties: 
Neither $K$ or $J$ is the product of any $(2p+1)$-knot and the Hopf link. 
$J$ is obtained from $K$ by a twist-move.  

Reason: It is well-known that there are countably infinitely many nonsimple  
$(2p+5)$-knots for any $p$ such that the fundamental group of the complement of each knot is not $\Z$.  
Let $K$ be such a knot.
It is trivial that there is a nontrivial knot $A$ 
which is obtained from the trivial knot by a twist-move.   
 Let $J=K\sharp A$.   
Then $J$ is a nonsimple knot. 
By Theorem \ref{apple},  
neither $K$ or $J$ is the product of any $(2p+1)$-knot and the Hopf link.

\section{Proof of Theorems in \S \ref{thepa}}\label{propa} 
\noindent{\bf Proof of Theorem \ref{Waltham}.}    
There is a Seifert matrix $X$ (resp. $Y$) for the 1-knot $J$ (resp. $K$). 

\smallbreak
\noindent 
(1) $X$ and $Y$ are $c\times c$ matrices  for a natural number $c$.  

\noindent
(2) Let $x_{ij}$ denote $(i,j)$-element of $X$.  
Let $y_{ij}$ denote $(i,j)$-element of $Y$.  
There are natural numbers $a, b$  $\leqq c$ such that  $a\neq b$ 
and that 
$$
\left\{
\begin{array}{ll}
x_{ij}=y_{ij}-1 & \text{ if $(i,j)=(a, b)$}\\
x_{ij}=y_{ij}    & \text{ if $(i,j)\neq(a, b)$ and if $(i,j)\neq(b, a). $}\\
\end{array}
\right.
$$
\noindent 
The $(b, a)$-element is determined by the $(a, b)$-element. 

Note that we can take Seifert matrices which satisfy the above conditions.
If necessarily, carry out surgeries on Seifert hypersurfaces 
by 2-dimensional 1-handles embedded in $S^3$. 

\smallbreak





By Proposition \ref{KKN}, 
$J\otimes^\mu($the Hopf link$)$ 
(resp. $K\otimes^\mu($the Hopf link$)$)
is a spherical knot and is a simple knot.

\smallbreak

We can make 
a $(4\mu+1)$-dimensional spherical  knot $J'$ (resp. $K'$) with the following properties. 

\smallbreak
\noindent 
(1)A handle decomposition of a Seifert hyper surface is a set of 
a $(4\mu+2)$-dimensional $(2\mu+1)$-handle  
and 
$(4\mu+2)$-dimensional $(2\mu+1)$-handles, 
where there may not be $(4\mu+2)$-dimensional $(2\mu+1)$-handles. 

\noindent
(2) A Seifert matrix associated with the above Seifert hypersurface is 
$(-1)^\mu X$ (resp. $(-1)^\mu Y$). 

\noindent
(3)$J'$ is obtained from $K'$ by a $(2\mu+1, 2\mu+1)$-pass-move. 

\noindent
Reason: Since $J$ (resp. $K$) is  diffeomorphic to the single circle, 
we can realize (1)(2)(3)  by using $(2\mu+1, 2\mu+1)$-pass-moves.

\vskip3mm
By \cite{Levinesimp}, 
simple $(4\mu+1)$-dimensional spherical knot $J'$ (resp. $K'$)
 is equivalent to 
$J\otimes^\mu($the Hopf link$)$ 
(resp. $K\otimes^\mu($the Hopf link$)$).

By the construction of 
$J'$ (resp. $K'$), 
$J\otimes^\mu($the Hopf link$)$ 
and 
 $K\otimes^\mu($the Hopf link$)$
differ by one $(2\mu+1, 2\mu+1)$-pass-move. 
\qed

\vskip3mm
\noindent
{\bf Proof of Theorem \ref{sky}.} 
The pass-move does not change diffeomorphism type of submanifolds. 
However  
$J\otimes^\mu[2]$ 
and
$K\otimes^\mu[2]$ 
do not have the same homeomorphism type in general by 
\cite{Kauffman, KauffmanNeumann}.  
Example: (The trivial 1-knot)$\otimes[2]$ and  (the trefoil knot)$\otimes[2]$. 
\qed

\vskip3mm
\noindent
{\bf Proof of Theorem \ref{mountain}.}   
Note $l\geqq1$. 
There is a Seifert matrix $X$  (resp. $Y$) for 
the simple $(2l+1)$-knot $J$ (resp. $K$) with the following properties.

\smallbreak
\noindent 
(1) $X$ and $Y$ are $c\times c$ matrices  for a natural number $c$.  

\noindent
(2) Let $x_{ij}$ denote $(i,j)$-element of $X$.  
Let $y_{ij}$ denote $(i,j)$-element of $Y$.  
There are integers $a, b$  $\leqq c$ such that $a\neq b$. 
We have 
$$
\left\{
\begin{array}{ll}
x_{ij}=y_{ij}-1 & \text{ if $(i,j)=(a, b)$}\\
x_{ij}=y_{ij}    & \text{ if $(i,j)\neq(a, b)$ and if $(i,j)\neq(b, a) $}\\
\end{array}
\right.
$$
\noindent 
The $(b, a)$-element is determined by the $(a, b)$-element. 
\smallbreak

Note that we can take Seifert matrices which satisfy the above conditions.
If necessarily, carry out surgeries on Seifert hypersurfaces 
by $(2l+2)$-dimensional $(l+1)$-handles embedded in $S^{2l+3}$.

By Proposition \ref{KKN2da}, 
$J\otimes^\mu($the Hopf link$)$ 
(resp. $K\otimes^\mu($the Hopf link$)$)
is a spherical knot and a simple knot. 
\smallbreak



We can make 
a $(2l+4\mu+1)$-dimensional spherical knot $J'$ (resp. $K'$) with the following properties. 

\smallbreak
\noindent 
(1) A handle decomposition of a Seifert hypersurface is a set of  
a $(2l+4\mu+2)$-dimensional 0-handle
and 
$(2l+4\mu+2)$-dimensional $(l+2\mu+1)$-handles. 

\noindent
(2) A Seifert matrix associated with the above Seifert hypersurface is 
$(-1)^\mu X$ (resp. $(-1)^\mu Y$).

\noindent
(3) $J'$ is obtained from $K'$ 
by a $(l+2\mu+1, l+2\mu+1)$-pass-move. 

\noindent 
Reason:  Since $J$ (resp. $K$) is  
homeomorphic to the standard sphere,   
we can realize (1)(2)(3) by using $(l+2\mu+1, l+2\mu+1)$-pass-moves.

\vskip3mm
By \cite{Levinesimp}, 
simple $(2l+4\mu+1)$-dimensional spherical knot $J'$ (resp. $K'$)
 is equivalent to 
$J\otimes^\mu($the Hopf link$)$ 
(resp. $K\otimes^\mu($the Hopf link$)$).

By the construction of 
$J'$ (resp. $K'$), 
$J\otimes^\mu($the Hopf link$)$ and  
$K\otimes^\mu($the Hopf link$)$ differ by one $(l+2\mu+1, l+2\mu+1)$-pass-move.  
\qed

\smallbreak
\noindent{\bf Proof of Theorem \ref{grape}.}   
Take the nontrivial 1-knot $K$ in the figure in Proof of Theorem \ref{Sunday}. 
Note that $K$ is obtained from the trivial knot by one pass-move. \qed

\smallbreak

\noindent {\bf Proof of $(*)$ in Note \ref{bird}.}  
It is same as  Proof of $(*)$ in Note \ref{tori}. \qed

\smallbreak

\noindent{\bf Proof of Theorem \ref{ki}.}  
We need a proposition.

\begin{prop}\label{getsuyo}   
Let $p\in\N$.   
Let $K$ be a simple $(2p+1)$-knot. 
Let $J$ be a $(2p+1)$-submanifold in $S^{2p+3}$ such that 
$J$ is obtained from $K$ by one $(p+1,p+1)$-pass-move. 
Then $J$ is a simple knot. 
\end{prop}

\noindent{\bf Proof of Proposition \ref{getsuyo}.}
By the definition of the $(p+1,p+1)$-pass-move, 
$J$ is a spherical $(2p+1)$-knot.   

See the following figure. 
We can take two copies of the $(p+1)$-sphere $Y_1$ and $Y_2$ in $B^{2p+3}\subset S^{2p+3}$ with the following property.    
$Y_1$ (resp. $Y_2$) is embedded trivially. 
The linking number of $Y_1$ and $Y_2$ is one. 
Carry out surgeries 
along two  $(p+1)$-spheres  $Y_1$ and $Y_2$   
by two $(2p+3)$-dimensional $(p+2)$-handles with the trivial framing 
on $B^{2p+3}$.  
Then $B^{2p+3}$ becomes the $(2p+3)$-ball again and 
$S^{2p+3}$ becomes the $(2p+3)$-sphere again. 
Furthermore the $(p+1,p+1)$-pass-move in the $(2p+3)$-ball $B^{2p+3}$ is done. 


Since  $p\geqq1$ holds and $K$ is a simple knot, 
$
\pi_i(S^{2p+3}-N(K))=
\pi_i(S^{2p+3}-N(K))
$ 
for $1\leqq i\leqq p)$. (Use the van Kampen's theorem and the Meyer-Vietoris theorem on the complements and the infinite cyclic covering spaces.) 
Therefore $J$ is a simple knot. 
This completes the proof of  Proposition \ref{getsuyo}.\qed

\bigbreak
\input surgery.tex
\bigbreak

By Theorem \ref{apple}, $K$ is a simple knot. 
By Proposition \ref{getsuyo},  
$J$ is a simple knot. 
By Proposition \ref{shuiro}   and Theorem \ref{apple}, 
there is a 1-knot  $J'$ to satisfy (i).

By \cite{Ogasa98n},   
$K'\otimes^\mu$(the Hopf link) 
and 
$J'\otimes^\mu$(the Hopf link) 
have the same Arf invariant.

By \cite{Kauffman, KauffmanNeumann}, 
a $(2\mu+1)$-Seifert matrix of 
$
\begin{cases}
\text{$K'\otimes^\mu$(the Hopf link) }\\
\text{$J'\otimes^\mu$(the Hopf link) }
\end{cases}
$\newline
is $(\pm1)\x$
a Seifert matrix of 
$
\begin{cases}
\text{$K'$}\\
\text{$J'$}
\end{cases}
$. 
Hence Arf $(K')=$Arf$(J')$.

By Theorem \ref{passth},  $K'$ is pass-move-equivalent to $J'$. 
Hence $J'$ satisfies (i) and (ii). 
This completes the proof of Theorem \ref{ki}. \qed

\smallbreak\noindent{\bf Note.}  
Proof of Note to Theorem \ref{ki}. 
Use the nontrivial knot $K$ in Proof of Theorem \ref{Sunday}. 
Take a knot-sum as many times as we need.

\smallbreak
\noindent{\bf Proof of Theorem \ref{aka}.}  
By Theorem \ref{apple}, $K$ is a simple knot. 
By Proposition \ref{getsuyo},   
$J$ is a simple knot. 

Then there is 
a $(p+3)$-Seifert matrix $X$  (resp. $Y$) for $J$ (resp. $K$)   
with the following properties. 

\smallbreak
\noindent 
(1) $X$ and $Y$ are $c\times c$ matrices  for a natural number $c$.

\noindent 
(2) Let $x_{ij}$ denote $(i,j)$-element of $X$.  
Let $y_{ij}$ denote $(i,j)$-element of $Y$.  
There are integers $a, b$  $\leqq c$ such that $a\neq b$. 
We have 
$$
\left\{
\begin{array}{ll}
x_{ij}=y_{ij}-1 & \text{ if $(i,j)=(a, b)$}\\
x_{ij}=y_{ij}    & \text{ if $(i,j)\neq(a, b)$ and if $(i,j)\neq(b, a). $}\\
\end{array}
\right.
$$
The $(b, a)$-element is determined by the $(a, b)$-element. 

\noindent(3)   
$X$ and $Y$ are not S-equivalent.  See \cite{Levinesimp} for S-equivalent. 
\smallbreak
Note that we can take Seifert matrices which satisfy the above conditions.
If necessarily, carry out surgeries on Seifert hypersurfaces  
by handles embedded in $S^{2p+7}$. 

\smallbreak

By Proposition \ref{aiiro} and Theorem \ref{apple}, 
there are simple $(2p+1)$-knots $K'$ and $J'$ 
with the following properties. 

\smallbreak
\noindent $\mathrm{(i)}$
$J$ is equivalent to $J'\otimes($the Hopf link$)$.     
$-X$ is a $(p+1)$-Seifert matrix for $J'$. 

\noindent $\mathrm{(ii)}$
$K$ is equivalent to $K'\otimes($the Hopf link$)$.   
$-Y$ is a $(p+1)$-Seifert matrix for $K'$.

\noindent $\mathrm{(iii)}$
The equivalence class of such a knot is unique.

\smallbreak
Therefore we can make 
a $(2p+1)$-dimensional simple knot $\bar J$ (resp. $\bar K$) with the following properties. 

\smallbreak
\noindent 
(I) A handle decomposition of a Seifert hypersurface is a set of  
a $(2p+2)$-dimensional 0-handle
and 
$(2p+2)$-dimensional $(p+1)$-handles, 
where there may not be a $(2p+2)$-dimensional $(p+1)$-handle. 

\noindent
(II) A Seifert matrix associated with the above Seifert hypersurface is $-X$ (resp. $-Y$). 

\noindent
(III) $\bar J$ and $\bar K$ differ by a single $(p+1,p+1)$-pass-move 
and are nonequivalent.

\noindent  
Reason: 
Since $J$ (resp. $K$) is  
homeomorphic to the standard sphere,   
we can realize (I) (II)(III) by using $(p+1, p+1)$-pass-moves. 

\smallbreak
By \cite{Levinesimp}, 
simple $(2p+1)$-dimensional spherical knot $\bar J$ (resp. $\bar K$) 
is equivalent to $J'$ (resp. $K'$).

This completes the proof. \qed
\smallbreak 

\noindent{\bf {Note.}} 
There are countably infinitely many pair $(J, K)$ of $(2p+5)$-dimensional knots $J$ and $K$ ($p\in\N$)
with the following properties: 
Neither $K$ or $J$ is the product of any $(2p+1)$-knot and the Hopf link. 
$J$ is obtained from $K$ by a $(p+1,p+1)$-pass-move.

Reason: It is well-known that there are countably infinitely many nonsimple knots 
such that the fundamental group of the complement of each knot is not $\Z$. 
Let $K$ be such a knot. 
It is trivial that there is a nontrivial knot $A$ 
which is obtained from the trivial knot  by a $(p+1,p+1)$-pass-move.  
Let $K$ be this $K$. Let $J=K\sharp A$.   
Then $J$ is a nonsimple knot. 
By Theorem \ref{apple},  
neither $K$ or $J$ is the product of any $(2p+1)$-knot and the Hopf link.

\smallbreak
\noindent{\bf Proof of Theorem \ref{spun}.}   
\cite{Kauffman, KauffmanNeumann} proved the following. 
Let $V$ be a Seifert hypersurface for an 
$n$-dimensional closed oriented submanifold 
$K\subset S^{n+2}$.
Then there is a Seifert hypersurface $W$ for the $(n+2)$-submanifold 
$K\otimes[2]$  $\subset S^{n+4}$ 
such that 
$W$ is diffeomorphic to $B^{n+3}\cup B^{n+3}$ and that 
$B^{n+3}\cap B^{n+3}$ is diffeomorphic to $V\times[-1,1]$.

Note that $P$ has a Seifert hypersurface which is diffeomorphic to 
the punctured Poincar\'e sphere.

Therefore, by the above theorem in \cite{Kauffman, KauffmanNeumann},  
the $(2\nu+2)$-submanifold 
$P\otimes^\nu[2]$ 
has a Seifert hypersurface $Q$ which consists of 
a $(2\nu+3)$-dimensional 0-handle,   
$(2\nu+3)$-dimensional $(\nu+1)$-handles,  and 
$(2\nu+3)$-dimensional $(\nu+2)$-handles. 
Note $\nu\geqq2$. 
By the van Kampen's theorem, $\pi_1(\partial Q)=1$. 
By using the Meyer-Vietoris exact sequence, 
$Q$ is a homology $(2\nu+3)$-ball. 
By  the van Kampen's theorem, $Q$ is simply-connected. 

By \cite{Smale}, $Q$ is diffeomorphic to  the $(2\nu+3)$-ball.

Hence 
$P\otimes^\nu[2]$  
has a Seifert hypersurface which is diffeomorphic to  the $(2\nu+3)$-ball. 

Hence 
$P\otimes^\nu[2]$ 
is equivalent to the trivial knot. 
\qed

\section{Proof of Theorems in \S \ref{thepol}}\label{propol}   

\noindent 
{\bf  Proof of Theorem \ref{skein1}.}  
By Theorem \ref{Hopf},  Theorem \ref{skein1} follows from Theorem \ref{skein2}. 
\qed

\smallbreak\noindent 
{\bf Proof of Theorem \ref{skein2}.}  
By \cite{KauffmanNeumann}, 
$\>K_*  \otimes^\nu
[2]\>$ 
has a Seifert hypersurface   $V_*$ with a handle decomposition of 
one (2$\nu$+2)-dimensional 0-handle and 
 (2$\nu$+2)-dimensional ($\nu$+1)-handles, 
where  there may not be (2$\nu$+2)-dimensional ($\nu$+1)-handles ($*=+,-,0$). 
%
Therefore 
the $\nu$-Alexander matrix 
associated with $V_*$ 
is `empty'. See \S\ref{Alex}.

By Proposition \ref{KLS}, 
a $(\nu+1)$-positive Seifert matrix
$$S_{\nu+1}(K_*\otimes^\nu[2])$$ 
and 
a $(\nu+1)$-negative Seifert matrix
$$N_{\nu+1}(K_*  \otimes^\nu[2])$$ 
are square matrices. 
By Proposition \ref{square},   
the $(\nu+1)$-$\Q[t, t^{-1}]$-Alexander polynomial 
of 
$K_*  \otimes^\nu[2]$      
is the $\Q[t, t^{-1}]$-balanced class of 
the determinant of 
the $(\nu+1)$-Alexander matrix 

\vskip3mm
$P_{\nu+1}(K_*   \otimes^\nu[2])=(p^*_{ij})
=t\cdot S_{\nu+1}(K_*\otimes^\nu[2])-N_{\nu+1}(K_*  \otimes^\nu[2])
$
\vskip3mm

%

\noindent
$=
t\cdot S_{\nu+1}(K_*\otimes^\nu[2])
+{(-1)^{\nu+1}}\cdot  ^t\hskip-1mm S_{\nu+1}(K_*\otimes^\nu[2]) 
$

\hskip5cm(Reason:    
$S_{\nu+1}(K_*  \otimes^\nu[2])=$
${(-1)^\nu}\hskip1mm ^t\hskip-1mm N_{\nu+1}(K_*\otimes^\nu[2])$ by \S\ref{Alex}.)

\vskip3mm
\noindent
$=
(-1)^{\frac{\nu(\nu-1)}{2}}\quad
(t\cdot S_1(K_*)
+{(-1)^{\nu+1}}\cdot  ^t\hskip-1mm S_1(K_*)).
$

\hskip3cm(Reason:    
$S_{\nu+1}(K_*\otimes^\nu[2])=(-1)^{\frac{\nu(\nu-1)}{2}}\quad S_1(K_*)$ 
by Proposition \ref{KLS}.(7).)

\vskip5mm
Since $(K_+,K_-,K_0)$ is a crossing-change-triple of 1-links,  
we can suppose that \newline
$S_1(K_*)=(s^*_{ij})$ 
 $(*=+,-,0)$
has the following properties.

\smallbreak
\noindent
(1) 
$(s^+_{ij})$ 
and 
$(s^-_{ij})$ 
are $\rho\x \rho$ matrices. 


\hskip3mm $(s^0_{ij})$ is a $(\rho-1)\x (\rho-1)$ matrix 
($\rho\geqq2, \quad \rho\in\N$). 

\smallbreak
\noindent
(2) 
\hskip4cm$s^+_{\rho,\rho}-s^-_{\rho,\rho}=1.$

\smallbreak
\noindent
(3)  
\hskip3cm$s^+_{ij}=s^-_{ij}=s^0_{ij}    (1\leqq i\leqq \rho-1, 1\leqq j\leqq \rho-1).$
$$s^+_{ij}=s^-_{ij} ((i,j)\neq(\rho,\rho)).$$ 
\smallbreak

\noindent
Note we can suppose that $\rho-1\geqq1$ 
by using surgeries of Seifert hypersurfaces, if necessary.

\vskip5mm 
 Therefore we can suppose that 
$P_{\nu+1}(K_*\otimes^\nu[2])=(p^*_{ij})$  
has the following properties ($*=+,-,0$). 

\smallbreak
\noindent
(1) 
$(p^+_{ij})$ 
and 
$(p^-_{ij})$  
are $\rho\x \rho$ matrices. 


\hskip3mm $(p^0_{ij})$ 
is a $(\rho-1)\x (\rho-1)$ matrix. 
($\rho\geqq2, \quad \rho\in\N$).

\smallbreak
\noindent
(2) 
\hskip4cm  $p^+_{\rho,\rho}-p^-_{\rho,\rho}=c(t+(-1)^{\nu+1}),$

\hskip5mm
where $c=\pm1$. 

\smallbreak
\noindent(3)     
\hskip3cm$p^+_{ij}=p^-_{ij}=p^0_{ij} 
(1\leqq i\leqq \rho-1,  1\leqq j\leqq \rho-1).$
$$p^+_{ij}=p^-_{ij}  ((i,j)\neq(\rho,\rho)).$$
\smallbreak

\vskip5mm
By calculus of determinants,   
$${\mathrm{det}}
P_{\nu+1}
(K_+  \otimes^\nu[2])
-
{\mathrm{det}}
P_{\nu+1}
(K_-
\otimes^\nu[2])
= c
(t+(-1)^{\nu+1})\cdot{\mathrm{det}}P_{\nu+1}(K_0\otimes^\nu[2]), 
$$

\noindent 
where $c=\pm1$.  

Hence Theorem \ref{skein2} holds. 
Note that the $\Q[t, t^{-1}]$-Alexander polynomial is a $\Q[t, t^{-1}]$-balanced class. 
\qed

\bigbreak\noindent 
{\bf Proof of Theorem \ref{pepper}.}  
The `$l=$even' case is proved in \cite{Ogasa09}. 

We prove the `$l=$odd' case.   

There is a Seifert hypersurface $V_*$ for $K_*$ ($*=+,-,0$) such that 
$(V_+, V_-, V_0)$ is related by a twist-move in $B^{2l+3}$.  





Take a positive Seifert matrix $S_\sharp(K_*)$  
and  a negative Seifert matrix $N_\sharp(K_*)$
for $K_*$ associated with $V_*$ ($*=+,-,0$). 
We suppose that 
$S_\sharp(K_*)$ and $N_\sharp(K_*)$ are defined by using the same ordered sets of the cycles. 
See \S\ref{Alex}. 

Therefore, we can suppose that 
$S_l(K_+)=S_l(K_-)=S_l(K_0)$ and that  
$N_l(K_+)=N_l(K_-)=N_l(K_0)$.

Since $K_+$ is a spherical knot, 
the linear map which is defined by $S_l(K_+)-N_l(K_+)$ is injective. 
Hence 
the linear map which is defined by $S_l(K_*)-N_l(K_*)$ is injective ($*=+,-,0$).

Hence 
the $l$-Alexander matrix  $t\cdot S_l(K_*)- N_l(K_*)$ 
associated with $V_*$
is injective. 
By Proposition \ref{square}, 
the $(l+1)$-$\Q[t, t^{-1}]$-Alexander polynomial 
of 
$K_*$      
is the $\Q[t, t^{-1}]$-balanced class of 
the determinant of 
the $(l+1)$-Alexander matrix 
%
%
%
\noindent
$P_{l+1}(K_*)=t\cdot S_{l+1}(K_*)-N_{l+1}(K_*)$.

Since $(K_+,K_-,K_0)$ is a twist-move-triple,  we can suppose that \newline
$S_{l+1}(K_*)=(s^*_{ij}) (*=+,-,0)$
has the following property.

\smallbreak
\noindent
(1) 
$(s^+_{ij})$ 
and 
$(s^-_{ij})$ 
are $\rho\x \rho$ matrices.

\hskip3mm $(s^0_{ij})$ is a $(\rho-1)\x (\rho-1)$ matrix. 
($\rho\geqq2, \quad \rho\in\N$).

\smallbreak\noindent
(2) 
\hskip4cm$s^+_{\rho,\rho}-s^-_{\rho,\rho}=1.$

\smallbreak
\noindent(3)  
\hskip3cm
$s^+_{ij}=s^-_{ij}=s^0_{ij}  (1\leqq i\leqq \rho-1, 1\leqq j\leqq \rho-1).$ 
$$s^+_{ij}=s^-_{ij}   ((i,j)\neq(\rho,\rho)).$$

\smallbreak
\noindent
Note we can suppose that $\rho-1\geqq1$ 
by using surgeries of Seifert hypersurfaces, if necessary.

\vskip5mm

We have 
$$P_{l+1}(K_*)=t\cdot S_{l+1}(K_*)-N_{l+1}(K_*)
=t\cdot S_{l+1}(K_*)-(-1)^l \cdot^t\hskip-1mm S_{l+1}(K_*).$$
Since $l$ is odd, 
$$P_{l+1}(K_*)=t\cdot S_{l+1}(K_*)+^t\hskip-1mm S_{l+1}(K_*).$$ 
Therefore   
$P_{l+1}(K_*)=(p^*_{ij})$ satisfies the following ($*=+,-,0$).

\smallbreak
\noindent
(1) 
$(p^+_{ij})$ 
and 
$(p^-_{ij})$  
are $\rho\x \rho$ matrices.

\hskip3mm $(p^0_{ij})$ is a $(\rho-1)\x (\rho-1)$ matrix 
($\rho\geqq2, \quad \rho\in\N$).

\smallbreak
\noindent
(2) 
\hskip4cm$p^+_{\rho,\rho}-p^-_{\rho,\rho}=t+1.$

\smallbreak\noindent(3)  
\hskip3cm$p^+_{ij}=p^-_{ij}=p^0_{ij} (1\leqq i\leqq \rho-1, 1\leqq j\leqq \rho-1).$
$$p^+_{ij}=p^-_{ij}  ((i,j)\neq(\rho,\rho)).$$ 

\vskip3mm

By calculus of determinants,   
$$
{\mathrm{det}}P_{l+1}(K_+)
-
{\mathrm{det}}P_{l+1}(K_-)
=(t+1)\cdot{\mathrm{det}}P_{l+1}(K_0).$$

Hence Theorem \ref{pepper}   
holds. 
\qed

\bigbreak
\noindent 
{\bf  Proof of Theorem \ref{skeincro}.}  

Take $V_*$,  $S_\sharp(K_*)$,   $N_\sharp(K_*)$ in Proof of Theorem \ref{pepper}.

Since $K_+$ is a spherical knot, 
we have that $S_l(K_+)$ and $N_l(K_+)$ are square matrix and furthermore 
we can suppose that 
the determinant of $S_l(K_+)-N_l(K_+)$ is $\pm1$. 
Hence 
we can suppose that 
the determinant of $S_l(K_*)-N_l(K_*)$ is $\pm1$ ($*=+,-,0$).

By \cite{KauffmanNeumann}, 
an $(l+\nu)$-Alexander matrix $P_{l+\nu}(K_*\otimes^\nu[2])$ is 
$\pm(t\cdot S_{l+\nu}(K_*)\pm ^t\hskip-1mm N_{l+\nu}(K_*))$. 
If we let $t=1$ or $t=-1$, then 
the determinant of $P_{l+\nu}(K_*\otimes^\nu[2])$ is not zero. 
Note that  $P_{l+\nu}(K_*\otimes^\nu[2])$ is a square matrix because 
$K_+$ is a spherical knot. 
Hence 
the linear map which is defined by $P_{l+\nu}(K_*\otimes^\nu[2])$ is injective. 

\vskip3mm
Hence 
the $(l+1+\nu)$-$\Q[t, t^{-1}]$-Alexander polynomial 
$A_{l+1+\nu}(K_*  \otimes^\nu[2])$ 
is the $\Q[t, t^{-1}]$-balanced class of 
the determinant of  
an $(l+1+\nu)$-Alexander matrix  
$P_{l+1+\nu}(K_*\otimes^\nu[2])$.

\vskip3mm
Proposition \ref{KLS} implies the following. 
$$S_{l+1+\nu}(K_* \otimes^\nu[2])=(-1)^\xi S_{l+1}(K_*).$$ 
$$N_{l+1+\nu}(K_*\otimes^\nu[2])
=(-1)^{\xi+\nu} N_{l+1}(K_*).$$  

\noindent ($*=+,-,0$.  \quad $\xi$ is a constant integer.)

\vskip3mm
Hence we can suppose that 

\noindent\hskip2cm
$P_{l+1+\nu}(K_*\otimes^\nu[2])$

\noindent\hskip2cm
$=t\cdot(-1)^\xi S_{l+1}(K_*)-(-1)^{\xi+\nu} N_{l+1}(K_*)$

\noindent\hskip2cm
$=(-1)^\xi (t\cdot S_{l+1}(K_*)+(-1)^{\nu+1} N_{l+1}(K_*))$ 

\noindent\hskip2cm
$=(-1)^\xi (t\cdot S_{l+1}(K_*)+(-1)^{\nu+1+l} \cdot  ^t\hskip-1mm S_{l+1}(K_*))$.  

\hskip5cm(Reason:    
$S_{\nu+1}(K_*)={(-1)^l}\hskip1mm ^t\hskip-1mm N_{\nu+1}(K_*)$ by \S\ref{Alex}.)

\vskip3mm
Hence we have the following. 
Let $P_{l+1+\nu}(K_*  \otimes^\nu[2])=(p^*_{ij})$.

\smallbreak
\noindent
(1) 
$(p^+_{ij})$ 
and 
$(p^-_{ij})$  
are $\rho\x \rho$ matrices.

\hskip3mm $(p^0_{ij})$  
is a $(\rho-1)\x (\rho-1)$ matrix. 
($\rho\geqq2, \quad \rho\in\N$).

\smallbreak
\noindent
(2) \hskip4cm
$p^+_{\rho,\rho}-p^-_{\rho,\rho}=c(t+(-1)^{l+1+\nu}),$

\hskip5mm where $c=\pm1$. 

\smallbreak
\noindent
(3)  
\hskip3cm$p^+_{ij}=p^-_{ij}=p^0_{ij}  (1\leqq i\leqq \rho-1, 1\leqq j\leqq \rho-1).$ 
$$p^+_{ij}=p^-_{ij}   ((i,j)\neq(\rho,\rho)).$$ 

\bigbreak

By calculus of determinants,   
$${\mathrm{det}}
P_{l+1+\nu}
(K_+ \otimes^\nu[2])
-
{\mathrm{det}}
P_{l+1+\nu}
(K_-\otimes^\nu[2])
= c (t+(-1)^{l+1+\nu})\cdot{\mathrm{det}}P_{l+1+\nu}(K_0\otimes^\nu[2]),$$

\noindent 
where $c=\pm1$.

Hence Theorem \ref{skeincro} holds. 
Note that the $(l+1+\nu)$-$\Q[t, t^{-1}]$-Alexander polynomial is 
a $\Q[t, t^{-1}]$-balanced class. 
\qed

\bigbreak
\noindent 
{\bf  Proof of Theorem \ref{proinertia}.}  
%
%
%
%
%
%
By \cite{KauffmanNeumann}, 
$\>K_*  
\otimes^\mu
(\text{the Hopf link})\>$ 
has a Seifert hypersurface   $V_*$ with a handle decomposition of 
one (4$\mu$+2)-dimensional 0-handle and 
 (4$\mu$+2)-dimensional (2$\mu$+1)-handles, 
where  there may not be (4$\mu$+2)-dimensional (2$\mu$+1)-handles ($*=+,-,0$).


By \cite{KauffmanNeumann} and Proposition \ref{KLS}, 
there is a Seifert matrix
$S_{2\mu+1} 
(\>K_* \otimes^\mu(\text{the Hopf link})\>)$ associated with $V_*$   
which is equal to  
a Seifert matrix $(-1)^\mu S_1(K_*)$ for the 1-knot $K_*$.   
(Recall 
$S_*(\quad)$,  
$N_*(\quad)$ 
in \S\ref{Alex}.)

Since 
$
\>K_*
\otimes^\mu
(\text{the Hopf link})\>
$
is a $(4\mu+1)$-submanifold,     
$N_{2\mu+1} 
(\>K_*
\otimes^\mu
(\text{the Hopf link})\>)$
$=^t\hskip-1mm S_{2\mu+1} 
(\>K_*
\otimes^\mu
(\text{the Hopf link})\>)$. 
Hence we can suppose that 
the absolute value of the determinant of a matrix 
which represents the intersection product \newline 
$H_{2\mu+1}(V_+;\Z)\x H_{2\mu+1}(V_+;\Z)\to\Z$ 
(resp. $H_{2\mu+1}(V_-;\Z)\x H_{2\mu+1}(V_-;\Z)\to\Z$) 
is $\pm1$. 
Hence  
$K_+
\otimes^\mu
(\text{the Hopf link})$ 
(resp. $K_-
\otimes^\mu
(\text{the Hopf link})$)  
is a homology sphere.

By the above handle decomposition of $V_*$ and \cite{Smale},    
$K_+
\otimes^\mu
(\text{the Hopf link})$ \newline
(resp. $K_-
\otimes^\mu
(\text{the Hopf link})$)  
is 
homeomorphic to the standard sphere.

By Proposition \ref{zero},  \newline 
$(K_+
\otimes^\mu
(\text{the Hopf link})$,   
$K_-
\otimes^\mu
(\text{the Hopf link})$,  
$K_0
\otimes^\mu
(\text{the Hopf link})$
)  
is a twist-move-triple.

By Theorem \ref{inertia}, the proof is completed. 
\qed

\bigbreak
\noindent 
{\bf  Proof of Theorem \ref{skeinpass}.}  
Proposition \ref{KLS} implies the following. 
$$S_{p+\nu} (\>K_* \otimes^\nu[2]\>)=(-1)^c S_p(K_*).$$ 
$$N_{p+\nu} 
(\>K_*\otimes^\nu[2]\>)=(-1)^{c+\nu} N_p(K_*),$$ 
\noindent  where $*=+,-,0$ and  $c$ is a constant integer.

\smallbreak
Since $K_+$, $K_-$ are spherical knots and 
$(K_+,K_-,K_0)$ is a $(p, q)$-pass-move-triple, 
we have the following.

There is a $(p+\nu)$-Alexander matrix

\noindent
$P_{p+\nu}(K_*\otimes^\nu[2])
=(-1)^c\{t\cdot S_p(K_*\otimes^\nu[2])-(-1)^{\nu}\cdot N_p(K_*\otimes^\nu[2])\}
=(p^*_{ij})$,   
which is a square matrix,  with the following property. 

\smallbreak
\noindent
(1) 
$(p^+_{ij})$ 
and 
$(p^-_{ij})$ 
are $\rho\x \rho$ matrices.

\hskip3mm$(p^0_{ij})$ 
is a $(\rho-1)\x (\rho-1)$ matrix. 
($\rho\geqq2, \quad \rho\in\N$)

\smallbreak\noindent
(2) 
\hskip4cm$p^+_{\rho,\rho}-p^-_{\rho,\rho}=(-1)^c\{t+(-1)^{1+\nu}\}.$

\smallbreak\noindent(3)  
\hskip3cm$p^+_{ij}=p^-_{ij}=p^0_{ij} (1\leqq i\leqq \rho-1,   1\leqq j\leqq \rho-1).$
$$p^+_{ij}=p^-_{ij}   ((i,j)\neq(\rho,\rho)).$$

\smallbreak

\noindent
Note we can suppose that $\rho-1\geqq1$ 
by using a surgery of Seifert hypersurface, if necessary.

\smallbreak

By calculus of determinants,

\smallbreak
\noindent 
${\mathrm{det}}P_{p+\nu}(K_+\otimes^\nu[2])
-{\mathrm{det}}P_{p+\nu}(K_-\otimes^\nu[2])  
=
(-1)^\zeta\cdot(t+(-1)^{\nu+1})\cdot{\mathrm{det}}P_{p+\nu}(K_0\otimes^\nu[2])
\cdot\cdot\cdot\cdot\cdot(!), $

\smallbreak
\noindent
where $\zeta$ is a constant integer. 

\smallbreak

Take ${\mathrm{det}}P_{p+\nu}(K_*\otimes^\nu[2])$ for each $*$ ($*=+,-,0$). 
Here, there are the following three cases (i)(ii)(iii).   
 
\smallbreak
\noindent 
(i) 
Suppose that  
${\mathrm{det}}P_{p+\nu}(K_*\otimes^\nu[2])
\neq0$ for a $*$. 
Let $p-1\neq n+1-p$. 
Then  $K_+$, $K_-$ and $K_0$ has 
a same $(p-1)$-Alexander matrix $t\cdot S_{p-1}-N_{p-1}$.  
The $(p-1)$-Alexander matrix has the following properties. 
Note that $t\cdot S_{p-1}-N_{p-1}$ is a square matrix. 
$t\cdot S_{p-1}-N_{p-1}$ is a nonsingular square matrix.  
Reason:   
$K_+$ and $K_-$ are spherical knots. Hence 
If $t=1$, $t\cdot S_{p-1}-N_{p-1}=S_{p-1}-N_{p-1}$ is nonsingular.

By \cite{KauffmanNeumann}, 
a $(p+\nu-1)$-Alexander matrix $P_{p+\nu-1}$
for
$K_*\otimes^\nu[2]$
is one of $\pm\{t\cdot S_{p-1}\pm N_{p-1}\}$.

If  $t=1$ or $t=-1$, 
$P_{p+\nu-1}$ is nonsingular. 
Hence 
$P_{p+\nu-1}$ is nonsingular. 
 
Hence the $({p+\mu})$-$\Q[t, t^{-1}]$-Alexander polynomial for 
$K_*\otimes^\nu[2]$
is 
${\mathrm{det}}P_{p+\mu}(K_*\otimes^\nu[2])$ for the $*$.

\smallbreak
\noindent 
(ii)
Suppose  that ${\mathrm{det}}P_{p+\nu}(K_*\otimes^\nu[2])\neq0$ for a $*$.  
Let $p-1=n+1-p$. 
A $(p+\nu-1)$-Alexander matrix for $K_*\otimes^\nu[2]$ defines an injective map 
because 
a $(p-1)$-Alexander matrix for $K_*$ does. 
See the identity right above Proposition \ref{Mars}. 


Hence the $({p+\nu})$-$\Q[t, t^{-1}]$-Alexander polynomial for $K_*$
is 
${\mathrm{det}}P_{p+\nu}(K_*\otimes^\nu[2])$ for the $*$.

\smallbreak
\noindent 
(iii) Suppose that ${\mathrm{det}}P_{p+\nu}(K_*\otimes^\nu[2])=0$ for a $*$. 
Then the $(p+\nu)$-$\Q[t, t^{-1}]$-Alexander polynomial for 
$K_*\otimes^\nu[2]$ is 
$0={\mathrm{det}}P_{p+\nu}(K_*\otimes^\nu[2])$ for the $*$.     
Here, note that we do not need to consider whether 
the linear map defined by a $(p+\nu-1)$-Alexander matrix is injective or not. 

\smallbreak

By the above (i), (ii), (iii) and the identity (!) several lines above here, 
the proof is completed.  
Note that the $(p+\nu)$-$\Q[t, t^{-1}]$-Alexander polynomial is 
a $\Q[t, t^{-1}]$-balanced class. 
\qed

\section{A problem}\label{open}   

\smallbreak
\begin{prob}\label{difficult}   
Suppose that $n$-dimensional closed oriented submanifolds 
$K$ and $K'$ differ by one twist-move (resp. pass-move). 
Suppose that $m$-dimensional closed oriented submanifolds 
$J$ and $J'$ differ by one twist-move (resp. pass-move). 
Then how do we characterize a relation between 
$K\otimes J$  
and  
$K'\otimes J'$?  

\noindent
(Recall the following. 
Let $n=1$. 
$K$ and $K'$ differ by one twist-move if and only if 
$K$ and $K'$ differ by one crossing-change.)
\end{prob}



\bigbreak\noindent
Louis H. Kauffman: 
Department of Mathematics, Statistics, and Computer Science,  
University of Illinois at Chicago, 
851 South Morgan Street, 
Chicago, Illinois 60607-7045, USA  \quad
kauffman@uic.edu

\bigbreak\noindent
Eiji Ogasa:  Computer Science, Meijigakuin University, Yokohama, Kanagawa, 244-8539, Japan 
\quad pqr100pqr100@yahoo.co.jp  \quad
ogasa@mail1.meijigkakuin.ac.jp

\end{document}